\documentclass[preprint,12pt]{elsarticle}

\usepackage{amssymb}
\usepackage{amsmath}
\usepackage{longtable, makecell, multirow, multicol}
\usepackage{booktabs}
\usepackage{lscape}
\usepackage{array}
\usepackage{float}
\usepackage{adjustbox}
\usepackage{multirow}
\newcolumntype{C}[1]{>{\centering\arraybackslash}p{#1}}
\usepackage[verbose]{placeins}

\usepackage{xcolor}

\usepackage{graphicx}

\usepackage{url}
\usepackage[breaklinks]{hyperref}
\usepackage{hyperref}
\usepackage{breakurl}

\usepackage[symbol*]{footmisc}

\title{Optimal Experimental Design using Eigenvalue-Based Criteria with Pyomo.DoE}
\author{Daniel J. Laky${}^{1,2}$}
\author{Shammah Lilonfe${}^{1}$}
\author{Shawn B. Martin${}^{3}$}
\author{Katherine A. Klise${}^{4}$}
\author{Bethany L. Nicholson${}^{5}$}
\author{John D. Siirola${}^{5}$}
\author{Alexander W. Dowling\corref{adowling}$^{1}$}
\address{$^{1}$Department of Chemical and Biomolecular Engineering \\ University of Notre Dame, Notre Dame, IN 46556}
\address{$^{2}$Departmant of Chemical Engineering \\ Auburn University, Auburn, AL 36849}
\address{$^{3}$ Mission Analytics, Sandia National Laboratories \\ Albuquerque, NM 87185}
\address{$^{4}$ Energy Water Systems Integration, Sandia National Laboratories \\ Albuquerque, NM 87185}
\address{$^{5}$ Center for Computing Research, Sandia National Laboratories \\ Albuquerque, NM 87185}
\journal{arXiv}

\begin{document}

\begin{frontmatter}

\begin{abstract}
Leveraging digital twins to accelerate scientific discovery requires acquisition of high-quality data to ensure predictive power. Time and resource limitations motivate the deployment of model-based design of experiments to elucidate optimal experimental campaigns to build and refine digital twins that realize value while respecting resource budgets. Additionally, control and optimization tasks, which can be enhanced by using equation-oriented optimization with algebraic models, enable value-adding decision making with predictive digital twins. Pyomo.DoE is a software package for optimal experimental design to build high-fidelity, equation-oriented models. Oftentimes, these high-fidelity digital models suffer from numerical errors due to identifiability issues and poor model scaling. Optimal experimental design helps to address these issues with specific information-based optimal design metrics, such as minimum eigenvalue optimality (E-optimality) and condition number optimality (ME-optimality), combating these problems directly by focusing on the numerically problematic portions of the model. However, embedding the sophisticated linear algebra functions (e.g., matrix inversion, eigenvalue computation) required during optimal experimental design remains a challenge, especially in equation-oriented optimization frameworks that leverage state-of-the-art derivative-based optimization tools. This work extends Pyomo.DoE to include callbacks that allow rigorous computation of eigenvalue-based experimental design metrics, resulting in heightened focus on parameters that are difficult to identify in the model, especially using equation-oriented programming. In addition, a brief tutorial on experimental design metrics is given in the methodology and supplementary information. Finally, we propose a new experiment-creation modeling abstraction for intrusive uncertainty quantification in Pyomo, demonstrating that aligning model-to-software abstractions reduces user modeling time by harmonizing critical steps in the workflow for building and refining high-value digital twins. The work highlights that choosing a design metric, or metrics, that best aligns with the experimental objective is paramount to gaining desired information.
\end{abstract}


\end{frontmatter}

\section{Introduction} \label{sec:intro}
Digital twins are ``{\em... sets of virtual information constructs that mimic the structure, context, and behavior of a natural, engineered, or social system (or system-of-systems), are dynamically updated with data from their physical twins, have predictive capability, and inform decisions that realize value...}'' \cite{nationalacademy2024digitaltwin, aiaa2020digitaltwin}. Developing and maintaining the predictive capability of digital twins requires high-quality data \cite{macias2024digital-twin-data-quality, yan2025digital-twin-data-reqs} to perpetuate iterative refinement between digital and physical systems. Especially during exploratory analysis and the earliest iterations of digital twin development, gathering the most information from limited resources (e.g., time, personnel, material availability, computational capacity) is critical.

Since digital twins have a mathematical structure (or collection of candidate mathematical structures), exploiting model-based design of experiments (MBDoE) to design objectively meaningful experimental campaigns is paramount to efficiently utilize limited resources. Maximizing or minimizing an objective function related to the informational content of a potential experiment is typically referred to as optimal experimental design. Many MBDoE frameworks compute optimal experimental campaigns by maximizing the information content of the next best experiment(s) while directly considering the underlying model structure utilizing the Fisher Information Matrix (FIM) \cite{fisher1971design}. As such, MBDoE (among other design of experiments strategies) has facilitated the development of automated and self-driving laboratories in many materials discovery~\cite{abolhasani2023automation, tom2024automationreview, seifrid2022semiconductorautomation} and reaction kinetics determination~\cite{quaglio2019automation, waldron2020automation, agunyole2024automation} applications. MBDoE has also been important to enable automatic model detection (e.g., symbolic regression \cite{gomes2019doe-using-symbolic, can2011doe-symbolic, castillo2003doe-symbolic, rogers2024doe-symbolic}). Purely data-driven approaches, such as active learning (in applications such as catalysis \cite{saidi2017doe-ann}, injection molding \cite{heinisch2021doe-ann}, and chemical reactions \cite{sangoi2022doe-ann}) and Bayesian optimization (in applications such as biological networks \cite{imani2020bayesian-doe, imani2022bayesian-doe}, dynamic chemical and control processes \cite{cao2024bayesian-doe, tulsyan2012bayesian-doe}, material design \cite{lei2021bayesian-doe}, fluid flow \cite{attia2022bayesian-doe}, among others \cite{greenhill2020bayesian-doe-review}), have also benefited from and contributed to optimal experimental design. 

However, data-driven approaches typically do not utilize first-principles models to directly inform experimental design. Therefore, these data-driven methods often forgo important structural and physical insight from these science-based models. To this end, advancements in MBDoE for parameter precision \cite{franceschini2008MBDoE, galvanin2010MBDoE-uncertainty} and model discrimination tasks \cite{galvanin2016MBDoE-for-MD-and-PE} that directly incorporate first-principles models have reshaped the landscape for optimal experimental design. See the seminal review paper by Franceschini and Macchietto~\cite{franceschini2008MBDoE} and subsequent recent reviews \cite{geremia2026MBDoE-review} for more information.

The popularity of MBDoE has led to the development of many software tools with varying levels of generality in \texttt{Python} (PyOED \cite{chowdhary2024pyOED}, PyDOE \cite{pyDOE_documentation}, MIDDoE~\cite{tabrizi2025middoe-conference-paper}), and \texttt{R} (POPED \cite{foraccia2004POPED}, \texttt{odw} \cite{butler2022odw}). However, some applications require detailed models that result in large-sacle optimization problems which are only tractable with intrusive (e.g., equation-oriented) tools~\cite{agi2024computation}. One commercial tool capable of intrusive MBDoE is gPROMS \cite{gproms2025}, which utilizes equation-oriented large-scale nonlinear dynamic optimization to interrogate digital models. Pyomo \cite{bynum2021pyomo}, an open-source python package for optimization, recently released a contributed package for MBDoE (Pyomo.DoE \cite{wang2022pyomo-doe}) which also leverages nonlinear programming. gPROMS utilizes control vector parameterization for dynamic optimization; Pyomo discretizes a dynamic system using collocation or finite difference methods to explicitly represent the complete time horizon of the model as algebraic equations \cite{nicholson2018pyomo-dae}.

In the case of simultaneous equation-oriented programming (e.g. using Pyomo), these tools are usually limited to D-optimal designs, given the relative ease of implementation of the log-determinant of the FIM in an equation-oriented framework, or relaxations of experimental design criteria (e.g., E-optimality \cite{boyd2004convex}). D-optimal designs tend to focus on the largest eigenvalue of the FIM (the direction of highest potential information), leading to designs that may ignore improvements to parameters with low information \cite{montgomery2017design} (conversely, high uncertainty). Therefore, designs that mathematically focus on directions of lowest information (or highest uncertainty), such as A-optimal and E-optimal designs, may be desired. In equation-oriented frameworks, A-optimal designs require the information matrix to be inverted using fully defined algebraic constraints, which has only recently been shown \cite{duarte2022A-opt} and may suffer from numerical instabilities when the information matrix is (nearly) singular. 

In addition to focusing on the areas of the model with low information, balancing the ratio between the highest and lowest directions of information (e.g., ME-optimality) may ensure a more balanced magnitude of information, and subsequently uncertainty, among all unknown parameters. Metrics that require eigenvalues (E-optimality and ME-optimality) are more difficult (or impossible) to pose as algebraic equations, with no explicit equations available for the eigenvalues of a general square matrix (e.g., the FIM) on problems that have more than four unknown parameters \cite{ruffini1813degree5, abel1826degree5}. There is a reformulation for E-optimal designs originally shown in Boyd \cite{boyd2004convex} and used in some research applications \cite{telen2015E-opt, ye2016Eopt}. Although useful, the reformulation does not guarantee mathematically E-optimal designs for general nonlinear problems. To the authors' knowledge, no general method for determining E- and ME-optimal experimental designs in an equation-oriented framework exists.

To address these gaps and related issues, we present new capabilities in Pyomo.DoE that improve equation-oriented optimal experimental design and standardize the model development workflow, thus unifying the critical steps to build and refine a digital twin. We leverage callbacks in Pyomo \cite{rodriguez2020pynumero, laky2022callbacks, wang2024measurement} to directly compute statistical measures of information content (A-optimality, D-optimality, E-optimality, and ME-optimality) and their derivatives to facilitate large-scale optimization. Furthermore, a unifying model abstraction is shown to streamline the critical steps in building and refining a digital twin, significantly reducing the modeling burden and enabling a higher degree of standardization in digital workflows. This enables significant improvement over the existing workflow for model development in Pyomo.DoE \cite{wang2022pyomo-doe}, and adds three previously inaccessible optimality criteria (i.e., E-optimality, ME-optimality, A-optimality) to the existing optimal experimental design workflow.

The rest of this work is organized as follows. Section~\ref{sec:methodology} reviews the underlying mathematics for MBDoE and is presented with particular focus on how callbacks simplify the model (with more detailed math included in the supplementary information). Section~\ref{sec:experiment-abstraction} gives an overview of the new experiment modeling abstraction in Pyomo.DoE that unifies the workflow to build digital twins. Finally, Section~\ref{sec:results} demonstrates these capabilities using a series of case studies including a large-scale model for the development of critical minerals separations (Section~\ref{sec:res-membrane}) and concludes with a discussion and future opportunities for further improvement.

\section{Methodology} \label{sec:methodology}
\begin{figure}[H]
\includegraphics[width = \textwidth]{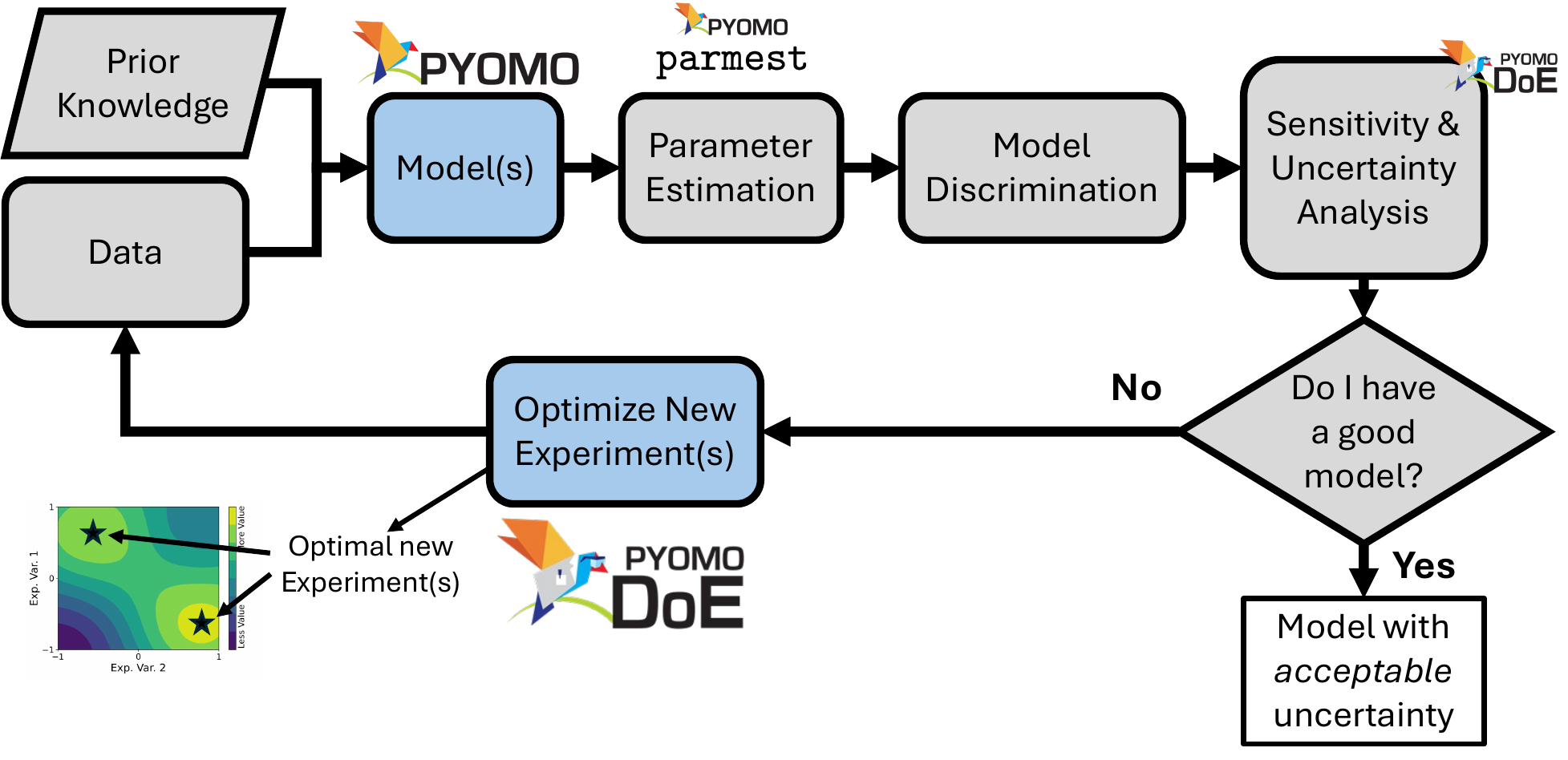}
\caption{General workflow for model building and identification using Pyomo, \texttt{parmest}, and Pyomo.DoE. This work focuses on equation-oriented modeling in Pyomo and optimally designing experiments for parameter precision (boxes highlighted in blue).} \label{fig:workflow}
\end{figure}

A predictive digital twin realizes the greatest value when closely emulating its physical counterpart. Building such digital twins requires iterative model building and refinement between the physical twin behavior (e.g., experimental data) and the digital twin (e.g., model prediction). Broadly, the workflow to build and refine digital twins is shown in Figure~\ref{fig:workflow} based on similar workflows in literature \cite{franceschini2008MBDoE, galvanin2010MBDoE-uncertainty, galvanin2016MBDoE-for-MD-and-PE, wang2022pyomo-doe, lynch2024crystallization-mbdoe, geremia2026MBDoE-review}. We consider the case where one or more candidate models are available to digitally describe the behavior of the physical system. Each model contains uncertain parameters that must be inferred to uniquely describe the physical system. Preliminary data is used to generate a best fit estimate of these parameters for the model(s). Then, the model(s) is (are) analyzed to understand sensitivity and uncertainty associated with the current estimate of unknown parameters. Often, the uncertainty in the model(s) is too high (low confidence in predictive power), leading the model builder to request more data to improve certainty with the model(s). The key question is the following: \textit{What experiment(s) should we perform next to improve the quality of the digital twin?}

\FloatBarrier

As mentioned in Section~\ref{sec:intro}, MBDoE facilitates designing optimal experiments in the context of the Fisher Information Matrix (FIM) \cite{fisher1971design}. This paper focuses on a special case of MBDoE where first-principles models are used to generate the FIM for nonlinear models \cite{franceschini2008MBDoE}. More specifically, we are interested in posing information metrics (linear algebra measures) of the FIM as algebra alongside these first-principles models, allowing the use of powerful, derivative-based optimization solvers (e.g., \texttt{IPOPT} \cite{ipopt}, \texttt{CONOPT} \cite{conopt}, \texttt{KNITRO} \cite{knitro}, \texttt{BARON} \cite{baron}).

In the following subsections, a brief overview of the FIM is provided, where readers are referred to several review/tutorial papers \cite{franceschini2008MBDoE, wang2022pyomo-doe, geremia2026MBDoE-review} and references within for more detailed information. Then, a tutorial-style description of the FIM metrics and how to compute them is given. We particularly focus on novelties that embed eigenvalue-based metrics in equation-oriented programs using callbacks.

\subsection{Parameter Estimation and the Fisher Information Matrix} \label{sec:parmest-and-fim}
In statistical inference, the model, $\boldsymbol{f}$, at a specific experimental condition, $\boldsymbol{\phi}_i$, is related to its physical response (experimental observation), $\boldsymbol{y}_i$, as shown below:

\allowdisplaybreaks
\begin{align}
    &\boldsymbol{y}_i = \boldsymbol{f}\left(\boldsymbol{\phi}_i, \boldsymbol{x}_i, \boldsymbol{\theta}\right)\;+ \boldsymbol{\varepsilon}_i, & \boldsymbol{\varepsilon}_i \, {\sim} \, \mathcal{N}\left(\boldsymbol{0}, \boldsymbol{\Sigma}_{\boldsymbol{y}}\right), \; \forall \; i \in \left\{1, \ldots, N_\text{exp}\right\} \label{eq:general-model}
\end{align}
\allowdisplaybreaks[0]

\noindent where $\boldsymbol{x}_i$ constitute any state or algebraic variables that will change with the choice of $\boldsymbol{\phi}_i$ and model structure of $\boldsymbol{f}$. We also assume that the data are corrupted by some observation error, $\boldsymbol{\varepsilon}_i$, which is correlated with a normal distribution having mean 0 and covariance matrix, $\boldsymbol{\Sigma}_{\boldsymbol{y}}$, known \emph{a priori} for any given $\boldsymbol{y}$ from $\left\{1, \ldots, N_\text{exp}\right\}$, defined below:

\begin{align}
    \boldsymbol{\Sigma}_{\boldsymbol{y}} = 
        \begin{bmatrix}
            \sigma^2_{y_1,y_1} & 0  & \ldots & 0   \\[0.35cm]
            0 & \sigma^2_{y_2,y_2}  & \ldots & 0   \\[0.1cm]
            \vdots & \vdots & \ddots & \vdots  \\[0.1cm]
            0 & 0  & \ldots & \sigma^2_{y_{N_\text{meas}},y_{N_\text{meas}}}   \\
        \end{bmatrix} \label{eq:measurement-error-matrix-def}
\end{align}
\noindent where $\sigma_{\hat{y}_i, \hat{y}_j}$ is the pairwise variance between measurement $i$ and measurement $j$. Under certain conditions (e.g., assuming independent measurements) this covariance matrix is reduced to a diagonal matrix (see \cite{wang2022pyomo-doe, franceschini2008MBDoE}), as shown above. However, these equations are valid for any general multivariate Gaussian distribution when the elements indicated as zero in Eq.~\ref{eq:measurement-error-matrix-def} are replaced by $\sigma_{\hat{y}_i, \hat{y}_j}^2$.

These data, $\boldsymbol{y}_i$, and experimental condition sets, $\boldsymbol{\phi}_i$, are indexed by experiment, $i=1,\dots,N_\text{exp}$. Each $\boldsymbol{y}_i$ is a vector of measurements for experiment $i$ with length $N_{\text{meas}}$. The goal of model-building and refinement is to determine the values of the unknown parameters, $\boldsymbol{\theta}$, that best fit the model. To obtain best-fit estimates, it is common to construct a likelihood function and find the parameter estimates that maximize the likelihood (see Section~\ref{sec:si-parmest} for details). The log likelihood function, $\ell$, is shown below:

\allowdisplaybreaks
\begin{align}
    &\ell\left(\boldsymbol{\theta}, \boldsymbol{y}_1, \ldots, \boldsymbol{y}_{N_\text{exp}}\right) = \frac{-N_\text{exp}}{2}\left(\text{ln}\left(2\pi\right) - \ln\left|\text{det}\left(\boldsymbol{\Sigma}_{\boldsymbol{y}}^{-1}\right)\right|\right) \notag\\
    &\qquad\qquad-\frac{1}{2}\sum_{i \in \left\{1, \ldots, N_\text{exp}\right\}}\left(\boldsymbol{y}_i - \boldsymbol{f}\left(\boldsymbol{\phi}_i, \boldsymbol{x}_i, \boldsymbol{\theta}\right)\right)^{\text{T}} \boldsymbol{\Sigma}_{\boldsymbol{y}}^{-1} \left(\boldsymbol{y}_i - \boldsymbol{f}\left(\boldsymbol{\phi}_i, \boldsymbol{x}_i, \boldsymbol{\theta}\right)\right) \label{eq:log-likelihood}
\end{align}
\allowdisplaybreaks[0]

\noindent where $N_\text{exp}$ is the total number of observations, $\boldsymbol{y}$, and $\boldsymbol{\Sigma}_{\boldsymbol{y}}$ is the covariance of the measurements made during experiment $i$, from Eq.~\ref{eq:measurement-error-matrix-def}.

When $\boldsymbol{\Sigma}_y$ is known \emph{a priori}, it directly follows that the maximization of this likelihood function is the minimization of a weighted sum of squared errors function (see Section~\ref{sec:si-parmest} for details):

\allowdisplaybreaks
\begin{align}
    \hat{\boldsymbol{\theta}} =\text{arg}\min_{\boldsymbol{\theta} \in \boldsymbol{\Theta}} & \sum_{i \in \left\{1, \ldots, N_\text{exp}\right\}}\left(\boldsymbol{y}_i - \boldsymbol{f}\left(\boldsymbol{\phi}_i, \boldsymbol{x}_i, \boldsymbol{\theta}\right)\right)^{\text{T}} \boldsymbol{\Sigma}_y^{-1} \left(\boldsymbol{y}_i - \boldsymbol{f}\left(\boldsymbol{\phi}_i, \boldsymbol{x}_i, \boldsymbol{\theta}\right)\right) \label{eq:param-estimation} \\
    \text{s.t.} & \;\;\; \hat{\boldsymbol{y}}_i = \boldsymbol{f}\left(\boldsymbol{\phi}_i, \boldsymbol{x}_i, \boldsymbol{\theta}\right) \qquad \qquad \qquad \forall \; i \in \left\{1, \ldots, N_\text{exp}\right\} \label{eq:y-hat}
\end{align}
\allowdisplaybreaks[0]

\noindent where the goal is to find the unknown parameters, $\hat{\boldsymbol{\theta}}$, which minimize the weighted sum of the squared difference (maximize the likelihood of fit) between the measurement, $\boldsymbol{y}_i$, and the model prediction, $\hat{\boldsymbol{y}}_i$ given by Eq.~\ref{eq:y-hat}, summed over all $N_\text{exp}$ experiments. The feasible bounds for uncertain parameters $\boldsymbol{\theta}$ are represented using $\boldsymbol{\Theta}$.

A natural next step in the model building process is to understand the sensitivity of the model predictions, $\hat{\boldsymbol{y}}$, to changes in the optimally fitted parameters, $\hat{\boldsymbol{\theta}}$. This is also often referred to as the sensitivity matrix:

\begin{align}
    \boldsymbol{Q} = 
        \begin{bmatrix}
            \frac{\partial\hat{y}_1}{\partial\theta_1} & \frac{\partial\hat{y}_1}{\partial\theta_2} & \ldots & \frac{\partial\hat{y}_1}{\partial\theta_p}  \\[0.35cm]
            \frac{\partial\hat{y}_2}{\partial\theta_1} & \frac{\partial\hat{y}_2}{\partial\theta_2} & \ldots & \frac{\partial\hat{y}_2}{\partial\theta_p}  \\[0.1cm]
            \vdots & \vdots & \vdots & \vdots  \\[0.1cm]
            \frac{\partial\hat{y}_{N_\text{meas}}}{\partial\theta_1} & \frac{\partial\hat{y}_{N_\text{meas}}}{\partial\theta_2} & \ldots & \frac{\partial\hat{y}_{N_\text{meas}}}{\partial\theta_p}  \\
        \end{bmatrix}
\end{align}
\noindent where the dimensions of $\boldsymbol{Q}$ are the number of observations for the experiment, $N_{\text{meas}}$, by the number of unknown parameters, $p$. The method Pyomo.DoE uses to gather this information is to write finite differencing equations on $\hat{\boldsymbol{y}}$ with respect to the parameters $\boldsymbol{\theta}$ at the current estimate $\hat{\boldsymbol{\theta}}$. An example with the central difference equation is shown below:

\begin{align}
    \boldsymbol{q}_{j}\left(\boldsymbol{\phi}, \boldsymbol{\theta}\right) &= \frac{\hat{\boldsymbol{y}}\left(\boldsymbol{\phi}, \hat{\boldsymbol{\theta}} \;+ h_j\boldsymbol{e}_j\right) - \hat{\boldsymbol{y}}\left(\boldsymbol{\phi},\hat{\boldsymbol{\theta}} \;- h_j\boldsymbol{e}_j\right)}{2h_j} &\forall \; j \in \left\{1, \ldots, p\right\} \label{eq:fim-fd} \\
    \boldsymbol{Q}\left(\boldsymbol{\phi}, \hat{\boldsymbol{\theta}}\right) &= 
        \begin{bmatrix}
            q_{1,1} & q_{1,2}  & \ldots & q_{1,p}   \\[0.35cm]
            q_{2,1} & q_{2,2}  & \ldots & q_{2,p}   \\[0.1cm]
            \vdots & \vdots & \vdots & \vdots  \\[0.1cm]
            q_{N_{\text{meas}},1} & q_{N_{\text{meas}},2}  & \ldots & q_{N_{\text{meas}},p}  \\
        \end{bmatrix} \label{eq:sens-mat}
\end{align}
\noindent where $h_j$ and $\boldsymbol{e}_j$ are the perturbation step size of parameter $\theta_j$ and unit direction vector for element $j$, respectively. Here, $\boldsymbol{q}_{j}$ represents the finite difference derivative for model prediction vector, $\hat{\boldsymbol{y}}$, with respect to parameter $\theta_j$ at the current best fit estimate $\hat{\boldsymbol{\theta}}$, and an arbitrary experimental design, $\boldsymbol{\phi}$. Additionally, we define $\boldsymbol{q}_{j}$ as a column vector with dimensions $N_{\text{Meas}}$ by 1. This means the elements of $\boldsymbol{Q}$, $q_{i,j}$ are the finite difference derivative of the $i$-th measurement with respect to the $j$-th parameter.

As shown in previous literature \cite{franceschini2008MBDoE, wang2022pyomo-doe}, this sensitivity matrix, $\boldsymbol{Q}\left(\boldsymbol{\phi}, \hat{\boldsymbol{\theta}}\right)$, can be used to approximate the covariance matrix of the unknown parameters at current estimate, $\hat{\boldsymbol{\theta}}$, using the following relationship:

\allowdisplaybreaks
\begin{align}
    \boldsymbol{V}\left(\boldsymbol{\phi}, \hat{\boldsymbol{\theta}}\right) & \approx \left[\boldsymbol{Q}^T\left(\boldsymbol{\phi}, \hat{\boldsymbol{\theta}}\right)\boldsymbol{\Sigma}_{y}^{-1}\boldsymbol{Q}\left(\boldsymbol{\phi}, \hat{\boldsymbol{\theta}}\right) + \boldsymbol{V}_{\boldsymbol{\theta}}\left(\hat{\boldsymbol{\theta}}\right)^{-1}\right]^{-1} \label{eq:cov-def}
\end{align}
\allowdisplaybreaks[0]

\noindent where $\boldsymbol{V}$ is the parameter covariance matrix from experiment conditions $\boldsymbol{\phi}$ with prior $\boldsymbol{V}_{\boldsymbol{\theta}}$ from the experiments used to find best estimate, $\hat{\boldsymbol{\theta}}$. The sensitivity matrix, $\boldsymbol{Q}$, comes from Eq.~\ref{eq:sens-mat} and measurement covariance matrix, $\boldsymbol{\Sigma}_{\boldsymbol{y}}$, comes from Eq.~\ref{eq:measurement-error-matrix-def}.

Finally, MBDoE typically utilizes the FIM \cite{fisher1971design}, which measures the amount of information about the best estimate, $\hat{\boldsymbol{\theta}}$, for given model responses $\hat{\boldsymbol{y}}$ at experimental conditions $\boldsymbol{\phi}$. The FIM represents the curvature (second derivative) of the log-likelihood function with respect to unknown parameters, $\boldsymbol{\theta}$ \cite{bard1974nonlinear}. Fortunately, for a consistent and asymptotically normal estimator, $\hat{\boldsymbol{\theta}}$, obtained by maximizing the likelihood function, the FIM is related to the covariance matrix approximately via a matrix inverse, as shown below:

\allowdisplaybreaks
\begin{align}
    \boldsymbol{V}\left(\boldsymbol{x}, \hat{\boldsymbol{\theta}}\right) & \approx \left[\boldsymbol{M}\left(\boldsymbol{x}, \hat{\boldsymbol{\theta}}\right)\right]^{-1} \label{eq:cov-fim-relationship}
\end{align}
\allowdisplaybreaks[0]

\noindent where $\boldsymbol{M}$ is the FIM. This result is from the seminal work of Rao in 1945 \cite{rao1945rao-cramer-bound}, who notes that the covariance matrix is bounded below by the inverse of the FIM, sometimes referred to as the ``Cram$\acute{\text{e}}$r-Rao'' bound. Those interested in learning more about this assumption should read the original work of Rao \cite{rao1945rao-cramer-bound}, or see a modernized description from Nielsen \cite{nielsen2013rao-cramer-chapter}.

Scalarization of the FIM is required to optimize the information in a potential experiment. The following section gives a brief description of the optimal design criteria available.

\subsection{Information-based Design Criteria} \label{sec:design-criteria}
Using these formulas, we can pose the optimal experimental design problem, as shown below:

\allowdisplaybreaks
\begin{subequations} \label{eq:formulation-1}
\begin{align}
    \min_{\boldsymbol{\phi}\in \boldsymbol{\Phi}} & \;\; \Psi\left[\boldsymbol{M}\left(\boldsymbol{\phi}, \hat{\boldsymbol{\theta}}\right)\right] \\
    \text{s.t.} & \;\; \boldsymbol{M}= \boldsymbol{Q}^T\boldsymbol{\Sigma}_{\boldsymbol{y}}^{-1}\boldsymbol{Q} + \boldsymbol{M}_{\hat{\boldsymbol{\theta}}} \\
    &\boldsymbol{q}_{j}\left(\boldsymbol{\phi}, \boldsymbol{\theta}\right) = \frac{\hat{\boldsymbol{y}}\left(\boldsymbol{\phi}, \hat{\boldsymbol{\theta}} \;+ h_j\boldsymbol{e}_j\right) - \hat{\boldsymbol{y}}\left(\boldsymbol{\phi},\hat{\boldsymbol{\theta}} \;- h_j\boldsymbol{e}_j\right)}{2h_j} \\
    & \qquad\qquad\qquad\qquad\qquad\qquad\qquad\qquad\forall \; j \in \left\{1, \ldots, p\right\} \notag \\
    &\boldsymbol{Q}\left(\boldsymbol{\phi}, \hat{\boldsymbol{\theta}}\right) = 
        \begin{bmatrix}
            q_{1,1} & q_{1,2}  & \ldots & q_{1,p}   \\[0.35cm]
            q_{2,1} & q_{2,2}  & \ldots & q_{2,p}   \\[0.1cm]
            \vdots & \vdots & \vdots & \vdots  \\[0.1cm]
            q_{N_{\text{meas}},1} & q_{N_{\text{meas}},2}  & \ldots & q_{N_{\text{meas}},p}  \\
        \end{bmatrix}
\end{align}
\end{subequations}
\allowdisplaybreaks[0]

\noindent where $\boldsymbol{M}_{\hat{\boldsymbol{\theta}}}$ is the prior information matrix from previous experiments used to optimize the best fit estimate, $\hat{\boldsymbol{\theta}}$, and $\boldsymbol{\Phi}$ is the experimental design space. In these problems, $\Psi$ represents a scalarized objective function of the FIM, $\boldsymbol{M}$. The goal is to find experimental conditions, $\hat{\boldsymbol{\phi}}$, that when executed will add the most information from the corresponding measurements. Oftentimes, these objective functions have to do with the eigenvalues (and eigenvectors) of the information matrix. Eigenvalues indicate the magnitude of information and eigenvectors indicate the direction of that information with respect to the unknown parameters. This means that a large eigenvalue of the FIM indicates a large amount of information, which in turn, because the covariance matrix is related via its inverse, corresponds to low uncertainty in that direction. There are many different scalarized objective functions for FIM-based optimal experimental design; however, those that are relevant to this paper are listed below:

\begin{itemize}
\item \textbf{A-optimality} minimizes the trace of the covariance matrix, $\boldsymbol{V}$. This is equivalent to minimizing the trace of the inverse of the FIM, as shown below:

\begin{align}
    \min\text{trace}\left(\boldsymbol{M}^{-1}\right) = \sum_{i}^p \frac{1}{\lambda_i} \label{eq:obj-trace}
\end{align}
where $\lambda_i$ is the $i$-th eigenvalue of the FIM. As shown from the formula, increasing a large eigenvalue will have less impact than increasing a small eigenvalue by the same amount (because the metric is related to the \textit{inverse} of the FIM). This means that A-optimal designs will tend to focus on improving the smallest eigenvalues, or the directions of lowest information.

\item \textbf{D-optimality} minimizes the determinant of the covariance matrix, $\boldsymbol{V}$. In turn, this is equivalent to maximizing the determinant of the FIM, shown in the equation below:

\begin{align}
    \max\left|\boldsymbol{M}\right| = \prod_{i}^p \lambda_i \label{eq:obj-log-det}
\end{align}
Here, since there is a product of eigenvalues, increasing the eigenvalue by the highest factor dominates the metric. Typically, these are the directions that are already of high information, sometimes biasing D-optimal designs to the directions of highest certainty. In the same vein, when minimizing the determinant of $\boldsymbol{V}$, shrinking one dimension of uncertainty may drastically reduce the determinant overall but retain high uncertainty in other dimensions.

\item \textbf{E-optimality} - minimizes the maximum eigenvalue of the covariance matrix, $\boldsymbol{V}$. This is the same as maximizing the minimum eigenvalue of the FIM, as shown below:

\begin{align}
    \max \;\min_{i\in\left\{1, \ldots, p\right\}}\lambda_i \label{eq:obj-min-eig}
\end{align}
In this case, the minimum eigenvalue is directly operated on, meaning we are finding the experiment that improves the direction of lowest information (highest uncertainty). This metric theoretically works complimentarily with metrics like D-optimality.

\item \textbf{ME-optimality} minimizes the condition number of the covariance matrix, $\boldsymbol{V}$. This is equivalent to minimizing the condition number of the FIM, shown in the equation below:

\begin{align}
    \min \; \frac{\max_{i\in\left\{1, \ldots, p\right\}}\lambda_i}{\min_{i\in\left\{1, \ldots, p\right\}}\lambda_i} \equiv \min \kappa \label{eq:obj-cond}
\end{align}
It is important to note that Eq.~\ref{eq:obj-cond} is valid to compute the condition number, $\kappa$, because the FIM (and covariance matrix $\boldsymbol{V}$) are \textit{normal} matrices as they are real-symmetric matrices. The goal here is to ensure that the relative magnitude of certainty in unknown parameters, $\boldsymbol{\theta}$, is as balanced as possible. Geometrically, the covariance matrix, or the FIM, can be represented as a hyperellipsoid using the formula below:

\begin{align}
    &\boldsymbol{a}^T\boldsymbol{V}\boldsymbol{a} \\
    &\boldsymbol{a} \in \mathbb{R}^p
\end{align}
where $\boldsymbol{a}$ is a vector in $\mathbb{R}^p$. The interpretation here is that the fundamental shape when minimizing the condition number is that $\boldsymbol{V}$, or $\boldsymbol{M}$, will be pushed closer to that of a ball (sphere in 3 dimensions, or balanced certainty in all parameter directions) than an elongated ellipsoid (more of an egg-like, or rugby-ball shape in 3 dimensions; meaning unbalanced certainty in some parameter directions). Focusing on a non-minor direction of the FIM (e.g., D-optimality) may elongate the ellipsoid, effectively increasing the disparity of certainty between the parameters even though information is ``increased'' overall. It should be noted that ME-optimality is agnostic to the magnitude of specific eigenvalues, as only the ratio is desired. Therefore, it is recommended to use ME-optimality to supplement another design metric that is improving information content by increasing eigenvalue magnitude. There is brief mention of this idea in Case Study 1.

\item \textbf{Pseudo-A-optimality} maximizes the trace of the FIM:

\begin{align}
    \max\text{trace}\left(\boldsymbol{M}\right) = \sum_{i}^p \lambda_i \label{eq:obj-pseudo-trace}
\end{align}

This objective has been mistakenly identified (and continues to be identified) as equivalent to A-optimality~\cite{franceschini2008MBDoE, wang2022pyomo-doe}. Shown later in Case Study 1, the trace of the FIM does not emulate the behavior of A-optimality. However, the trace of the FIM does not require inversion or eigenvalues, making it numerically attractive in near-singular cases. Therefore, in some instances, pseudo-A-optimality provides some value when recommending an experiment but should not be confused with A-optimality.

\end{itemize}

Throughout this work, we assume that $\boldsymbol{M}$ is full rank, i.e., the minimum eigenvalue is greater than 0. If the minimum eigenvalue is zero, which indicates the model is not identifiable, the metrics proposed above (Eq.~\ref{eq:obj-trace} to~\ref{eq:obj-cond}) will be ill-posed. We recommend that the reader explore other methodologies to improve the conditioning of $\boldsymbol{M}$ via reformulation, choosing a different model, or leveraging additional prior information. 

Table~\ref{tab:criteria-summary} summarizes the practical interpretation of these criteria and when each objective may be preferred. Case Study 1 visualizes how these metrics change over an experimental design space and provides geometric understanding utilizing the eigenvalues of the system. This geometric understanding should help guide those performing experimental design to select the right metric(s) to be used during optimal experimental design.

\begin{table}[t]
\centering
\caption{Summary of the optimal design criteria discussed in this work.}
\label{tab:criteria-summary}
\small
\begin{adjustbox}{max width=\textwidth}
\begin{tabular}{p{2.3cm} p{4.4cm} p{7.0cm}}
\toprule
\textbf{Objective} & \makecell[l]{\textbf{Mathematical}\\\textbf{Interpretation}} & \textbf{Practical Interpretation and Use} \\
\midrule
A-optimality & Minimize $\sum_{i}^p 1/\lambda_i$ (Eq.~\ref{eq:obj-trace}), i.e., minimize $\text{trace}\left(\boldsymbol{M}^{-1}\right)$. & Because the inverse eigenvalues are summed, small eigenvalues contribute most strongly to the objective. In practice, this criterion tends to prioritize directions with low information (high uncertainty), and is useful when the primary goal is to improve poorly informed parameters. \\
\midrule
D-optimality & Maximize $\prod_i^p \lambda_i$ (Eq.~\ref{eq:obj-log-det}), equivalent to maximizing $\left|\boldsymbol{M}\right|$. & This criterion increases information volume overall; however, the multiplicative structure can favor directions that are already informative. It is often appropriate when broad information gain is desired, even if improvements are not evenly distributed across parameter directions. \\
\midrule
E-optimality & Maximize $\min_{i\in\left\{1, \ldots, p\right\}} \lambda_i$ (Eq.~\ref{eq:obj-min-eig}), i.e., maximize the minimum eigenvalue of $\boldsymbol{M}$. & The objective acts directly on the smallest eigenvalue of the FIM, so it targets the least informative direction. This is particularly useful when one or more parameter combinations remain weakly informed and the design objective is to improve worst-direction information content. \\
\midrule
ME-optimality & Minimize $\frac{\max_{i\in\left\{1, \ldots, p\right\}} \lambda_i}{\min_{i\in\left\{1, \ldots, p\right\}} \lambda_i}$ (Eq.~\ref{eq:obj-cond}), i.e., minimize the condition number of $\boldsymbol{M}$. & This criterion emphasizes balance in information across directions by reducing disparity between the largest and smallest eigenvalues. Since it does not directly maximize the overall magnitude of information, it is often most effective when used alongside another criterion that increases information content. \\
\midrule
Pseudo-A-optimality & Maximize $\sum_i^p \lambda_i$ (Eq.~\ref{eq:obj-pseudo-trace}), i.e., maximize $\text{trace}\left(\boldsymbol{M}\right)$. & Although this objective is computationally convenient, it is not equivalent to A-optimality and can emphasize already informative directions. It may still be useful as a surrogate in settings where numerical robustness is a dominant concern. \\
\bottomrule
\end{tabular}
\end{adjustbox}
\normalsize
\end{table}

\subsection{Challenging Design Criteria in Equation-Oriented Programs}
Knowledge of the form of these design criteria (Eq.~\ref{eq:obj-trace} through~\ref{eq:obj-cond}) is not sufficient to compute these metrics within a simultaneous, equation-oriented framework. We have equations to define the FIM (Eq.~\ref{eq:fim-fd} through~\ref{eq:cov-fim-relationship}), but we can not explicitly pose these design metrics in all cases. For instance, there is no general formula for inverting a matrix (e.g., to determine proper A-optimality, Eq.~\ref{eq:obj-trace}) or finding the eigenvalues of a matrix (e.g., to determine D-, E-, and ME-optimality).

One trick for posing D-optimality is to use a Cholesky factorization, which can be explicitly posed as algebra, to decompose the original FIM into an upper and lower triangular component~\cite{wang2022pyomo-doe}. This allows easy computation of the determinant, as the determinant of triangular matrices is trivial, by computing the product of the diagonal. Oftentimes, the log of the determinant is taken to improve numerical stability. Fortunately, the log does not complicate the algebra and can still be posed easily in equation-oriented programming. A-optimality can also be posed utilizing the Cholesky factorization, as shown recently in \cite{duarte2022A-opt}, although this A-optimality formulation is not currently available in Pyomo.DoE.

However, we cannot pose any algebraic constraints to determine the eigenvalues of a general matrix. Given that the eigenvalues can be found by solving a $p$-th order polynomial, root formulas exist for orders of $p$ up to four. Beyond four, however, there are proofs showing a general form does not exist \cite{abel1826degree5, ruffini1813degree5}. Also, adding a large system of equations to perform Cholesky factorization may be deleterious to the computational efficiency of solving the optimal experimental design problem due to the addition of potentially problematic nonlinear equality constraints to an already complicated, nonlinear model. To address these challenges, we utilize callbacks in Pyomo (via the pynumero package \cite{rodriguez2020pynumero}) to explicitly determine the value of each criteria (A-, D-, E-, ME-optimality) at each iteration. The only requirement beyond evaluating the criteria is that the derivative information---first derivative, Jacobian, and preferably second derivative, Hessian, as well---be supplied to the algebraic solver, in this case \texttt{IPOPT}, at each iteration. In Pyomo, we refer to these callbacks as an external input-output model, or an external \emph{Grey Box} model. Throughout the rest of this article, the term Grey Box refers to a functional callback (and other operations performed alongside the callback such as derivative evaluation) to embed a challenging (linear algebra) function within an equation-oriented optimization formulation.

To this end, mathematical program~\ref{eq:formulation-1} is altered to use callbacks in the objective function:

\begin{subequations} \label{eq:formulation-gb}
\begin{align}
    \min_{\boldsymbol{\phi} \in \boldsymbol{\Phi}} & \;\; \xrightarrow{\boldsymbol{M}_{\text{GB}}}\fcolorbox{black}{gray}{$\text{\color{white}\bf{Grey Box}}$}\xrightarrow{\Psi} \\
    \text{s.t.} & \;\; \boldsymbol{M}= \boldsymbol{Q}^T\boldsymbol{\Sigma}_{\hat{y}}^{-1}\boldsymbol{Q} + \boldsymbol{M}_{\hat{\boldsymbol{\theta}}} \\
    & \;\; \boldsymbol{M}= \boldsymbol{M}_{\text{GB}} \label{eq:FIM-greybox-constraint} \\
    &\;\;\boldsymbol{q}_{j}\left(\boldsymbol{\phi}, \boldsymbol{\theta}\right) = \frac{\hat{\boldsymbol{y}}\left(\boldsymbol{\phi}, \hat{\boldsymbol{\theta}} \;+ h_j\boldsymbol{e}_j\right) - \hat{\boldsymbol{y}}\left(\boldsymbol{\phi},\hat{\boldsymbol{\theta}} \;- h_j\boldsymbol{e}_j\right)}{2h_j} \\
    & \qquad\qquad\qquad\qquad\qquad\qquad\qquad\qquad\forall \; j \in \left\{1, \ldots, p\right\} \notag \\
    &\;\;\boldsymbol{Q}\left(\boldsymbol{\phi}, \hat{\boldsymbol{\theta}}\right) = 
        \begin{bmatrix}
            q_{1,1} & q_{1,2}  & \ldots & q_{1,p}   \\[0.35cm]
            q_{2,1} & q_{2,2}  & \ldots & q_{2,p}   \\[0.1cm]
            \vdots & \vdots & \vdots & \vdots  \\[0.1cm]
            q_{N_{\text{meas}},1} & q_{N_{\text{meas}},2}  & \ldots & q_{N_{\text{meas}},p}  \\
        \end{bmatrix}
\end{align}
\end{subequations}
\noindent All equations in formulation~\ref{eq:formulation-gb} are identical to formulation~\ref{eq:formulation-1}, except for the objective function and the constraint~\ref{eq:FIM-greybox-constraint}. In this case, we have a Grey Box model with the FIM as an input and the scalarized objective, $\Psi$, as an output. Constraint~\ref{eq:FIM-greybox-constraint} is required to algebraically link variables in the Pyomo model ($\boldsymbol{M}$) to the inputs ($\boldsymbol{M}_{\text{GB}}$) of the Grey Box object. We know exactly what is being computed and how it is being computed within the callback, but we have no way of representing it compactly as algebra (especially for E- and ME-optimality). Therefore, the information we need to supply to \texttt{IPOPT} at each iteration is the value of the Grey Box output, $\Psi$, and the first (and optionally second) derivative of the Grey Box output, $\Psi$, with respect to the Grey Box input(s), $\boldsymbol{M}$. As with our previous implementation, optimization formulation~\ref{eq:formulation-gb} is automatically generated if the user specifies that they wish to use the Grey Box objective (i.e., \texttt{use\_grey\_box\_objective = \color{blue}{True}}).

\subsubsection{First Derivative and Second Derivative of Design Criteria}
This section will briefly outline the formulas for the first and second derivatives, the Jacobian, and Hessian, respectively, of each scalarized objective function, $\Psi$, with respect to the FIM, $\boldsymbol{M}$ (note that since constraint~\ref{eq:FIM-greybox-constraint} is enforced, the derivative of $\Psi$ is equivalent with respect to $\boldsymbol{M}$ or $\boldsymbol{M}_{\text{GB}}$). The supplementary information provides a detailed derivation. 

Two additional notes: (i) the Hessian has indices $\left(i, j, k, l\right)$ as it is a 4-th order tensor, and (ii) symmetry is enforced by modeling only the upper triangular portion of the FIM in constraint~\ref{eq:FIM-greybox-constraint}. All implementation details are included in the supplementary information.

For A-optimality, the first derivative can be found using the following:

\begin{align}
    \frac{\partial\,\text{trace}\left(\boldsymbol{M}^{-1}\right)}{\partial \boldsymbol{M}} & = - \left(\boldsymbol{M}^{-1} \boldsymbol{M}^{-1}\right)^T \label{eq:A-opt-first-deriv}
\end{align}

Similarly, the second derivative is shown below:

\begin{align}
    \frac{\partial}{\partial M_{kl}}\left(\frac{\partial\,\text{trace}\left(\boldsymbol{M}^{-1}\right)}{\partial \boldsymbol{M}}\right)_{ij} & = -\left( M_{il}^{-1} \left(M^{-1} M^{-1}\right)_{kj} + \left(M^{-1} M^{-1}\right)_{il} M_{kj}^{-1}\right) \label{eq:A-opt-second-deriv}
\end{align}
where $M_{ij}$ represents the element of the FIM, $\boldsymbol{M}$, in the $i$-th row and $j$-th column. In these equations, we utilize the \texttt{linalg.pinv} functionality in \texttt{NumPy} \cite{numpy} to compute the pseudo-inverse of the FIM. The pseudo-inverse ensures numerical stability because at any given optimization iteration, the FIM may be rank deficient.

For D-optimality, we use the log of the determinant. The first derivative can be found using the following, requiring the chain rule:

\begin{align}
    \frac{\partial\,\text{ln}\left|\boldsymbol{M}\right|}{\partial \boldsymbol{M}} & = \frac{1}{2} \left(\boldsymbol{M}^{-1} + \boldsymbol{M}^{-T}\right) =\boldsymbol{M}^{-1}\label{eq:D-opt-first-deriv}
\end{align}

Taking the second derivative results in the formula below:
\begin{align}
    \frac{\partial}{\partial M_{kl}}\left(\frac{\partial\,\text{ln}\left|\boldsymbol{M}\right|}{\partial \boldsymbol{M}}\right)_{ij} & = M_{il}^{-1}M_{kj}^{-1} \label{eq:D-opt-second-deriv}
\end{align}
Similarly to A-optimality, we also compute the psuedo-inverse of the FIM with \texttt{NumPy} for the first and second derivative of D-optimality. Also, the log-determinant of the FIM is computed within \texttt{NumPy} using the signed log-determinant function (\texttt{linalg.slogdet}) to avoid a negative determinant value within the natural log function.

For E-optimality, we start with the eigenvalue problem on $\boldsymbol{M}$ and utilize this to compute the first derivative:

\begin{align}
    \boldsymbol{M}\boldsymbol{v}_s& = \lambda_s\boldsymbol{v}_s &\forall s \in \left\{1, \ldots, p\right\} \label{eq:eigenvalue_problem} \\
    \boldsymbol{v}_s\cdot \boldsymbol{v}_s & = 1 &\forall s \in \left\{1, \ldots, p\right\}\label{eq:eigenvector_orthonormality} \\
    \boldsymbol{v}_i\cdot \boldsymbol{v}_j & = 0 &\forall i \neq j; \left(i, j\right) \in \left\{1, \ldots, p\right\}^2\label{eq:eigenvector_orthonogonality} \\
    \frac{\partial\,\boldsymbol{M}\boldsymbol{v}_\text{min}}{\partial \boldsymbol{M}} & = \frac{\partial\,\lambda_\text{min}\boldsymbol{v}_\text{min}}{\partial \boldsymbol{M}}\label{eq:E-opt-first-deriv-1} \\
    \frac{\partial\,\lambda_\text{min}}{\partial \boldsymbol{M}} & = \boldsymbol{v}_\text{min} \boldsymbol{v}_\text{min}^T\label{eq:E-opt-first-deriv-2}
\end{align}
where $s$ is the index of some eigenvalue and eigenvector of the FIM. Note, the number of eigenvalues matches the square dimension of the FIM, or the total number of parameters, $p$. In addition, the subscript ``min'' refers to the index of the minimum eigenvalue and the corresponding eigenvector.

The second derivative can then be found using the following equation: 
\begin{align}
    \frac{\partial}{\partial M_{kl}}\left(\frac{\partial\,\lambda_\text{min}}{\partial \boldsymbol{M}}\right)_{ij} = & \sum_{s \neq \text{min}} \left(\frac{1}{\lambda_\text{min} - \lambda_s} v_{s, k} v_{\text{min}, l} v_{s, j} \right) v_{\text{min}, i} \notag\\ & + v_{\text{min}, j} \sum_{s \neq \text{min}} \left(\frac{1}{\lambda_\text{min} - \lambda_s} v_{s, l} v_{\text{min}, k} v_{s, i}\right)\label{eq:E-opt-second-deriv}
\end{align}
See the supplementary material for a step-by-step derivation. These derivative expressions require both the eigenvalues and the eigenvectors of the FIM, which are calculated using the \texttt{linalg.eig} functionality in \texttt{NumPy}. Also, to compute $\Psi$, we directly find the minimum eigenvalue using \texttt{NumPy}.

For ME-optimality, we follow a similar path to E-optimality with the added chain rule. We note here that we once again use a log transformation for the condition number, which makes the objective better scaled. It follows that the first derivative is:

\begin{align}
    \frac{\partial\,\text{ln}\left(\frac{\lambda_\text{max}}{\lambda_\text{min}}\right)}{\partial \boldsymbol{M}} & = \frac{1}{\lambda_\text{max}} \frac{\partial \lambda_\text{max}}{\partial \boldsymbol{M}} - \frac{1}{\lambda_\text{min}} \frac{\partial \lambda_\text{min}}{\partial \boldsymbol{M}} \label{eq:ME-opt-first-deriv}
\end{align}
where the subscript ``max'' refers to the index of maximum eigenvalue and corresponding eigenvector. Eq.~\ref{eq:E-opt-first-deriv-2} is valid for any eigenvalue and eigenvector combination, $s$, and we have the pieces to define Eq.~\ref{eq:ME-opt-first-deriv} using Eq.~\ref{eq:E-opt-first-deriv-2} for both the ``min'' and ``max'' eigenvalue/eigenvector pairs. Here, we note that the FIM is positive semi-definite, so the typical absolute value included in the definition of matrix condition number is not required; however, we ensure positivity of the condition number in our implementation.

The second derivative utilizes the chain rule again and results in the following:

\begin{align}
    \frac{\partial}{\partial M_{kl}}\left(\frac{\partial\,\text{ln}\left(\frac{\lambda_\text{max}}{\lambda_\text{min}}\right)}{\partial \boldsymbol{M}}\right)_{ij} = &\;\;\;\;\frac{1}{\lambda_\text{max}}  \frac{\partial}{\partial M_{kl}}\left(\frac{\partial \lambda_\text{max}}{\partial \boldsymbol{M}}\right)_{ij} \notag\\
    & - \frac{1}{\lambda_\text{max}^2} \left(\frac{\partial \lambda_\text{max}}{\partial \boldsymbol{M}}\right)_{kl}\left(\frac{\partial \lambda_\text{max}}{\partial \boldsymbol{M}}\right)_{ij}\notag \\ 
    & + \frac{1}{\lambda_\text{min}^2} \left(\frac{\partial \lambda_\text{min}}{\partial \boldsymbol{M}}\right)_{kl} \left(\frac{\partial \lambda_\text{min}}{\partial \boldsymbol{M}}\right)_{ij} \label{eq:ME-opt-second-deriv} \\ 
    & - \frac{1}{\lambda_\text{min} ^ 2} \frac{\partial}{\partial M_{kl}}\left(\frac{\partial \lambda_\text{min}}{\partial \boldsymbol{M}}\right)_{ij} \notag
\end{align}
Eq.~\ref{eq:E-opt-second-deriv} is also valid for any unique eigenvalue; thus, Eq.~\ref{eq:E-opt-first-deriv-2} and Eq.~\ref{eq:E-opt-second-deriv} can be used to find an expression for the first and second derivatives of the maximum eigenvalue to compute all terms in Eq.~\ref{eq:ME-opt-second-deriv}. Once again, these derivative expressions require both the eigenvalues and eigenvectors, which are computed using the \texttt{linalg.eig} functionality in \texttt{NumPy}. Finally, the condition number is computed using the maximum and minimum eigenvalues, also using \texttt{NumPy}. 

As a brief aside we note that using the \texttt{linalg.cond} function has slightly different numerics than using the maximum and minimum eigenvalues from \texttt{linalg.eig} because the functions follow different numerical pathways. It is important to be consistent with which numerical method is used when computing these metrics; otherwise, the numerical error introduced while computing the analytical derivative may exceed the threshold of using a more crude, finite difference scheme to find a numerical derivative which ultimately may lead the optimizer astray.

Analytical derivatives are highly favored over numerical derivatives in this situation as: (i) the linear algebra system needs to be solved only once instead of multiple times, and (ii) the numerical error introduced from a finite difference scheme may be too large, hindering the optimizer (e.g., \texttt{IPOPT}).

Given that we have a method to compute both the value of $\Psi$ and the derivative(s) of these four experimental design criteria (A-, D-, E-, and ME-optimality), we can utilize the \texttt{ExternalGreyBox} feature within Pyomo \cite{rodriguez2020pynumero, bynum2021pyomo} to embed our Grey Box objective within an equation-oriented programming architecture. This requires using \texttt{cyipopt} \cite{cyipopt} to interface with \texttt{IPOPT}. \texttt{cyipopt} allows for functional callbacks  (e.g., to \texttt{NumPy}) to be evaluated and embedded within a solver call to \texttt{IPOPT} so long as the first derivative (Jacobian) information, and, optionally, the second derivative (Hessian) is provided.

More information on reproducing the results in this document, including the versions of software used and what environments are required, is included in the supplementary information. 

\section{Unifying Software Abstraction - The Experiment Class} \label{sec:experiment-abstraction}
An important characteristic when considering open-source software development is ease of use for the target audience. In this case, bringing optimal experimental design with minimal effort to scientists and engineers (end users) is key. In this section, we will describe a unifying model abstraction to: (i) reduce the barrier to entry for end users trying to perform optimal experimental design on their own models, (ii) standardize the workflow within Pyomo to connect all elements of the model building workflow (Figure~\ref{fig:workflow}) with one, convenient object, and (iii) provide an easy avenue to promote external development from community members for open-source software contribution.

\begin{figure}[t]
\includegraphics[width = \textwidth]{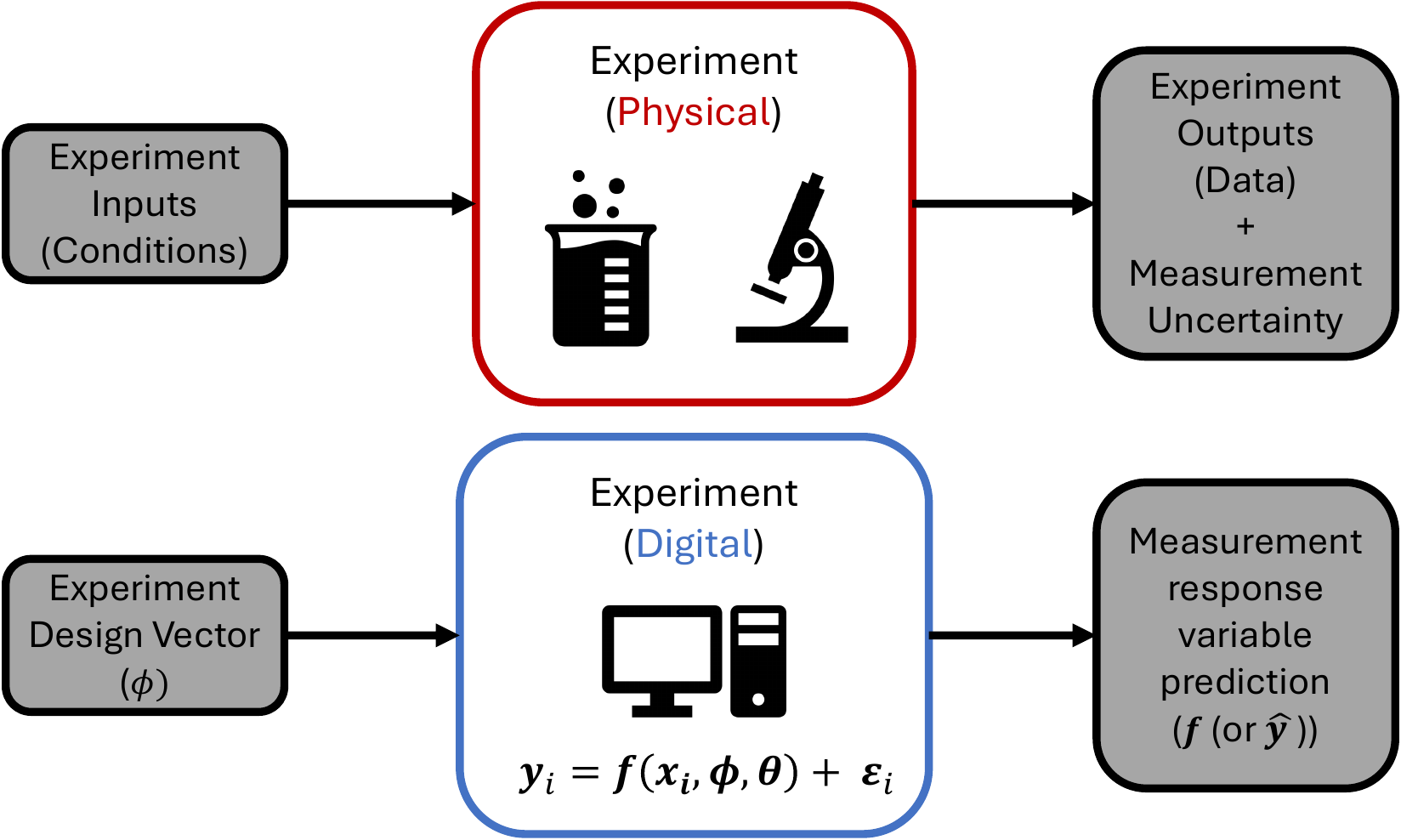}
\caption{Parallelism between the physical experiment that a scientist performs with the digital \texttt{Experiment} that enables the model building workflow from Figure~\ref{fig:workflow} within Pyomo.} \label{fig:experiment-class}
\end{figure}

For this purpose, we developed an abstraction that we call the \texttt{Experiment} class. The idea is to connect the physical experiment to the digital model simulating the physical system that is interrogated during model-building tasks. In Figure~\ref{fig:experiment-class}, the connection between the physical decisions made and the digital decisions required to define an \texttt{Experiment} object are through the parallel idea of the inputs to the process (left side of Figure~\ref{fig:experiment-class}). Similarly, the model prediction and the experimental data are the outputs of the digital and physical systems, respectively (right side of Figure~\ref{fig:experiment-class}). The digital \texttt{Experiment} requires a candidate model representation of the physical process that is labeled by the scientist to point the automated modeling tools in Pyomo to perform each step of the model building and validation workflow in Figure~\ref{fig:workflow}.

\FloatBarrier

\begin{figure}[H]
\includegraphics[width = \textwidth]{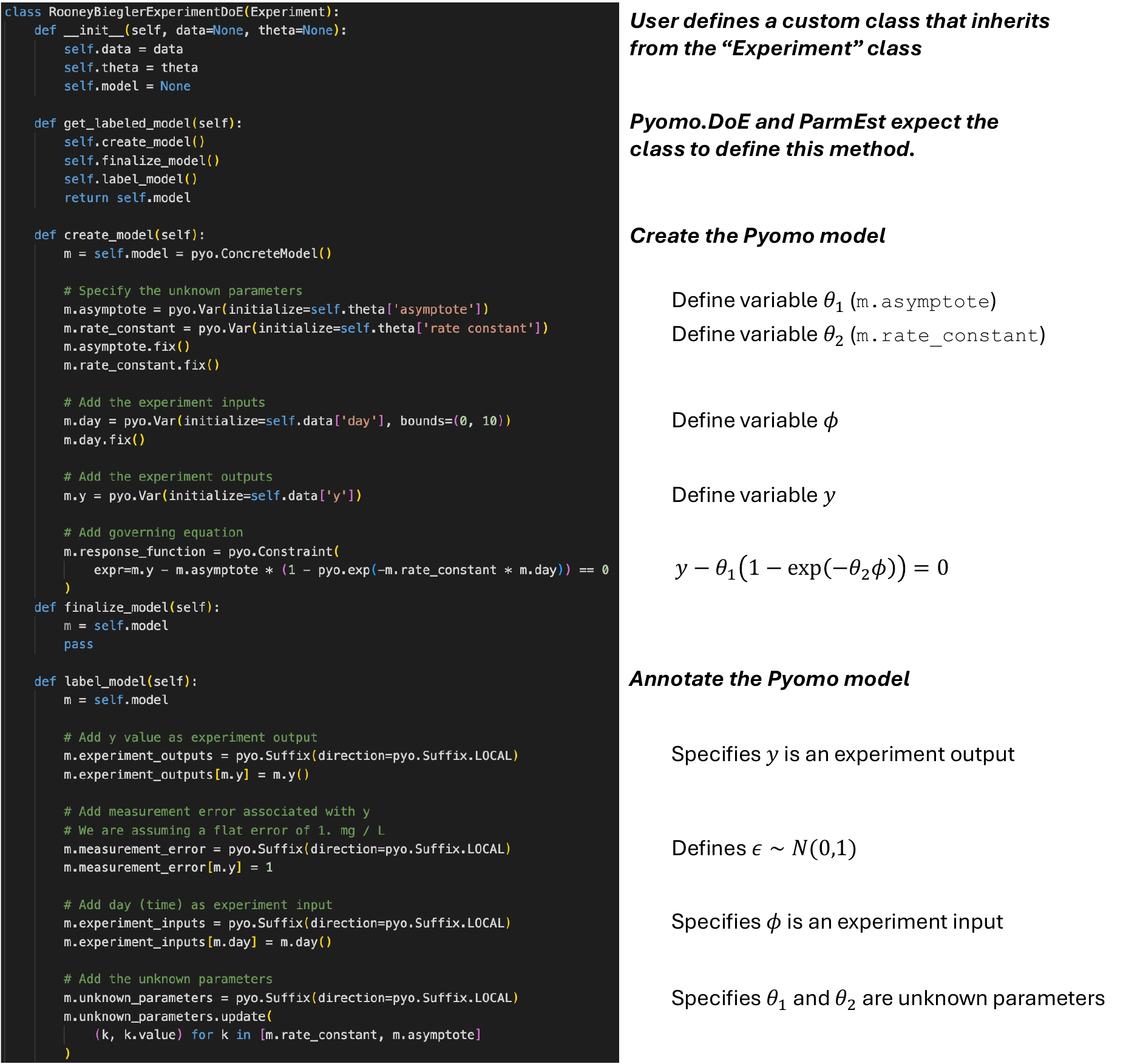}
\caption{Example connecting the elements of a small model (in Case Study 1, Eq.~\ref{eq:batch-reaction}) to the code required to label the important components of the model using the \texttt{label\_model} function of an \texttt{Experiment} object.}
\label{fig:general-experiment-code-snippet}
\end{figure}

When a scientist chooses the conditions under which they carry out their physical experiment, they must consider how the numerical counterpart of these conditions must be specified in the model. We call these conditions \texttt{experiment\_inputs} ($\boldsymbol{\phi}_i$ in Eq.~\ref{eq:general-model}). Similarly, the scientist also needs to specify which digital components of the model are the predicted values of the measured data in the physical system, which we call \texttt{experiment\_outputs} ($\boldsymbol{y}_i$ in Eq.~\ref{eq:general-model}). Importantly, the scientist must also provide a measurement error (which we call \texttt{measurement\_error}; $\boldsymbol{\varepsilon}_i$ in Eq.~\ref{eq:general-model}) associated with these outputs as these are a key component in generating adequate information and covariance calculations (Eqs. \ref{eq:measurement-error-matrix-def}, \ref{eq:cov-def}, and \ref{eq:cov-fim-relationship}). Finally, the unknown parameters in the model, $\boldsymbol{\theta}$, must also be identified by the scientist when defining the model, which we call \texttt{unknown\_parameters}. An example of the additional code required to label the model is given in Figure~\ref{fig:general-experiment-code-snippet} to illustrate ease of use.

In this way, the scientist has a labeled model which may be reused in any Pyomo contributed packages which expect this set of model labels. For instance, someone performing least squares parameter estimation must build an objective function to minimize the difference between the model prediction and the physical data. However, since the scientist has already labeled the relevant information for the unknown parameters (\texttt{unknown\_parameters}), the output data (\texttt{experiment\_outputs}), and the measurement error associated with the output data (\texttt{measurement\_error}), Pyomo's parameter estimation package, \texttt{parmest}, can automatically create the objective function (Eq.~\ref{eq:param-estimation}) and formulate the optimization problem to estimate the unknown parameters that best fit the model to the data. Similarly, Pyomo.DoE formulates the optimal experimental design problem with the preferred objective function (formulation~\ref{eq:formulation-1} or formulation~\ref{eq:formulation-gb}) and optimizes experimental design using a prior estimate of the parameters and associated prior covariance matrix. By simply adding these labels to the model, as shown in Figure~\ref{fig:general-experiment-code-snippet}, the user can now formulate equation-oriented mathematical programs for parameter estimation and optimal experimental design problems without changing the model structure to include Eq.~\ref{eq:param-estimation} and formulations~\ref{eq:formulation-1} and~\ref{eq:formulation-gb}, respectively.

This modeling abstraction is standard for intrusive uncertainty quantification and experiment design in Pyomo. It was used exclusively to generate the results shown in the following section and has been a monumental component in getting others to adopt Pyomo.DoE in their respective model-building workflows. This illustrates how dedicating a small amount of time to the software design elements of open-source tools goes a long way toward user adoption and impact in the broader scientific community. For those curious, we have an open-source tutorial for those wanting to use Pyomo.DoE and \texttt{parmest} in their own workflows \url{https://dowlinglab.github.io/pyomo-doe}. The tutorial covers how to pose models within Pyomo and how to appropriately label these models to be used in Pyomo's model building workflow (Figure~\ref{fig:workflow}).

\section{Results and Discussion} \label{sec:results}
In this section, we will show three case studies that demonstrate the new capabilities in Pyomo.DoE for complex design metrics using callbacks. First, a small example with two unknown parameters to design the time of reaction for a batch-style reaction system. Second, we utilize a linear control system. Finally, we showcase the new capabilities on a design problem for a membrane cascade system to recover critical minerals from recycled battery waste.

\subsection{Case Study 1: Batch Reaction System} \label{sec:res-case-1}
The first example is adapted from data found in Bates and Watts~\cite{bates1988nonlinear} (section A1.4), which was originally found in Marske~\cite{marske1967data}. In this problem, samples are prepared and periodically examined to gather dissolved oxygen concentration to understand biological oxygen demand in a physical system. The system can be described with the following equation:

\begin{align}
    f\left(\phi, \boldsymbol{\theta}\right) = \theta_1\left(1 - \text{exp}\left(-\theta_2\phi\right)\right) \label{eq:batch-reaction}
\end{align}
where $\phi$ is the sample time and $\theta_1$ and $\theta_2$ are unknown parameters which need to be estimated from data. 

We consider that there are two experiments that have already been run: those listed in Table~\ref{tab:batch-reaction-data}. For optimal experimental design, we wish to design an experiment for the optimal time to measure the sample.

\begin{table}[!ht]
\centering
\caption{Preliminary data for batch reaction case study taken from A1.4 in Bates and Watts \cite{bates1988nonlinear}.}
\begin{tabular}{|c|c|}
    \hline
    \bfseries Sample Time (day) & \bfseries Output data (mg $\cdot$ L$^{-1}$) \\ \hline \hline 
    1 & 8.3 \\ \hline 
    7 & 19.8 \\ \hline 
\end{tabular} 
\label{tab:batch-reaction-data}
\end{table}
\FloatBarrier

Moreover, we assume that these measurements are independent and have known measurement error of 1 (mg $\cdot$ L$^{-1}$). Using these experimental data, we use nonlinear regression within Pyomo (the \texttt{parmest} contributed package) to find the initial best estimates for $\theta_1$ and $\theta_2$ as 20.3 mg $\cdot$ L$^{-1}$ and 0.53 day$^{-1}$, respectively. From the nonlinear regression problem, we obtain a covariance matrix for these unknown parameters (shown below in Eq.~\ref{eq:cov-batch-reaction}), which becomes the prior for optimal experimental design (Eq.~\ref{eq:cov-def}).

\begin{align}
&\;\;\boldsymbol{V}_{\boldsymbol{\theta}} = 
    \begin{bmatrix}
        4.03 \cdot 10^{-3} & -8.27 \cdot 10^{-3}   \\[0.35cm]
        -8.27 \cdot 10^{-3}  & 42.8 \cdot 10^{-3}   \\
    \end{bmatrix} \label{eq:cov-batch-reaction}
\end{align}

\begin{figure}[H]
\includegraphics[width = \textwidth]{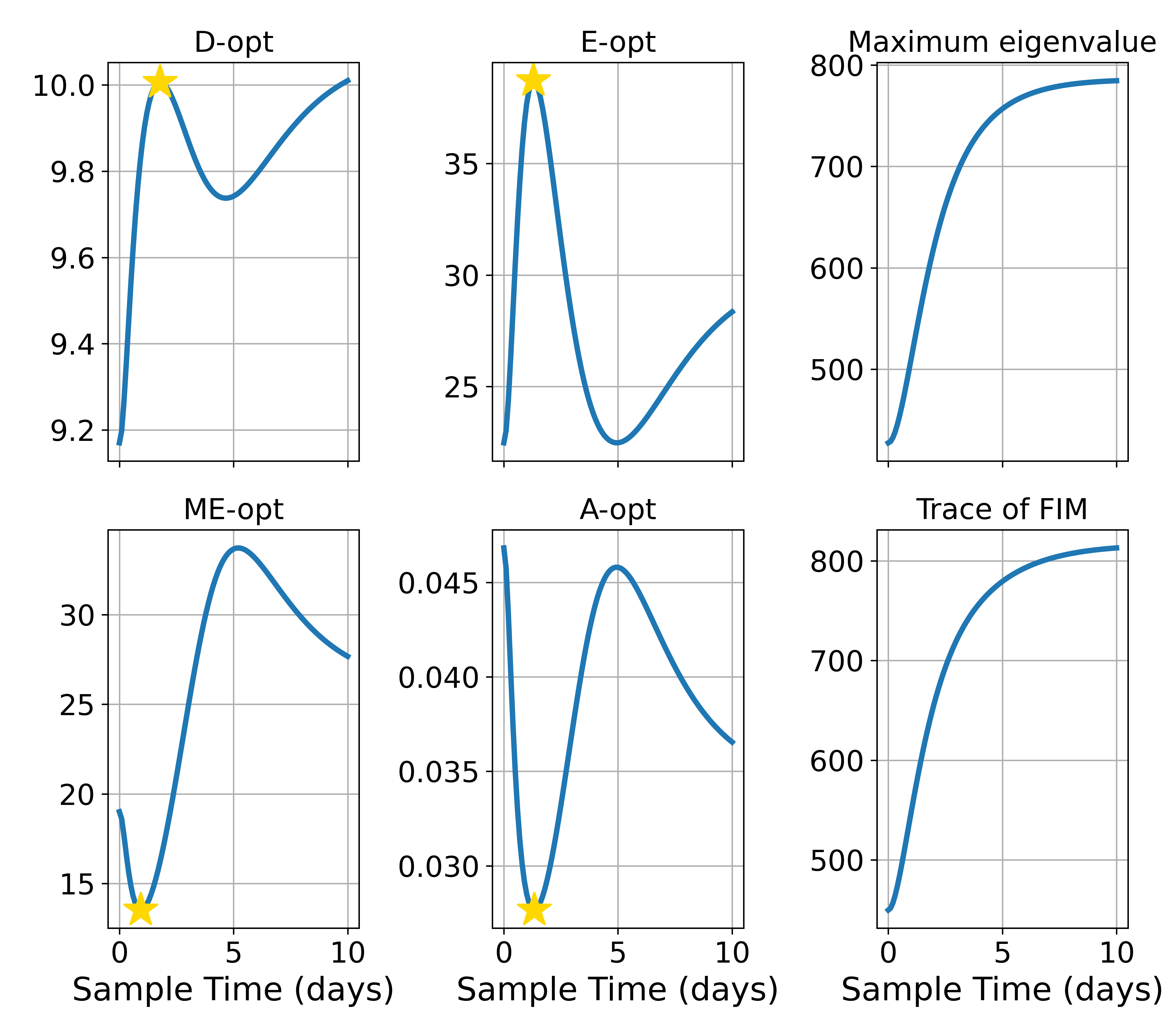}
\caption{Objective function values, $\Psi$ (D-, E-, ME-, A-opt), and related values (maximum eigenvalue and trace of the FIM) plotted as a function of the experimental design decision, sample time (days). The star represents the optimal point found for each criterion using Pyomo.DoE with a Grey Box objective.} \label{fig:batch-reaction-comparison}
\end{figure}

To improve the quality of fit, we will design an additional experiment that takes between 1 and 10 days. For this, we utilize Pyomo.DoE to compute optimal designs for each scalarized objective. Since this experimental design problem has only 1 decision variable, we can visualize how each information metric (e.g., A-, D-, E-, and ME-optimality) and the eigenvalues of the FIM change with respect to the experimental design. These results (both the optimal design point and sensitivity analysis) are displayed in Figure~\ref{fig:batch-reaction-comparison}, with optimal points listed in Table~\ref{tab:batch-reaction-optima}.

\FloatBarrier

\begin{table}[!ht]
\centering
\caption{Optimal sample time (days) for each criteria.}
\begin{tabular}{|c|c|}
    \hline
    \bfseries Optimality Criteria & \bfseries Sample Time (day) \\ \hline \hline 
    D-opt   & 1.78     \\ \hline
    A-opt   & 1.33     \\ \hline
    E-opt   & 1.30     \\ \hline
    ME-opt  & 0.94    \\ \hline
\end{tabular} 
\label{tab:batch-reaction-optima}
\end{table}

As shown in Table~\ref{tab:batch-reaction-optima}, the optimal design for each criterion is within a range of 1 day (from about 0.9 to 1.8 days) of the other designs. D-optimality typically focuses on the largest eigenvalue \cite{franceschini2008MBDoE, montgomery2017design} (maximum at $\phi=10$ days), but in this case, the largest increase in determinant has more to do with the sharp increase in the minimum eigenvalue. Importantly, the A-optimality curve closely follows that of E-optimality in Figure~\ref{fig:batch-reaction-comparison}. This is because the sum of the inverse of the eigenvalues is most impacted when the smallest eigenvalue changes. Given that the minimum eigenvalue here is one order of magnitude lower than the maximum, it dominates the summation in Eq.~\ref{eq:obj-trace}.

Furthermore, the condition number (ME-optimality), or ratio between the two eigenvalues, is shown to be lowest at $\phi=0.94$ days. For D-optimality and E-optimality, a time of 0.94 days is suboptimal, leading to less information. This indicates that ME-optimality should be used carefully, for instance if some previous knowledge of the system indicates matrix conditioning is important, or when the condition number exceeds values of $10^{w}$ where $w$ is a known numerical stability or sloppiness threshold (e.g., $w=3$, \cite{gutenkunst2007sloppy-params}).

We also show the behavior of the trace of the FIM (pseudo-A-optimality) in the bottom right panel of Figure~\ref{fig:batch-reaction-comparison}. As shown in Figure~\ref{fig:batch-reaction-comparison}, pseudo-A-optimality does not follow the same trend as the proper definition of A-optimality (Eq.~\ref{eq:obj-trace}, bottom middle panel of Figure~\ref{fig:batch-reaction-comparison}). Although the trace of the FIM does not emulate the behavior of A-optimality, the trace of the FIM does not require inversion or eigenvalues, making it numerically attractive in near-singular cases. Also, as one might expect, the trace of the FIM is dominated by the behavior of the maximum eigenvalue, as the curves (maximum eigenvalue: top right panel, trace of FIM: bottom right panel) are almost indistinguishable. We also note that even a small system such as this (one-equation model with only one design decision) is non-convex and displays multiple local optima in all four of the desired criteria. This is important in practice because gradient-based solvers may converge to different locally optimal experiments depending on initialization. Therefore, robust OED workflows should use multi-start strategies (or other globalization safeguards) to avoid over-interpreting any single local solution.

\subsection{Case Study 2: Linear Control System, TCLab}
The second example is a linear control system; specifically, the temperature control lab (TCLab) \cite{oliveira2019tclab, park2020tclab, hedengren2020tclab, oliveira2020tclab}. In the undergraduate controls course at Notre Dame, we utilize a hands-on system known as the TCLab to teach control theory using state-space models in Python \cite{dowling2025teaching}. However, the device is useful beyond education within research as a small-scale system that can easily generate data and has well-known equations that dictate system physics. The device can be modeled as a two-body heating system in the equations below:

\begin{subequations}
\begin{align}
    C_p^H \frac{\text{d}T_H}{\text{d}t} = & \; U_a \left(T_\text{amb} - T_H\right) + U_b \left(T_S - T_H\right) + \alpha P u\left(t\right) \\
    C_p^S \frac{\text{d}T_S}{\text{d}t} = & \; U_a \left(T_H - T_S\right)
\end{align}
\end{subequations}

\noindent where the temperature of the sensor, $T_S$, and the heater, $T_H$, are modeled as the two states. Here, the ambient temperature of the system is $T_\text{amb}$. We can measure $T_S$ using the device, and therefore we utilize $T_S$ as our measurements and assume that the measurement error associated with the measurements is 0.25$^{\circ}$C (or equivalently, K). 

Energy is transferred to the system at time $t$ through the control function $u$ where electrical energy heats the system dictated by material specific coefficients, $\alpha$ and $P$, which are known beforehand ($\alpha = 0.00016, P=200$). It should be noted that our control decisions $u$ are the experimental decisions during optimal experimental design. Finally, there are four unknown parameters in the system. First, the heat transfer coefficients for transfer between the heater and the sensor, $U_b$, and between the ambient environment and the heater, $U_a$. And second, the specific heat of the heater and sensor represented by $C_p^H$ and $C_p^S$, respectively. Generally, the script $H$ denotes variables and parameters related to the heater and the script $S$ denotes variables and parameters related to the sensor.

To improve scaling, we reparameterize using $\beta_i$ as shown in the equations below:

\begin{subequations}
\begin{align}
    \frac{\text{d}T_H}{\text{d}t} = & \; \beta_1 \left(T_\text{amb} - T_H\right) + \beta_2 \left(T_S - T_H\right) + \beta_4 u\left(t\right) \\
    \frac{\text{d}T_S}{\text{d}t} = & \; \beta_3 \left(T_H - T_S\right) \\
    \beta_1 = & \; \frac{U_a}{C_p^H}, \qquad \beta_2 = \; \frac{U_b}{C_p^H} \qquad 
    \beta_3 = \; \frac{U_b}{C_p^S}, \qquad \beta_4 = \; \frac{\alpha P}{C_p^H}
\end{align}
\end{subequations}
This reparameterization shows that the TCLab can be modeled as a linear time-invariant (LTI) system.

\begin{figure}[H]
\includegraphics[width = \textwidth]{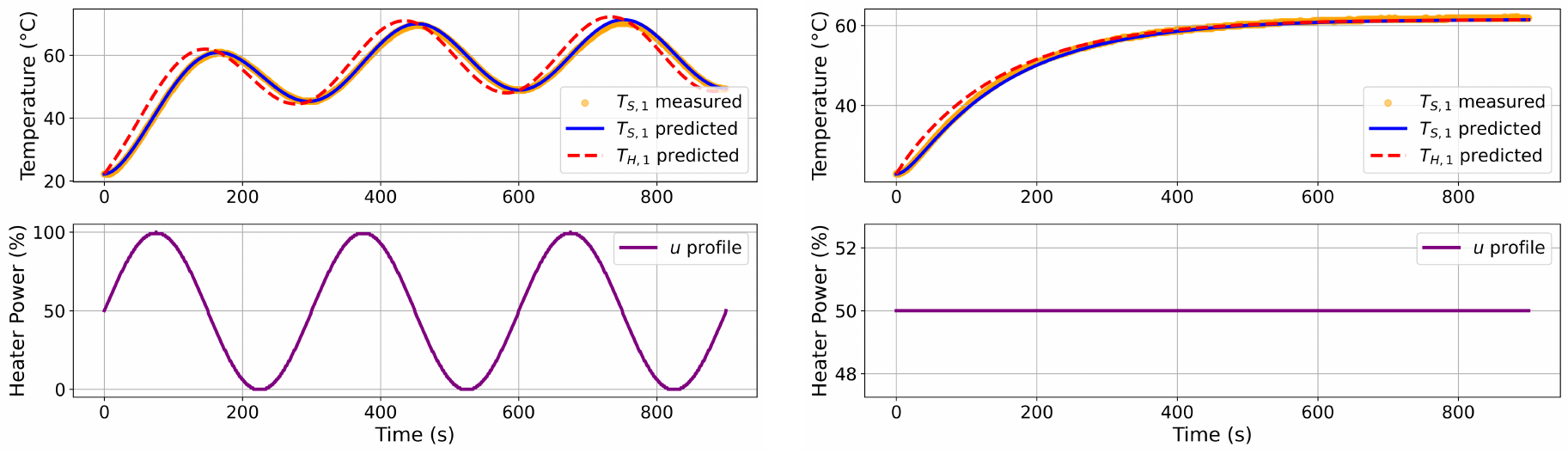}
\caption{Data for the sine wave (left) and step test (right) experiments using the TCLab system. Measured data (orange points), predicted profile using optimal parameters (solid blue line), and predicted heater temperature (dashed red line) are shown on the top of each subplot where the heater input is shown as a solid red line on the bottom of each subplot.} \label{fig:experiment-data-tclab}
\end{figure}

To emulate the model building workflow in Figure~\ref{fig:workflow}, we begin by defining a couple of initial experiments. First, we utilize a sine wave experiment in which the values of $u\left(t\right)$ are a sine wave with a period of 5 minutes and amplitude of 0.5 (heater power varies from 0 to 1 in a sinusoidal manner, as shown on the left of Figure~\ref{fig:experiment-data-tclab}). The second experiment is a simple step test, where the device is held at 50\% power for the full experiment duration (right side of Figure~\ref{fig:experiment-data-tclab}). Both experiments have common exploratory control profiles, are 15 minutes in duration (or 900 seconds on the plots), and were run on the same TCLab device. Also, the device must not exceed 85$^{\circ}$C to avoid being damaged. When performing nonlinear regression with \texttt{parmest}, the optimal values for the original parameters are shown in Table~\ref{tab:optimal-parameter-values-tclab}.

\FloatBarrier

\begin{table}[!ht]
\centering
\caption{Optimal parameter values for the TCLab heat transfer model using experimental data from Figure~\ref{fig:experiment-data-tclab} where the uncertainty from the covariance matrix is shown as using pair-wise covariance in Figure~\ref{fig:tclab-pairwise-uncertainty-combined}a).}
\begin{tabular}{|c|c|}
    \hline
    \bfseries Parameter & \bfseries Optimal Value \\ \hline \hline 
    $U_a$   & 0.0418    [W$\cdot^\circ$C$^{-1}$]  \\ \hline
    $U_b$   & 0.0303   [W$\cdot^\circ$C$^{-1}$]   \\ \hline
    $C_p^H$   & 5.487  [J$\cdot^\circ$C$^{-1}$]    \\ \hline
    $C_p^S$  & 0.588 [J$\cdot^\circ$C$^{-1}$]    \\ \hline
\end{tabular} 
\label{tab:optimal-parameter-values-tclab}
\end{table}
\FloatBarrier

\begin{figure}[H]
\includegraphics[width = \textwidth]{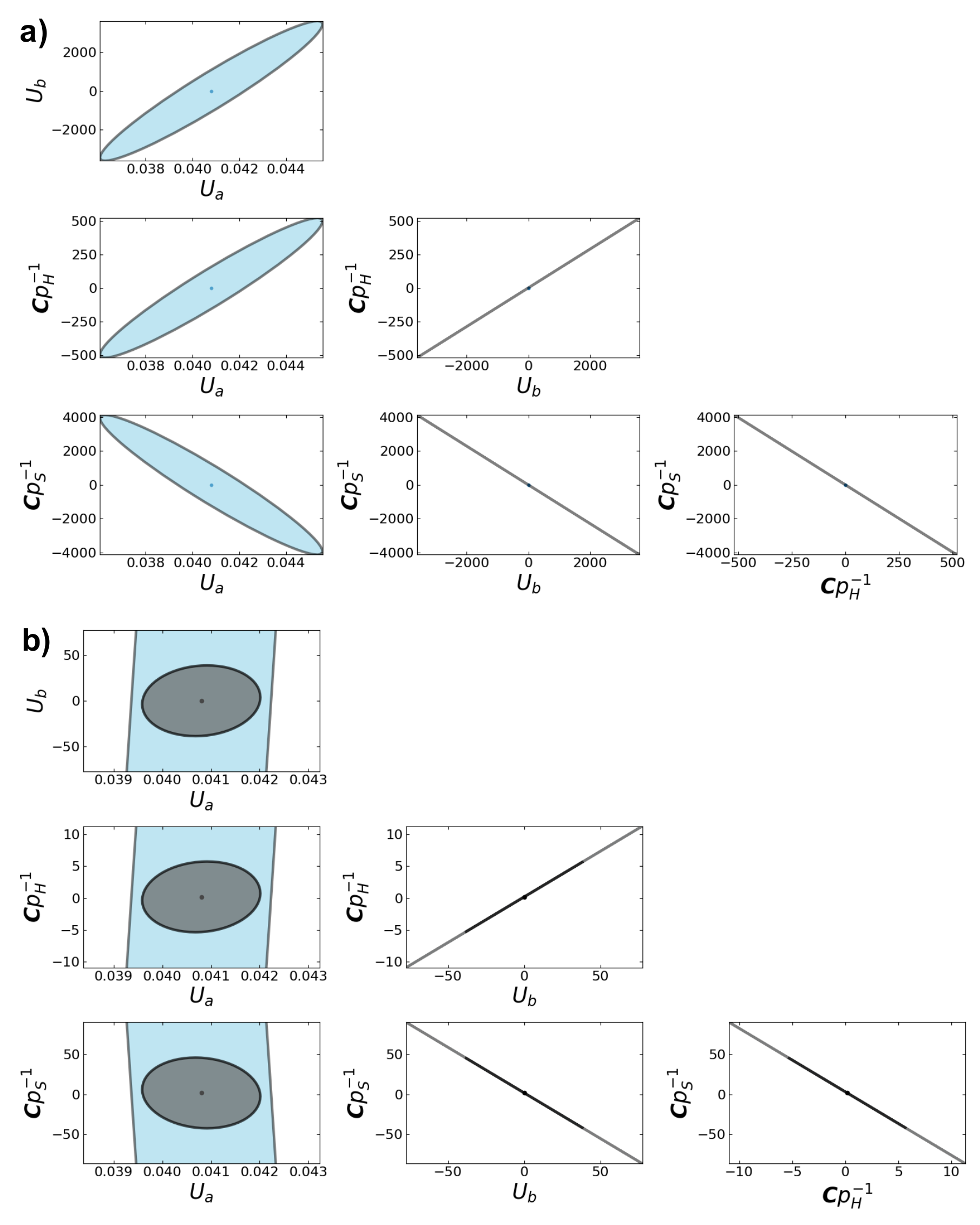}
\caption{Pairwise covariance of parameters with preliminary data only (a) and additional data (b) in the TCLab case study. With additional data (b), uncertainty after the optimal experiment (gray shading with black border) drastically reduces uncertainty (2 to 4 orders of magnitude).} \label{fig:tclab-pairwise-uncertainty-combined}
\end{figure}
\FloatBarrier

\begin{figure}[H]
\includegraphics[width = \textwidth]{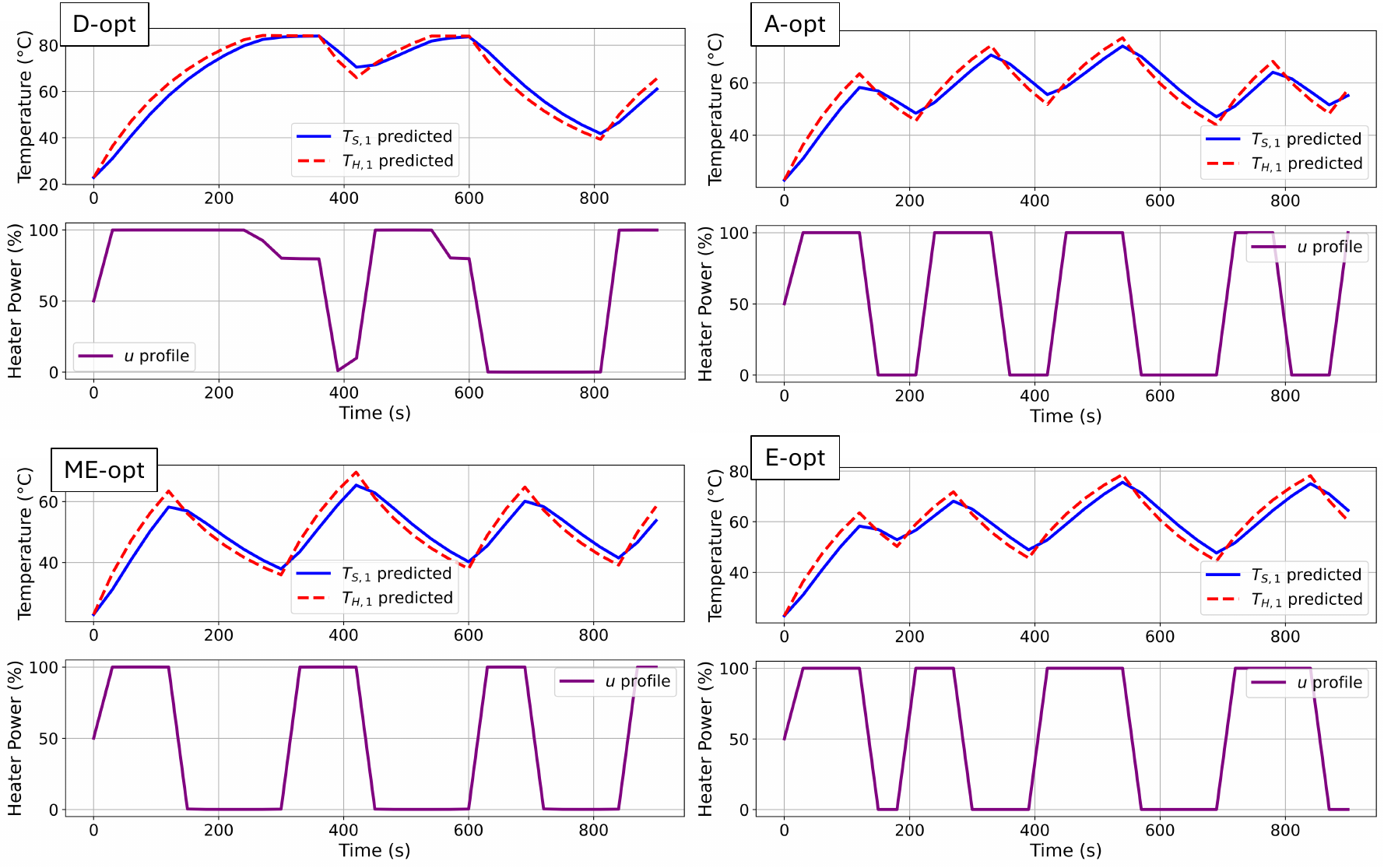}
\caption{Optimal profiles for D-optimality (top left), A-optimality (top right), ME-optimality (bottom left), and E-optimality (bottom right). In each subplot, the predicted system behavior is shown in the top panel and the experimental design of the control profile is shown in the bottom panel.} \label{fig:tclab-results}
\end{figure}

After parameter estimation, we wish to design an experiment to reduce parametric uncertainty. The experimental design decisions are control inputs at 30 second intervals starting at $t=0$ for a 900 second experiment (30 design decisions). Figure~\ref{fig:tclab-results} shows the optimal experimental design for each of the four design criteria. Each profile is similar in that a series of ``on-off'' step-tests are performed; however, the profiles are unique due to varying frequency of these steps and different maximum temperature achieved during the experiment.

Figure~\ref{fig:tclab-results} provides the control profiles, while Table~\ref{tab:optimal-tclab-criteria} quantifies their information content. To construct Table~\ref{tab:optimal-tclab-criteria}, we first solve four separate OED problems (one each for D-, A-, E-, and ME-optimality), which yields four optimized input profiles. For each profile, we then compute the resulting FIM and evaluate all four criteria on that same FIM. Therefore, each row is a single optimized experiment, each column is a metric used to score that experiment, and the table is a cross-evaluation matrix. All entries are reported in $\log_{10}$ scale. Bold values indicate the best entry in each column according to the objective sense: maximum for D- and E-optimality, minimum for A- and ME-optimality.

\FloatBarrier

\begin{table}[!ht]
\centering
\caption{Cross-evaluation of optimal experimental designs (rows) using all four criteria (columns). For instance, in row D-opt, the entry under the E-opt column is the $\log_{10}$ E-optimality value of the FIM generated by the D-optimal profile. Thus, diagonal values represent self-performance and off-diagonal values represent tradeoffs across criteria. Bold values denote the best value in each column using the correct objective sense (maximize D- and E-optimality; minimize A- and ME-optimality). All values are reported in $\log_{10}$ form (e.g., D-optimality is $\log_{10}\left(\left|\boldsymbol{M}\right|\right)$). *A-optimality required a scaling factor ($10^6$) to converge to a profile different from the initial guess.}
\begin{tabular}{|c|c|c|c|c|}
    \hline
    & \bfseries D-opt & \bfseries A-opt & \bfseries E-opt & \bfseries ME-opt \\ \hline \hline 
    \bfseries D-opt   & \bfseries 20.28 & -3.34 & 3.36 & 3.00    \\ \hline
    \bfseries A-opt*  & 19.74 & -3.33 & 3.35 & 2.97 \\ \hline
    \bfseries E-opt   & 19.81 & \bfseries -3.36 & \bfseries 3.3834 & 2.94     \\ \hline
    \bfseries ME-opt & 19.52 & -3.35 & 3.3831 & \bfseries 2.91    \\ \hline
\end{tabular} 
\label{tab:optimal-tclab-criteria}
\end{table}
\FloatBarrier

From the D-opt column, the D-optimal design is best (20.28), indicating the largest volumetric information gain in this set. From the A-opt column, the E-optimal design is best (-3.36), indicating better improvement of low-information directions than the converged A-optimal profile in this run. From the E-opt column, the E-optimal design is best (3.3834), which confirms the expected focus on the weakest information direction. From the ME-opt column, the ME-optimal design is best (2.91), indicating the most balanced information across parameter directions.

Table~\ref{tab:optimal-tclab-criteria} therefore shows that the MBDoE formulation is able to find its respective best value experiment for D-, E-, and ME-optimality. However, the value of A-optimality is better when utilizing the E-optimal solution than utilizing the A-optimal one. This may be explained because A-optimality required a scaling factor (in this case $10^6$) to converge a profile other than the initial guess, as shown in Figure~\ref{fig:tclab-results}. This is because the solver tolerance used is $10^{-5}$ and the eigenvalues of the covariance matrix had values ranging from approximately $10^{-4}$ to $10^{-7}$. Our choice of $10^6$ for the scaling factor scaled the eigenvalues of the covariance matrix closer to that of the model equations. Therefore, the solution depends on this scaling factor. A reasonable scaling factor for A-optimality in this equation-oriented format is the magnitude of the largest eigenvalue of the prior FIM (in this case on the order of $10^6$).

Another observation for this system is that most of the solutions render a very similar E-optimality value, changing only about 0.03 orders of magnitude over the 4 optimal designs, indicating that E-optimality (and also A-optimality) are not the best criteria for experimental design in this case. If balance in parametric certainty is desired, the ME-optimal solution is recommended, and if overall volumetric reduction in uncertainty is desired, D-optimality should be used. Also, the $\text{log}_{10}$ representation in Table~\ref{tab:optimal-tclab-criteria} makes it abundantly clear that A-optimality is dominated by the minimum eigenvalue, as most of the E-optimality values are almost exactly the inverse (or negative in the log sense) of A-optimality.

Likely, the best solutions from a D-optimal standpoint originate from exploring the highest temperatures that the system allows. The information scales with the value of the measurements and the error associated with these measurements (Eqs.~\ref{eq:cov-def} and ~\ref{eq:cov-fim-relationship}). In this case, the measurement error is assumed to be constant, leading to the higher information for higher values of measurements. Nevertheless, utilizing a scheme that uses 0\% to 100\%, ``on-off'' step functions with varying frequency combined with achieving the highest temperature in the system should lead to near-optimal information increase. The decrease in uncertainty in the unknown parameters is shown with the reduced pairwise covariance of the original parameters, as shown in Figure~\ref{fig:tclab-pairwise-uncertainty-combined}b). Here, we show the uncertainty of the original parameters using prior data (as shown in Figure~\ref{fig:tclab-pairwise-uncertainty-combined} a) compared to the uncertainty when the new data are included from the D-optimal experiment. As shown, there is still plenty of room for improvement, but the uncertainty has been greatly reduced, and the parameters have significantly fewer non-physical values within the region.

In summary, the four criteria provide distinct and interpretable tradeoffs for this system. D-optimality yields the strongest volumetric information gain, ME-optimality yields the best balance in parametric certainty, and E-optimality most directly improves the weakest information direction. For this case, A-optimality is sensitive to scaling and does not provide a clear practical advantage over E-optimality.

One can also gain a deeper understanding of what the experiments are numerically targeting utilizing the eigendecomposition of the FIM. Look into our previous work on the subject \cite{lynch2024crystallization-mbdoe}, or refer to the SI Section~\ref{sec:si-tclab-eigenanalysis} for analysis on the TCLab. The issues shown in this analysis emphasize why this system can be complicated to teach, as using standard control profiles (e.g., simple sine/step test) result in estimability issues. E-optimality shows that we can resolve some of that, but we need either additional data sources (analyzing structural identifiability) or reformulating of the model, which will be explored in future work. Additionally, future work may explore the use of global optimizers as these problems are inherently non-convex and the solutions provided here are only locally optimal. Finally, future plans include the analysis of the impact of discretization schemes on the optima to ensure the optimizer is not exploiting numerical structure (from a discretization or integration scheme) instead of mathematical structure (from the model itself).

\subsection{Characterizing a Three-Stage Membrane System for Critical Minerals and Materials Recovery} \label{sec:res-membrane}
The third case study is a process-scale model of membrane-based separations for the recovery of critical minerals and materials (CMM). Numerous studies including Kim et al. \cite{kim2021role}, Gaustad et al. \cite{gaustad2021rare}, Hamzat et al. \cite{hamzat2025rare}, and Liang et al. \cite{liang2024review} have identified secondary sources such as recycled materials (e.g., electronic waste) as promising pathways to increase domestic CMM production and meet the growing demand associated with digital technologies that rely heavily on stable supply of CMMs, such as advanced data-computing systems and semiconductor manufacturing. Additionally, Lair et al. \cite{lair2024critical}, Gebreslassie et al. \cite{gebreslassie2024advanced}, and Alemrajabi et al. \cite{alemrajabi2022separation} have identified membrane-based processes as promising technologies for the low-cost, environmentally sustainable, and energy-efficient recovery of CMMs.

\begin{figure}[h!]
    \centering
    \includegraphics[width=0.8\textwidth]{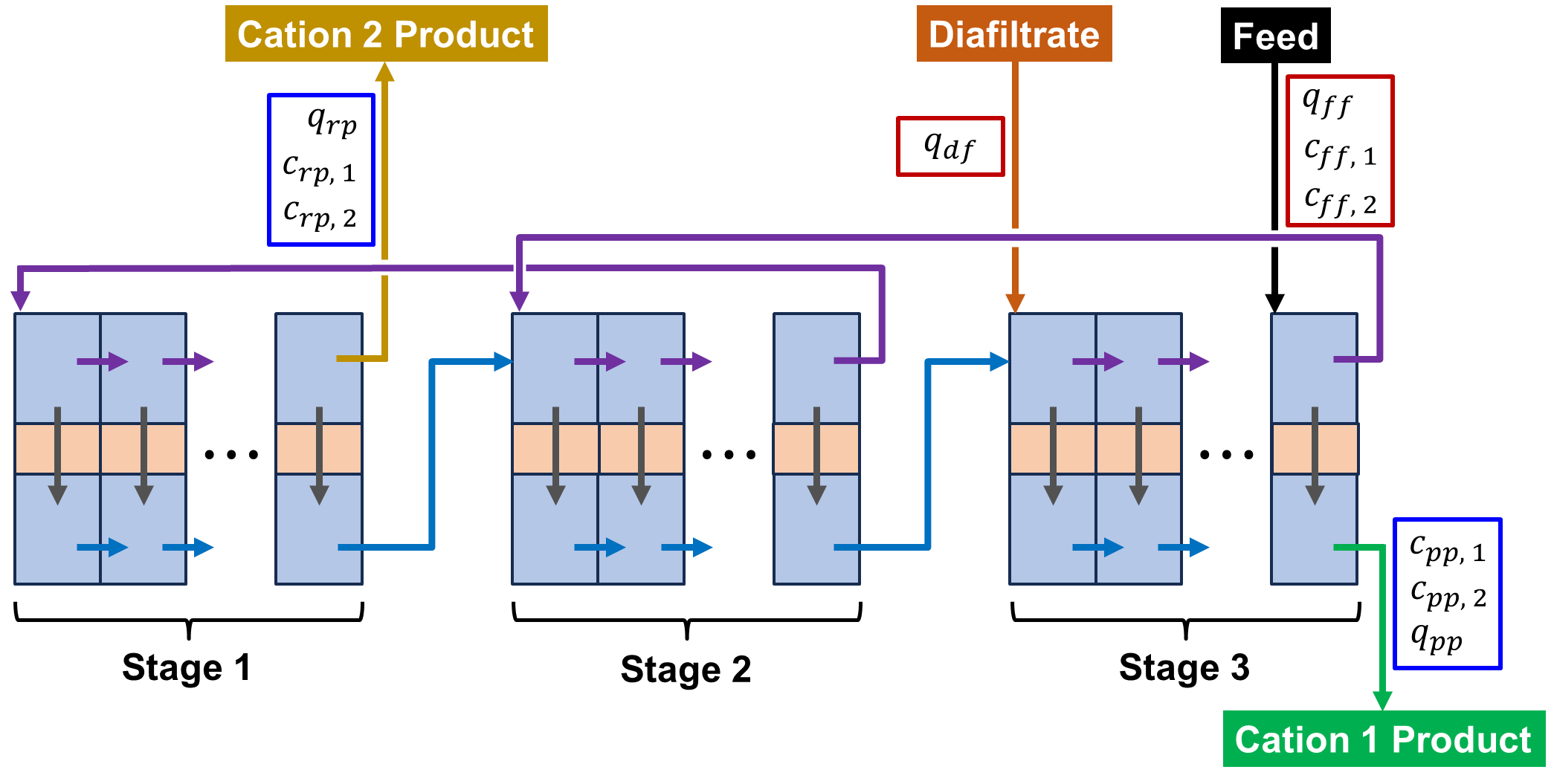}
    \caption{Proposed three-stage membrane cascade to recover critical minerals and materials from recycled sources. The text highlighted in red are the experimental design decisions, whereas those in blue are the variables that are measured in the system ($q$ denotes flowrate while $c$ denotes concentration).}
    \label{fig:membrane-cascade}
\end{figure}

Wamble et al. \cite{wamble2022optimal} proposed a three-stage diafiltration membrane cascade (see Figure~\ref{fig:membrane-cascade}) for the efficient recovery of CMMs from secondary sources (e.g., recycled materials such as electronic waste). This membrane system uses a diafiltrate, which is a dilute solution, to reduce the effects of concentration polarization (a condition caused by the accumulation of ions on the feed side of the membrane) which often leads to a decrease in the separation performance \cite{wamble2022optimal, muetzel2022device, ouimet2022data}. The separation efficiency of the membrane system is defined by parameters such as the water flux and the sieving coefficients. The water flux measures the flow of water through the membrane, whereas the sieving coefficients measure the ability of the membrane to allow or reject the passage of ions. Wamble et al. \cite{wamble2022optimal} assumed that the water flux and sieving coefficients are constant across all elements and stages of the membrane, which may not be a realistic assumption for many CMM separation applications.

Instead, in this study, we assume that the water flux is dependent on the applied pressure in the membrane system (Eq.~\ref{eq:basic-water-flux}) and the sieving coefficients vary linearly with the ionic strength of the feed solution (Eq.~\ref{eq:basic-sieving-coefficient}):
\begin{align}
     J_{w} = L_p\left(\Delta P - \Delta \pi\right) \label{eq:basic-water-flux}
\end{align}

\noindent with
\begin{align}
     \Delta \pi = RT\sum_{k=1}^{n}\left(c_{r,\, k} - c_{p,\, k}\right) \label{eq:basic-osmotic-pressure}
\end{align}

\noindent and
\begin{align}
    S_{i} = \bar{S}_{i} + \frac{\delta_{i}}{2}\sum_{k=1}^{n}c_{f,\, k,}\, z_{k}^2 
    \quad \forall \, i \in \{1,\,2,\cdots,\,n \} \label{eq:basic-sieving-coefficient}
\end{align}

\noindent where $J_{w}$ is the flux of water through the membrane, $\Delta P$ is the applied pressure, $\Delta \pi$ is the opposing osmotic pressure, $n$ is the number of ionic species in the solution, $S_{i}$ is the sieving coefficient of ionic species $i$, $z_{k}$ is valency of ionic species $k$, $R$ is the gas constant, $T$ is the temperature of the solution, $c_{f,\, k}$, $c_{r,\, k}$, and $c_{p,\, k}$ are the concentrations of ionic species $k$ on the feed, retentate, and permeate sides of the membrane elements, respectively.

The pressure-dependent water flux and concentration-dependent sieving coefficient models are assumed based on the findings of Baker and Wijmans \cite{wijmans1995solution}, Baker \cite{baker2023membrane}, Bartels et al. \cite{bartels2005effect}, and Bason et al. \cite{bason2010analysis} on solution diffusion and ion transport in membranes. Five parameters in the water flux and sieving coefficient models are unknown or are not directly measurable. These unknown model parameters are listed in Table~\ref{tab:unknown_parameters}.
\begin{table}[h!]
\centering
\caption{Unknown model parameters in the three-stage membrane cascade.}
\begin{tabular}{lll}
\hline
\textbf{Parameter} & \textbf{Unit} & \textbf{Description} \\
\hline
$L_p$ & $\text{m}\cdot\text{h}^{-1}\cdot\text{Pa}^{-1}$ & Water permeability constant \\
$\bar{S}_{1}$ & dimensionless & Constant cation 1 sieving coefficient \\
$\bar{S}_{2}$ & dimensionless & Constant cation 2 sieving coefficient \\
$\delta_{1}$ & $\text{m}^3\cdot\text{mol}^{-1}$ & Cation 1 ionic strength coefficient \\
$\delta_{2}$ & $\text{m}^3\cdot\text{mol}^{-1}$ & Cation 2 ionic strength coefficient \\
\hline
\end{tabular}
\label{tab:unknown_parameters}
\end{table}

The complete mathematical model of the three-stage membrane cascade showing detailed solvent and species flows with 444 variables and 444 equality constraints is presented in Section~\ref{sec:detailed-membrane-model}.

Based on these solvent flows, the experimental design decisions include the fresh feed flowrate, diafiltrate flowrate, and the concentration of ionic species in the fresh feed (see Table~\ref{tab:measured_variables}). The experimental measurements consist of the flowrates of the permeate and retentate products and their corresponding ionic species concentrations, as summarized in Table~\ref{tab:measured_variables}. The measurement error associated with the flowrate and concentration measurements is assumed to be 2 $\text{m}^{3}\cdot\text{h}^{-1}$ and 0.1 $\text{kg}\cdot\text{m}^{-3}$, respectively.
\begin{table}[h!]
\centering
\caption{Measurements and design decisions in the three-stage membrane cascade.}
\small
\begin{tabular}{lll}
\toprule
\textbf{Variable} & \textbf{Unit} & \textbf{Description} \\
\midrule
\multicolumn{3}{c}{\textbf{Measurements}} \\
\midrule
$q_{pp}$ & $\text{m}^{3}\cdot\text{h}^{-1}$ & Flowrate of cation 1 product \\
$c_{pp,\,1}$ & $\text{kg}\cdot\text{m}^{-3}$ & Concentration of cation 1 in cation 1 product \\
$c_{pp,\, 2}$ & $\text{kg}\cdot\text{m}^{-3}$ & Concentration of cation 2 in cation 1 product \\
$q_{rp}$ & $\text{m}^{3}\cdot\text{h}^{-1}$ & Flowrate of cation 2 product \\
$c_{rp,\, 1}$ & $\text{kg}\cdot\text{m}^{-3}$ & Concentration of cation 1 in cation 2 product \\
$c_{rp,\, 2}$ & $\text{kg}\cdot\text{m}^{-3}$ & Concentration of cation 2 in cation 2 product \\
\midrule
\multicolumn{3}{c}{\textbf{Design Decisions}} \\
\midrule
$q_{df}$ & $\text{m}^{3}\cdot\text{h}^{-1}$ & Flowrate of the fresh diafiltrate \\
$q_{ff}$ & $\text{m}^{3}\cdot\text{h}^{-1}$ & Flowrate of the fresh feed \\
$c_{ff,\,1}$ & $\text{kg}\cdot\text{m}^{-3}$ & Concentration of cation 1 in the fresh feed \\
$c_{ff,\,2}$ & $\text{kg}\cdot\text{m}^{-3}$ & Concentration of cation 2 in the fresh feed \\
\bottomrule
\end{tabular}
\label{tab:measured_variables}
\end{table}

Lastly, this study considers the binary separation of cation 1 ($\text{Li}^+$) and cation 2 ($\text{Co}^{2+}$) as the only ions present in the feed, although other ions can be treated. To fully characterize the membrane cascade and enable scale-up analyses (e.g., techno-economic analysis), the five unknown model parameters listed in Table~\ref{tab:unknown_parameters} must be estimated from experimental data.

As with the previous case studies, starting with some experimental data is required to provide prior information to the sequential experimental design process. In this case, the membrane system's prior data is obtained from synthetic experiments generated by performing a 2-level factorial design with low and high levels of fresh feed and diafiltrate flowrates in the ranges [90, 110] and [27, 33] $\text{m}^{3}\cdot \text{h}^{-1}$, respectively, given fresh feed cation 1 and cation 2 concentrations of 1.7 and 17 $\text{kg}\cdot \text{m}^{-3}$, respectively \cite{wamble2022optimal}. The ground truth values of the parameters used in the full-factorial simulations are shown in Table~\ref{tab:true-parameters}. The simulation results were corrupted with uncorrelated Gaussian measurement errors with zero mean and covariance matrix, $\boldsymbol{\Sigma_{y}}$, whose leading diagonal contains the square of the errors (2 $\text{m}^{3}\cdot\text{h}^{-1}$ and 0.1 $\text{kg}\cdot\text{m}^{-3}$ for flowrate and concentration measurements, respectively).

The first question is the following: \emph{based on the full-factorial data, can we reliably estimate the five model parameters to fully characterize the membrane system}?
\begin{table}[h!]
\centering
\caption{Ground truth values of the model parameters used in the full-factorial experiments.}
\begin{tabular}{lll}
\hline
\textbf{Parameter} & \textbf{Unit} & \textbf{True Value} \\
\hline
$L_p$ & $\text{m}\cdot\text{h}^{-1}\cdot\text{Pa}^{-1}$ & $3\times 10^{-7}$ \\
$\bar{S}_{1}$ & dimensionless & 1.3 \\
$\bar{S}_{2}$ & dimensionless & 0.5 \\
$\delta_{1}$ & $\text{m}^3\cdot\text{mol}^{-1}$ & $5\times 10^{-4}$ \\
$\delta_{2}$ & $\text{m}^3\cdot\text{mol}^{-1}$ & $1.5\times 10^{-4}$ \\
\hline
\end{tabular}
\label{tab:true-parameters}
\end{table}

To answer this question, we once again utilize $\texttt{parmest}$. Similar to the other cases, the model is labeled following the style shown in Figure~\ref{fig:general-experiment-code-snippet}. In this case, the experimental measurements of flowrates and concentrations have unit and error mismatch, making maximum likelihood estimation even more important as the weighted sum of squared errors (WSSE) is required to reconcile these differences. The $\texttt{parmest}$ toolbox automatically generates the WSSE objective function, as defined in Eq.~\ref{eq:param-estimation}, enabling accurate estimation of model parameters from these non-uniform measurements and solving the corresponding minimization problem. A detailed derivation of Eq.~\ref{eq:param-estimation} is provided in Section~\ref{sec:si-parmest}.

\begin{figure}[H]
\includegraphics[width = \textwidth]{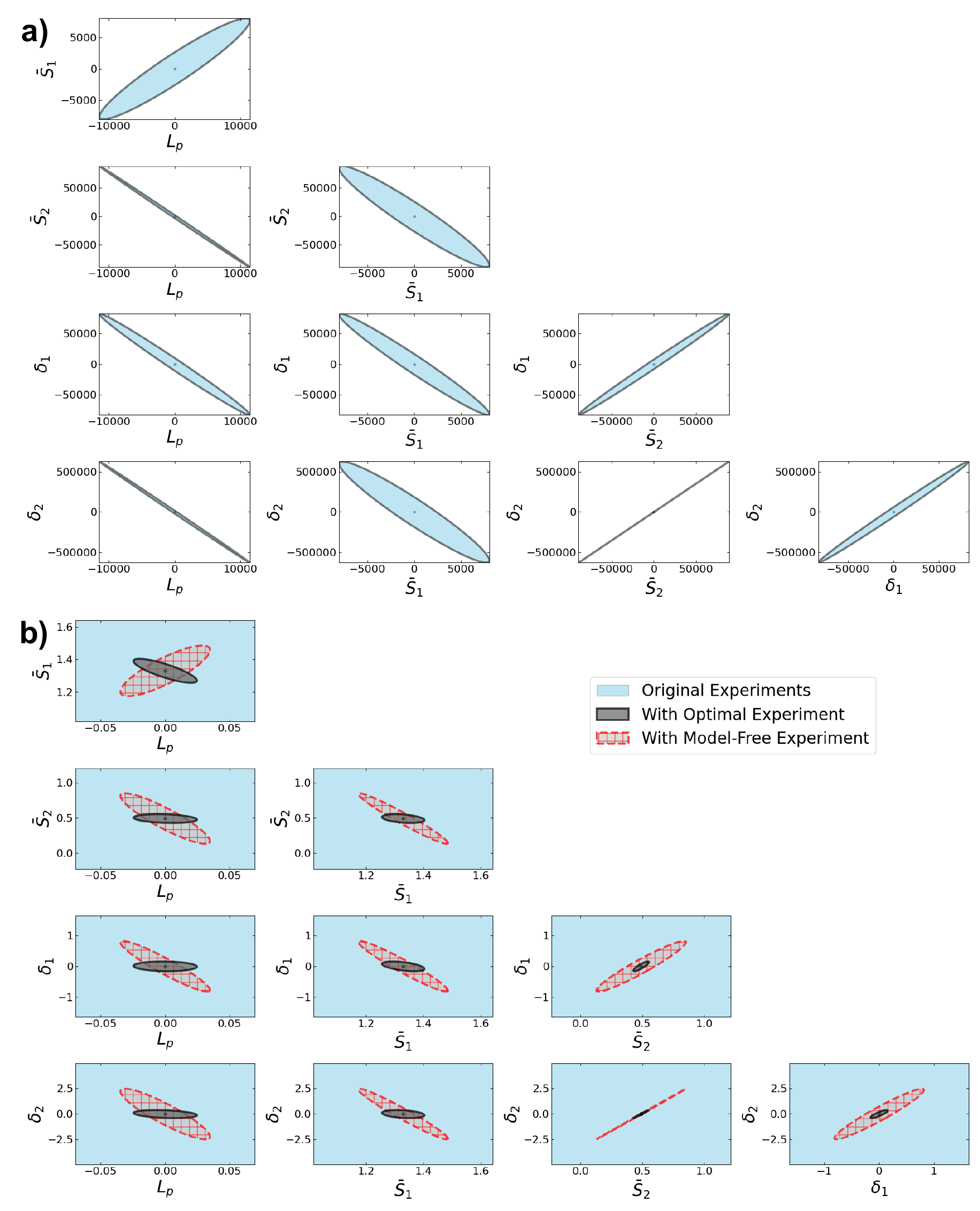}
\caption{Pairwise covariance of parameters with preliminary data only (a) and additional data (b) in the membrane case study. With additional data (b), uncertainty after the optimal experiment (gray shading with black border) and adding a model-free experiment (dashed red border) significantly reduces uncertainty (4 to 6 orders of magnitude).} \label{fig:membrane-uncertainty-parmest-combined}
\end{figure}

Using the preliminary data described in the previous section, the estimated values of the model parameters from $\texttt{parmest}$ are shown in Table~\ref{tab:optimal-parmest-values}. The precision of the parameter estimates was quantified by computing the covariance matrix, which is approximated as the inverse of the FIM, as defined in Eq.~\ref{eq:FIM-wsse}, assuming a consistent and asymptotically normal estimator. The parameter estimates are close to the true values in Table~\ref{tab:true-parameters}, but their associated uncertainty is high, as shown in Figure~\ref{fig:membrane-uncertainty-parmest-combined}a). This leads us to the next question: \emph{which experiment(s) should we perform to provide the most information (reduce uncertainty) about the parameters}?
\begin{table}[h!]
\centering
\caption{Estimated values of the model parameters where the uncertainty from the covariance matrix is shown using pair-wise posteriors in Figure~\ref{fig:membrane-uncertainty-parmest-combined}a. $\text{\textbf{Estimate 1}}$ are the parameter estimates using data from full-factorial simulation, whereas $\text{\textbf{Estimate 2}}$ are the values from data that includes the optimal experimental conditions in Table~\ref{tab:optimal-membrane-conditions}.}
\begin{adjustbox}{max width=\textwidth}
\begin{tabular}{lllll}
\hline
\textbf{Parameter} & \textbf{Unit} & \textbf{True Value} & \textbf{Estimate 1} & \textbf{Estimate 2} \\
\hline
$L_p$ & $\text{m}\cdot\text{hr}^{-1}\cdot\text{Pa}^{-1}$ & $3\times 10^{-7}$ & $2.97\times 10^{-7}$ $\pm$ $3.24\times 10^{-9}$ & $2.98\times 10^{-7}$ $\pm$ $2.90\times 10^{-9}$ \\
$\bar{S}_{1}$ & dimensionless & 1.3 & 1.33 $\pm$ 0.12 & 1.33 $\pm$ 0.07 \\
$\bar{S}_{2}$ & dimensionless & 0.5 & 0.49 $\pm$ 0.02 & 0.49 $\pm$ 0.01 \\
$\delta_{1}$ & $\text{m}^3\cdot\text{mol}^{-1}$ & $5\times 10^{-4}$ & $5.14\times 10^{-4}$ $\pm$ $1.58\times 10^{-4}$ & $5.01\times 10^{-4}$ $\pm$ $6.43\times 10^{-5}$ \\
$\delta_{2}$ & $\text{m}^3\cdot\text{mol}^{-1}$ & $1.5\times 10^{-4}$ & $1.34\times 10^{-4}$ $\pm$ $3.30\times 10^{-5}$ & $1.35\times 10^{-4}$ $\pm$ $2.08\times 10^{-5}$ \\
\hline
\end{tabular}
\end{adjustbox}
\label{tab:optimal-parmest-values}
\end{table}
\FloatBarrier

\FloatBarrier

Therefore, we wish to design an experiment to improve the uncertainty in the parameter estimates, shown in Figure~\ref{fig:membrane-uncertainty-parmest-combined}a). The experimental design decisions for this case are the feed flowrate, diafiltrate flowrate, and the concentrations of ionic species in the feed. In Figure~\ref{fig:membrane-uncertainty-parmest-combined}a) we see many directions of high uncertainty, especially in $\delta_{2}$, so we will utilize E-optimality to determine the next best experiment. The optimal design conditions identified by E-optimality are specified in Table~\ref{tab:optimal-membrane-conditions}. Here, the recommendation is to perform a new experiment with higher concentrations of cation 1 and cation 2 (cation 1 in [1.5, 2.0] and cation 2 in [15, 20] $\text{kg}\cdot\text{m}^{-3}$) and remain at the highest feed flowrate allowed and lowest diafiltrate flowrate defined by the bounds specified previously (feed and diafiltrate flowrates in the ranges [90, 110] and [27, 33] $\text{m}^{3}\cdot \text{h}^{-1}$, respectively).

High feed flowrate and ionic species concentration will produce higher concentrations of the ions in the permeate and retentate products as a result of the increased mass flow of the ions across the membrane stages. For the diafiltrate flowrate, the lowest condition (27 $\text{m}^{3}\cdot \text{h}^{-1}$) was selected as optimal because, compared to other higher conditions, the lowest flowrate will produce the smallest increase in the volume of the solvent on the feed side of the membrane stages, leading to higher concentration of ionic species in the permeate and retentate products. Additionally, the criteria that focus on eigenvalues (A-opt, E-opt, ME-opt) all recommend these experimental conditions.

\begin{table}[h!]
\centering
\caption{E-optimal experimental design for the membrane cascade system.}
\begin{tabular}{lll}
\hline
\textbf{Design Decision} & \textbf{Unit} & \textbf{Optimal Value} \\
\hline
$q_{df}$ & $\text{m}^{3}\cdot\text{h}^{-1}$ & 27 \\
$q_{ff}$ & $\text{m}^{3}\cdot\text{h}^{-1}$ & 110 \\
$c_{ff,\,1}$ & $\text{kg}\cdot\text{m}^{-3}$ & 2.0 \\
$c_{ff,\,2}$ & $\text{kg}\cdot\text{m}^{-3}$ & 20.0 \\
\hline
\end{tabular}
\label{tab:optimal-membrane-conditions}
\end{table}
\FloatBarrier

Including data from this experiment at the optimal parameter estimates from Table~\ref{tab:optimal-parmest-values} drastically improves confidence in the model predictions, as shown by the stark reduction in uncertainty in Figure~\ref{fig:membrane-uncertainty-parmest-combined}b). The parameters that have very small values (i.e., $L_p$, $\delta_{1}$, and $\delta_{2}$) still contain regions of non-physical parameter values in the uncertainty set, but are significantly better than those included in the original uncertainty regions in Figure~\ref{fig:membrane-uncertainty-parmest-combined}a). The uncertainty regions from Figure~\ref{fig:membrane-uncertainty-parmest-combined}a) are included in Figure~\ref{fig:membrane-uncertainty-parmest-combined}b) (blue shaded ellipses). However, these are so large in comparison to the improvement that they encompass the entire axes shown here. Quantitatively, the axes limits change by 4 to 6 orders of magnitude from Figure~\ref{fig:membrane-uncertainty-parmest-combined}a) to Figure~\ref{fig:membrane-uncertainty-parmest-combined}b). Additionally, Figure~\ref{fig:membrane-uncertainty-parmest-combined}b) shows a comparison of the reduction in uncertainty due to the optimal experiment (dark gray shading with black border) with the reduction using a reasonable next best experiment based on expert or model-free recommendations (light gray shading with dashed red border). Here, a reasonable addition to the experiment set would be the middle of both flow ranges, at a feed flowrate of 100 $\text{kg}\cdot\text{m}^{-3}$ and a diafiltrate flowrate of 30 $\text{kg}\cdot\text{m}^{-3}$ (center of the full 2-factor design). Although this experiment substantially reduces uncertainty relative to the prior data set, the reduction is substantially smaller than that achieved by using the E-optimal experiment. This is particularly noticeable when looking at the worst pair-wise uncertain regions (i.e., $\bar{S}_{2}$ vs. $\delta_{2}$), which by design is the focus of the E-optimal experiment. Nevertheless, additional experiments should be conducted alongside this E-optimal experiment to reduce uncertainty to an acceptable level (as indicated in Figure~\ref{fig:workflow}), enabling the use of this model for the techno-economic evaluation of new CMM recovery pathways.\cite{lair2024critical} Also, other uncertainty representations beyond the consideration of the FIM may be more appropriate to retain physical boundaries of the parameters alongside accurate parametric uncertainty.

\section{Conclusions and Future Directions}
This paper introduces a methodology for embedding complicated design criteria within equation-oriented solvers to make optimal experimental design more robust by offering a variety of different objective options (design criteria). To the authors' knowledge, embedding these eigenvalue calculations with more than four unknown parameters in these equation-oriented programs has not yet been demonstrated. Also, the literature is sparse on the derivatives of the condition number of a symmetric real matrix. Visually, we illustrated, using a single-variable problem, how the FIM design criteria are related and target different regions of the experimental design space. E-optimality, which was only made possible with the advancements herein, was used to target highly uncertain parameters and markedly improved the theoretical information of the experimental data set. In general, we improved equation-oriented optimal experimental design by incorporating eigenvalue-based metrics, enabling scientists to focus on information that addresses specific deficiencies in model confidence.

However, these methods are only as robust as the numerical routines that solve the eigenvalue problem (e.g., \texttt{NumPy.linalg.eig}) and the inverse of a matrix (e.g., \texttt{NumPy.linalg.pinv}). The robustness of E-optimality when eigenvalues are repeated must be considered, even though this is numerically rare. We hypothesize that a simple fix for repeated eigenvalues is to add a small, random diagonal matrix (similar to an inertia correction in numerical optimization) to prevent repetition; this approach may be explored in future work. Additionally, the design metrics presented within this paper, A-, D-, E-, and ME-optimality, do not cover the entire range of potential objective functions that could be used in optimal experimental design. This work outlines a platform for which more complicated objectives can be formulated and utilized within the same `Grey Box' structure (\texttt{ExternalGreyBox} feature within Pyomo). Therefore, there is a great opportunity to expand the Pyomo.DoE framework (and equation-oriented optimal experimental design in general) to include more exotic design criteria that would also be challenging, or impossible, to represent as explicit algebraic equations. Also, we plan to make the Pyomo.DoE framework more compact by including options to forgo the finite difference representation of the sensitivity matrix, $\boldsymbol{Q}$, by analytically computing sensitivity equations using symbolic differentiation

Optimal experimental design is a mature field; however, there remain opportunities to contribute, particularly in the context of complex first-principles models and simultaneous, equation-oriented programming. For example, we see extensive opportunities for Pyomo.DoE to accelerate the creation and validation of science-based models for complex, many-component mixtures important for CMM processing and other emerging application areas \cite{lair2024critical}.

\section*{Symbols and Nomenclature}
\begin{center}
\begin{longtable}{|c|p{8.25cm}|}
    \caption{Mathematical symbols used in models and functions described in Section~\ref{sec:parmest-and-fim}\/} \\

    \hline \bfseries Symbol & \bfseries Definition \\ 
    \hline \hline 
    \endfirsthead 

    \hline \bfseries Symbol & \bfseries Definition \\ 
    \hline \hline 
    \endhead 

    \hline \hline 
    \endfoot 

    \hline \hline 
    \endlastfoot 

    \multicolumn{2}{|l|}{\textbf{Sets}} \\ \hline 
    $\boldsymbol{\Phi} \in \mathbb{R}^{N_d}$ & Set of experimental design conditions \\ 
    $\boldsymbol{\Theta} \in \mathbb{R}^{p}$ & Set of feasible unknown parameter values\\ \hline 

    \multicolumn{2}{|l|}{\textbf{Functions}} \\ \hline 
    $\boldsymbol{f}$: $\mathbb{R}^{N_\text{states} + N_d + p} \rightarrow \mathbb{R}^{N_\text{meas}}$ & Mathematical model for predicting $\boldsymbol{y}$ \\
    $\ell: \mathbb{R}^{N_\text{meas} + p} \rightarrow \mathbb{R}$ & Log likelihood function \\ 
    $\Psi: \mathbb{R}^{p \times p} \rightarrow \mathbb{R}$ & Experimental design criterion \\ 
    $\text{trace}: \mathbb{R}^{p \times p} \rightarrow \mathbb{R}$ & Trace of a matrix \\ 
    $\text{det}: \mathbb{R}^{p \times p} \rightarrow \mathbb{R}$ & Determinant of a matrix \\ \hline 

    \multicolumn{2}{|l|}{\textbf{Constants}} \\ \hline 
    $\boldsymbol{e}_i \in \mathbb{R}^{p}$ & Unit vector in the direction of component $i$\\
    $N_\text{exp} \in \mathbb{Z}^{+}$ &  Number of experiments \\ 
    $N_\text{meas} \in \mathbb{Z}^{+}$ & Number of measurements in an experiment \\ 
    $N_\text{d} \in \mathbb{Z}^{+}$ & Number of experimental design decisions \\ 
    $p \in \mathbb{Z}^{+}$ & Number of unknown parameters in the model \\ \hline 

    \multicolumn{2}{|l|}{\textbf{Variables}} \\ \hline
    $\boldsymbol{y}_i \in \mathbb{R}^{N_\text{meas}}$ & Experimental measurements (data) for experiment $i$ \\
    $\hat{\boldsymbol{y}}_i \in \mathbb{R}^{N_\text{meas}}$ & Prediction of measurements for experiment $i$ using model $\boldsymbol{f}$ \\
    $\boldsymbol{\phi} \in \mathbb{R}^{N_\text{d}}$ & Experimental design decisions $i$ \\
    $\hat{\boldsymbol{\phi}} \in \mathbb{R}^{N_\text{d}}$ & Optimal experimental design decisions $i$ \\
    $\boldsymbol{x} \in \mathbb{R}^{N_\text{states}}$ & State variables in model $\boldsymbol{f}$ \\
    $\boldsymbol{\theta} \in \mathbb{R}^{p}$ & Unknown parameters in model $\boldsymbol{f}$ \\ 
    $\hat{\boldsymbol{\theta}} \in \mathbb{R}^{p}$ & Optimal value of unknown parameters in model $\boldsymbol{f}$ \\  
    $\boldsymbol{\varepsilon}_i \in \mathbb{R}^{N_\text{meas}}$ & Measurement error for experiment $i$ \\
    $\boldsymbol{\Sigma}_{\boldsymbol{y}} \in \mathbb{R}^{N_\text{meas} \times N_\text{meas}}$ & Measurement covariance matrix for an experiment \\
    $\boldsymbol{Q} \in \mathbb{R}^{N_\text{meas} \times p}$ & Sensitivity matrix for outputs $\hat{\boldsymbol{y}}$ with respect to unknown parameters $\boldsymbol{\theta}$\\ 
    $\boldsymbol{q}_i \in \mathbb{R}^{N_\text{meas}}$ & Sensitivity vector for outputs $\hat{\boldsymbol{y}}$ with respect to unknown parameter $i$\\ 
    $h_i \in \mathbb{R}$ & Finite difference perturbation to compute $\boldsymbol{q}_i$ for unknown parameter $i$\\ 
    $\boldsymbol{V} \in \mathbb{R}^{p \times p}$ & Parameter covariance matrix \\
    $\boldsymbol{V}_{\boldsymbol{\theta}} \in \mathbb{R}^{p \times p}$ & Prior parameter covariance matrix at $\boldsymbol{\theta}$ \\
    $\boldsymbol{M} \in \mathbb{R}^{p \times p}$ & Fisher Information Matrix \\
    $\boldsymbol{M}_{\boldsymbol{\theta}} \in \mathbb{R}^{p \times p}$ & Prior Fisher Information Matrix at $\boldsymbol{\theta}$ \\ 
    $\boldsymbol{M}_{\text{GB}} \in \mathbb{R}^{p \times p}$ & Grey Box Fisher Information Matrix \\
    $M_{ij} \in \mathbb{R}$ & The $i,j$-th element of the Fisher Information Matrix \\ 
    $\boldsymbol{\lambda} \in \mathbb{R}^p$ & Vector of eigenvalues for the FIM \\
    $\lambda_{\text{min}} \in \mathbb{R}$ & Minimum eigenvalue of the FIM \\
    $\lambda_{\text{max}} \in \mathbb{R}$ & Maximum eigenvalue of the FIM \\ 
    $\boldsymbol{v}_s \in \mathbb{R}^p$ & Eigenvector of the FIM corresponding to the $s$-th eigenvalue $\lambda_s$ \\ 
    $\boldsymbol{v}_{\text{min}} \in \mathbb{R}^p$ & Eigenvector of the FIM corresponding to the maximum eigenvalue \\ 
    $\boldsymbol{v}_{\text{min}} \in \mathbb{R}^p$ & Eigenvector of the FIM corresponding the minimum eigenvalue $\lambda_s$ \\ 
    $\kappa \in \mathbb{R}$ & Condition number of the FIM \\ 
    $\boldsymbol{a} \in \mathbb{R}^p$ & Arbitrary real vector to illustrate ellipsoidal interpretation of the condition number of the FIM \\ \hline

\end{longtable}
\end{center}

\section{Author Contributions}
\textbf{Daniel J. Laky:} Conceptualization, Methodology, Software, Formal analysis, Investigation, Data curation, Writing– original draft, Writing– review \& editing, Visualization.
\textbf{Shammah Lilonfe:} Formal analysis, Data curation, Writing– original draft, Software, Writing– review \& editing, Visualization.
\textbf{Shawn Martin:} Conceptualization, Methodology, Software -- redesign of \texttt{parmest} using the \texttt{Experiment} class, Writing- review \& editing.
\textbf{Katherine A. Klise:} Conceptualization, Methodology, Software Writing- review \& editing.
\textbf{Alexander W. Dowling:} Project administration, Conceptualization, Methodology, Software, Formal analysis, Investigation, Data curation, Writing– original draft, Writing– review \& editing, Visualization.

\section{Data Availability}
All results were run on a Mac mini using an Apple\copyright M4 chip and 32GB RAM. Pyomo.DoE and all features described herein are currently available in the main branch of Pyomo, version \texttt{6.9.4} at \url{https://github.com/Pyomo/pyomo}. More specifically, the code for the Grey Box implementation can be found at \url{https://github.com/Pyomo/pyomo/blob/main/pyomo/contrib/doe/grey\_box\_utilities.py}. In addition, all code used to generate the results, the resulting figures, and how to utilize this tool are included at \url{https://github.com/dowlinglab/doe-greybox-paper}. The results using GreyBox require the installation of \texttt{cyipopt} with HSL linear solvers. Guidance for installation and usage of these solvers is included in the \texttt{cyipopt} documentation \url{https://cyipopt.readthedocs.io/en/stable/}. Other packages include \texttt{SciPy}~\cite{2020SciPy-NMeth} version \texttt{1.15.3}, \texttt{NumPy}~\cite{numpy} version \texttt{2.2.6}, \texttt{pandas}~\cite{reback2020pandas} version \texttt{2.2.3} for data management, and \texttt{matplotlib} version \texttt{3.10.3} for plotting.

\section{Conflicts of Interest}
The authors declare no conflicts of interest.

\section{Acknowledgments}
The authors would like to acknowledge Michael Bynum for his efforts in helping identify the numerical errors discussed in Section~\ref{sec:si-symmetry}. Also, we thank Prof. Liviu Nicolaescu from Notre Dame for discussions and his help with confirming that the derivative formulas presented in this work are valid for real, symmetric matrices with full rank. In addition, we thank Prof. Bill Phillip and his group at Notre Dame for ideas on the water flux and sieving coefficients model of the membrane case study. Finally, we thank the PrOMMiS team for technical input on code developments for the implementation of the membrane model in the IDAES-PSE+ framework.

This effort was funded by the U.S. Department of Energy’s Process Optimization and Modeling for Minerals Sustainability (PrOMMiS) Initiative, supported by the Hydrocarbons and Geothermal Energy Office. For questions and comments, please contact our Technical Director, Aimee Curtright, \url{aimee.curtright@netl.doe.gov}.

Disclaimer: This project was funded by the Department of Energy, National Energy Technology Laboratory an agency of the United States Government, $\textbf{in part}$, through a $\textbf{site}$ support contract. Neither the United States Government nor any agency thereof, nor any of their employees, nor the support contractor, nor any of their employees, makes any warranty, express or implied, or assumes any legal liability or responsibility for the accuracy, completeness, or usefulness of any information, apparatus, product, or process disclosed, or represents that its use would not infringe privately owned rights. Reference herein to any specific commercial product, process, or service by trade name, trademark, manufacturer, or otherwise does not necessarily constitute or imply its endorsement, recommendation, or favoring by the United States Government or any agency thereof. The views and opinions of authors expressed herein do not necessarily state or reflect those of the United States Government or any agency thereof. Sandia National Laboratories is a multi-mission laboratory managed and operated by National Technology and Engineering Solutions of Sandia, LLC., a wholly owned subsidiary of Honeywell International, Inc., for the U.S. Department of Energy's National Nuclear Security Administration under contract DE-NA-0003525.

\bibliographystyle{elsarticle-num}
\bibliography{sources.bib}{}

\renewcommand\thesection{\arabic{section}}

\setcounter{equation}{0}
\renewcommand{\theequation}{S\arabic{equation}}

\setcounter{figure}{0}
\renewcommand{\thefigure}{S\arabic{figure}}

\setcounter{table}{0}
\renewcommand{\thetable}{S\arabic{table}}

\setcounter{page}{1}
\renewcommand{\thepage}{SI-\arabic{page}}

\setcounter{section}{0}
\renewcommand{\thesection}{S\arabic{section}}

\newpage

{\Large
\begin{center}
Supporting Information: Optimal Experimental Design using Eigenvalue-Based Criteria with Pyomo.DoE
\end{center}
}
\begin{center}
Daniel J. Laky${}^{1,2}$, Shammah Lilonfe${}^{1}$, Shawn B. Martin${}^{3}$, Katherine A. Klise${}^{4}$, Bethany L. Nicholson${}^{5}$, John D. Siirola${}^{5}$, Alexander Dowling${}^{1,*}$ corresponding email - adowling@nd.edu

\emph{${}^{1}$Department of Chemical and Biomolecular Engineering}\\
\emph{University of Notre Dame, Notre Dame, IN 46556, USA}\\
\emph{$^{2}$Departmant of Chemical Engineering} \\
\emph{Auburn University, Auburn, AL 36849} \\
\emph{$^{3}$ Mission Aanalytics, \\Sandia National Laboratories, Albuquerque, NM 87185}\\
\emph{$^{4}$ Energy Water Systems Integration, \\Sandia National Laboratories, Albuquerque, NM 87185}\\
\emph{$^{5}$ Center for Computing Research, \\Sandia National Laboratories, Albuquerque, NM 87185}\\
\end{center}

\section{Parameter estimation background: maximum likelihood and Fisher information matrix estimation} \label{sec:si-parmest}
Any experimental data can be expressed in the following mathematical form:
\begin{equation} \label{eq:general-model-form}
\begin{aligned}
    \boldsymbol{y}_i & = \boldsymbol{f}\left(\boldsymbol{x}_{i}, \boldsymbol{\theta}\right)\;+ \boldsymbol{\varepsilon}_i \\
    & = \hat{\boldsymbol{y}}_{i} + \boldsymbol{\varepsilon}_i, \quad \forall \; i \in \left\{1, \ldots, N_{\text{exp}}\right\}
\end{aligned}
\end{equation}

\noindent where $\boldsymbol{y}_{i} \in \mathbb{R}^n$ are observations of the measured or output variables, $\hat{\boldsymbol{y}}_{i} \in \mathbb{R}^n$ are model predictions of the measured variables, $\boldsymbol{x}_{i} \in \mathbb{R}^{q}$ are the decision or input variables, $\boldsymbol{\theta} \in \mathbb{R}^p$ are the model parameters, $\boldsymbol{\varepsilon}_{i} \in \mathbb{R}^n$ are the measurement errors, and $N_{\text{exp}}$ is the number of experiments.

If measurement errors are correlated with mean 0 and covariance matrix, $\boldsymbol{\Sigma}_{\boldsymbol{y}}$, known a priori, such that:
\begin{align}
   \boldsymbol{\varepsilon}_i &\sim N(0, \boldsymbol{\Sigma}_{\boldsymbol{y}})
\end{align}

\noindent The likelihood function of independent measurements, $\boldsymbol{y}_1,\cdots, \boldsymbol{y}_{N_{\text{exp}}}$, with a probability function, $f_{\boldsymbol{y}_{1}},\cdots, f_{\boldsymbol{y}_{N_{\text{exp}}}}$, is defined as:
\begin{align}
    L(\boldsymbol{\theta}; \boldsymbol{y}_{1},\cdots, \boldsymbol{y}_{N_{\text{exp}}}) = \prod_{i=1}^{N_{\text{exp}}} f_{\boldsymbol{y}_i}(\boldsymbol{y}_i;\boldsymbol{\theta})
\end{align}
\begin{equation} \label{eq:likelihood-derivation}
\begin{aligned}
    L(\boldsymbol{\theta}; \boldsymbol{y}_1,\cdots, \boldsymbol{y}_{N_{\text{exp}}}) & = (2\pi)^{-N_{\text{exp}}/2} |\det \boldsymbol{\Sigma}_{\boldsymbol{y}}|^{-N_{\text{exp}}/2}\\
    & \times \exp \left(-\frac{1}{2} \sum_{i \in \left\{1, \ldots, N_{\text{exp}}\right\}} \left(\boldsymbol{y}_{i} - \boldsymbol{f}(\boldsymbol{x}_{i};\boldsymbol{\theta})\right)^\text{T} \boldsymbol{\Sigma}_{\boldsymbol{y}}^{-1} \left(\boldsymbol{y}_{i} - \boldsymbol{f}(\boldsymbol{x}_{i};\boldsymbol{\theta})\right) \right)
\end{aligned}
\end{equation}

\noindent where the weighted sum of squared errors (WSSE) is:
\begin{align}
\text{WSSE} = \frac{1}{2} \sum_{i \in \left\{1, \ldots, N_{\text{exp}}\right\}} \left(\boldsymbol{y}_{i} - \boldsymbol{f}(\boldsymbol{x}_{i};\boldsymbol{\theta})\right)^\text{T} \boldsymbol{\Sigma}_{\boldsymbol{y}}^{-1} \left(\boldsymbol{y}_{i} - \boldsymbol{f}(\boldsymbol{x}_{i};\boldsymbol{\theta})\right) \label{eq:wsse}
\end{align}

\noindent The log-likelihood function is:
\begin{align}
    l(\boldsymbol{\theta}, \boldsymbol{y}_1,\cdots, \boldsymbol{y}_{N_{\text{exp}}}) = \log{L(\boldsymbol{\theta}; \boldsymbol{y}_1,\cdots, \boldsymbol{y}_{N_{\text{exp}}})} = -\text{WSSE} + c \label{eq:log-likelihood-wsse}
\end{align}

The maximum likelihood estimator, $\boldsymbol{\hat{\theta}}$, are the values of $\boldsymbol{\theta}$ that maximize Eq.~\ref{eq:log-likelihood-wsse}.
\begin{equation} \label{eq:mle_wsse}
\begin{aligned}
\hat{\boldsymbol{\theta}} & = \arg\max_{\boldsymbol{\theta}} \, l(\boldsymbol{\theta}, \boldsymbol{y}_1,\cdots, \boldsymbol{y}_{N_{\text{exp}}}) = \arg\max_{\boldsymbol{\theta}}\, (-\text{WSSE} + c) \\
& = \arg\max_{\boldsymbol{\theta}}\, (-\text{WSSE}) = \arg\min_{\boldsymbol{\theta}}\, \text{WSSE}
\end{aligned}
\end{equation}

The Fisher information matrix, $\boldsymbol{M}$, about the parameters, $\boldsymbol{\theta}$, contained in the measurements, $\boldsymbol{y}_1,\cdots, \boldsymbol{y}_{N_{\text{exp}}}$, is defined as:
\begin{align}
    \boldsymbol{M} = - \mathbb{E}_{\boldsymbol{\theta}}\left[\frac{\partial^2 l(\boldsymbol{\theta}, \boldsymbol{y}_1,\cdots, \boldsymbol{y}_{N_{\text{exp}}})}{\partial \boldsymbol{\theta} \partial \boldsymbol{\theta}}\right] \label{eq:FIM-definition}
\end{align}

\noindent Substituting Eq.~\ref{eq:log-likelihood-wsse} in Eq.~\ref{eq:FIM-definition} leads to:
\begin{align}
 \boldsymbol{M} = \left(\frac{\partial^2 \text{WSSE}}{\partial \boldsymbol{\theta} \partial \boldsymbol{\theta}}\right)_{\boldsymbol{\theta} = \hat{\boldsymbol{\theta}}} \label{eq:FIM-wsse}
\end{align}

\section{Detailed three-stage membrane model for critical minerals and materials separation} \label{sec:detailed-membrane-model}
This section presents the complete mathematical model of the proposed three-stage membrane cascade to recover critical minerals and materials from secondary sources, such as recycled materials. In this membrane system, we assume that the temperature of the feed solution is constant, so that the density remains unchanged, and that the diafiltrate is pure water with no dissolved ions. The membrane elements of all stages are assumed to be produced from the same material as discussed by Wamble et al. \cite{wamble2022optimal} and therefore are identical.

The following sets are used to define the ionic species, membrane elements, membrane stages, and membrane flows:
\begin{itemize}
    \item $\mathcal{J} = \{1,\, 2\}$ -- set of ionic species
    \item $\mathcal{E} = \{1,\, 2,\, \ldots, 10\}$ -- set of membrane elements
    \item $\mathcal{S} = \{1,\, 2,\, 3\}$ -- set of membrane stages
    \item $\mathcal{U} = \{\text{permeate},\, \text{retentate}\}$ -- set of membrane flows
\end{itemize}

To denote solvent flows in the membrane cascade, the following subscripts are used:
\begin{itemize}
    \item $ff$ -- fresh feed to the membrane system
    \item $df$ -- fresh diafiltrate added to the membrane system
    \item $p$ -- permeate side of the membrane elements
    \item $r$ -- retentate side of the membrane elements
    \item $f$ -- feed side of the membrane elements
\end{itemize}

Wamble et al. \cite{wamble2022optimal} developed the following equation for the mass flow of ionic species across the elements of the membrane stages:
\begin{equation} \label{eq:membrane-separation}
\begin{aligned}
    \frac{dc_{r,\, j,\, e,\, s}}{dz} & = \frac{-J_{w,\, e,\, s}\, w \, c_{r,\,j,\, e, \, s}\, (S_{j,\, e,\, s}\, - \, 1)}
    {q_{r,\, e,\, s}} \\
    & = 
    \frac{-J_{w,\, e,\, s}\, w \, c_{r,\,j,\, e, \, s}\, (S_{j,\, e,\, s}\, - \, 1)}
    {q_{f,\, e,\, s}\,- J_{w,\, e,\, s}\, w \,z}
    \quad \forall \, j \in \mathcal{J}, \, e \in \mathcal{E}, \, s \in \mathcal{S}
\end{aligned}
\end{equation}

In Eq.~\ref{eq:membrane-separation}, Wamble et al. \cite{wamble2022optimal} assumed that $J_{w,\, e,\, s}$ is constant at 0.1 $\text{m}\cdot\text{h}^{-1}$ and $S_{j,\, e,\, s}$ constant at 1.3 and 0.5 for $\text{Cation 1}$ and $\text{Cation 2}$, respectively. In this study, we consider a more detailed representation of the membrane system by implementing pressure-dependent water flux (Eq.~\ref{eq:water_flux}) and concentration-dependent sieving coefficients (Eq.~\ref{eq:sieving_coefficient}).
\begin{align}
     J_{w,\, e,\, s} = L_p\left(\Delta P - \Delta \pi_{e,\, s}\right) 
     \quad \forall \, e \in \mathcal{E}, \, s \in \mathcal{S}\label{eq:water_flux}
\end{align}
\begin{align}
    S_{j,\, e,\, s} = \bar{S}_{j} + \delta_{j}\sum_{k \in \mathcal{J}}c_{f,\, k,\, e,\, s}\, z_{k}^2 
    \quad \forall \, j \in \mathcal{J}, \, e \in \mathcal{E}, \, s \in \mathcal{S} \label{eq:sieving_coefficient}
\end{align}

\noindent where the opposing osmotic pressure is:
\begin{align}
     \Delta \pi_{e,\, s} = RT\sum_{k \in \mathcal{J}}\left(c_{r,\, k,\, e,\, s} - c_{p,\, k,\, e,\, s}\right)
     \quad \forall \, e \in \mathcal{E}, \, s \in \mathcal{S}\label{eq:osmotic_pressure}
\end{align}

\noindent Substituting Eq.~\ref{eq:water_flux} and Eq.~\ref{eq:sieving_coefficient} into Eq.~\ref{eq:membrane-separation} leads to:
\begin{align}
    \frac{dc_{r,\, j,\, e,\, s}}{dz} = \frac{\beta \, c_{r,\,j,\, e, \, s}\, (\alpha - \gamma \, c_{r,\,j,\, e, \, s})}
    {q_{f,\, e,\, s}\,- (\alpha - \gamma \, c_{r,\,j,\, e, \, s})z} \quad \forall \, j \in \mathcal{J}, \, e \in \mathcal{E}, \, s \in \mathcal{S}  \label{eq:lumped_membrane_model}
\end{align}

\noindent where
\begin{align*}
     \beta = 1 - S_{j,\, e,\, s}, \quad \gamma = L_pwRT(1 - S_{j,\, e,\, s})
     \quad \forall \, j \in \mathcal{J}, \, e \in \mathcal{E}, \, s \in \mathcal{S}
\end{align*}

\noindent and
\begin{align*}
     \alpha = L_pw\left(\Delta P - RT(1 - S_{k,\, e,\, s})c_{r,\, k,\, e,\, s}\right)
     \quad \text{for} \, k \in \mathcal{J} \, \text{and} \, k\neq j, \, \forall \, e \in \mathcal{E}, \, s \in \mathcal{S}
\end{align*}

Eq.~\ref{eq:lumped_membrane_model} does not have an analytical solution. Numerical methods such as the finite difference backward approximation were used to solve Eq.~\ref{eq:lumped_membrane_model} over the elements of the membrane stages, where each element is treated as a finite volume.

The mass flow of solvent across each element of the membrane stage is defined as:
\begin{align}
     M_{\text{solvent},\, u,\, e,\, s} = \frac{J_{w,\, e,\, s}\,L\,w\, \rho_{\text{solvent}}}{N}
\end{align}

As shown in Figure~\ref{fig:membrane-cascade}, the permeate outlet from stage 1 and the recycled retentate outlet from stage 3 are the feed streams to stage 2. The flowrate and concentrations of the combined feed to stage 2 are calculated as follows:
\begin{align}
     q_{f,\, 1,\,2} = q_{r,\,10,\,3} + q_{p,\,10,\,1}
\end{align}
\begin{align}
     c_{f,\,j,\,1,\,2} = \frac{q_{r,\,10,\,3}\,c_{r,\,j,\,10,\,3} + 
     q_{p,\,10,\,1}\, c_{p,\, j,\,10,\,1}}{q_{f,\, 1,\,2}}
     \quad \forall \, j \in \mathcal{J}
\end{align}

Similarly, the permeate outlet from stage 2, combined with the fresh diafiltrate, forms the feed to stage 3 (Eq.~\ref{eq:stage_3_flow} and~\ref{eq:stage_3_conc}). In addition, the retentate from element 9 of stage 3, when combined with the fresh feed, constitutes the feed stream to element 10 of stage 3 (Eq.~\ref{eq:stage_3_side_flow} and ~\ref{eq:stage_3_side_conc}).
\begin{align}
     q_{f,\, 1,\,3} = q_{df} + q_{p,\,10,\,2} \label{eq:stage_3_flow}
\end{align}
\begin{align}
     c_{f,\,j,\,1,\,3} = \frac{q_{df}\,c_{df,\,j} + 
     q_{p,\,10,\,2}\, c_{p,\, j,\,10,\,2}}{q_{f,\, 1,\,3}}
     \quad \forall \, j \in \mathcal{J} \label{eq:stage_3_conc}
\end{align}
\begin{align}
     q_{f,\, 10,\,3} = q_{ff} + q_{r,\,9,\,3} \label{eq:stage_3_side_flow}
\end{align}
\begin{align}
     c_{f,\,j,\,10,\,3} = \frac{q_{ff}\,c_{ff,\,j} + 
     q_{r,\,9,\,3}\, c_{r,\, j,\,9,\,3}}{q_{f,\, 10,\,3}}
     \quad \forall \, j \in \mathcal{J} \label{eq:stage_3_side_conc}
\end{align}

Lastly, the recycled retentate outlet from stage 2 is the feed to stage 1:
\begin{align}
     q_{f,\, 1,\,1} = q_{r,\,10,\,2} \label{eq:stage_1_flow}
\end{align}
\begin{align}
     c_{f,\,j,\,1,\,1} = c_{r,\,j,\,10,\,2}
     \quad \forall \, j \in \mathcal{J} \label{eq:stage_1_conc}
\end{align}

The known parameters in this model and their values are presented in Table~\ref{tab:known-model-parameters}.
\begin{table}[h!]
\centering
\caption{Known parameters in the three-stage membrane cascade.}
\begin{tabular}{llll}
\hline
\textbf{Parameter} & \textbf{Unit} & \textbf{Description} & \textbf{Value} \\
\hline
$R$ & $\text{J}\cdot\text{mol}^{-1}\cdot\text{K}^{-1}$ & Gas constant & 8.314 \\
$\rho_{\text{solvent}}$ & $\text{kg}\cdot\text{m}^{-3}$ & Density of water & 1000 \\
$\Delta P$ & $\text{Pa}$ & Applied membrane pressure & $10^6$ \\
$T$ & $\text{K}$ & Temperature & 298.15 \\
$w$ & $\text{m}$ & Width of membrane elements & 1.5 \\
$L$ & $\text{m}$ & Length of membrane stages & 200 \\
$N$ & dimensionless & Number of elements in a stage & 10 \\
$z_{1}$ & dimensionless & Valency of cation 1 & 1 \\
$z_{2}$ & dimensionless & Valency of cation 2 & 2 \\
\hline
\end{tabular}
\label{tab:known-model-parameters}
\end{table}

\section{Optimal design criteria first and second derivative derivations} \label{sec:in-depth-criteria-derivations}
The following section is concerned with the derivatives of optimality criteria for A-, D-, E-, and ME-optimality conditions which, as shown previously in the manuscript, are of the following forms:
\begin{itemize}
\item \textbf{A-optimality}: 

\begin{align}
    \min\text{trace}\left(\boldsymbol{M}^{-1}\right) = \sum_{i}^p \frac{1}{\lambda_i} \label{eq:obj-trace-SI}
\end{align}

\item \textbf{D-optimality}:

\begin{align}
    \max\left|\boldsymbol{M}\right| = \prod_{i}^p \lambda_i \label{eq:obj-log-det-SI}
\end{align}

\item \textbf{E-optimality}:

\begin{align}
    \max \;\min_{i\in\left\{1, \ldots, p\right\}}\lambda_i \label{eq:obj-min-eig-SI}
\end{align}

\item \textbf{ME-optimality}:

\begin{align}
    \min \; \frac{\max_{i\in\left\{1, \ldots, p\right\}}\lambda_i}{\min_{i\in\left\{1, \ldots, p\right\}}\lambda_i} \equiv \min \kappa \label{eq:obj-cond-SI}
\end{align}
\end{itemize}

\subsection{A-optimality derivatives}
The first challenge in the derivative of A-optimality is that we need the derivative of a matrix with respect to itself, as shown below:

\begin{align}
    \Psi_\text{A} = &\text{trace}\left(\boldsymbol{M}^{-1}\right) \label{eq:obj-trace-si} \\
    \frac{\partial \Psi_\text{A}}{\partial \boldsymbol{M}} = & \frac{\partial \,\text{trace}\left(\boldsymbol{M}^{-1}\right)}{\partial \boldsymbol{M}} \\
     = & \frac{\partial}{\partial \boldsymbol{M}} \left(\sum_i \left(\boldsymbol{M}^{-1}\right)_{ii}\right) \label{eq:trace-sum-deriv-si}
\end{align}

This begs the question: \emph{what are each of the terms of the sum in Eq.~\ref{eq:trace-sum-deriv-si}}? Park's paper from 2018 has material for learning some basic matrix/tensor derivatives \cite{park2018tensor-calc}. To find this, we utilize a trick to determine the derivative of the inverse of the matrix with respect to itself:

\begin{align}
    \delta_{ab} = & M_{ac}M^{-1}_{cb} \\
    \frac{\partial \delta_{ab}}{\partial M_{de}} = & \frac{\partial M_{ac}}{\partial M_{de}}M^{-1}_{cb} + M_{ac}\frac{\partial M^{-1}_{cb}}{\partial M_{de}} \\
    0 = & \frac{\partial M_{ac}}{\partial M_{de}}M^{-1}_{cb} + M_{ac}\frac{\partial M^{-1}_{cb}}{\partial M_{de}} \\
    M_{ac}\frac{\partial M^{-1}_{cb}}{\partial M_{de}} = & - \frac{\partial M_{ac}}{\partial M_{de}}M^{-1}_{cb} \\
    M^{-1}_{fa}M_{ac}\frac{\partial M^{-1}_{cb}}{\partial M_{de}} = & - M^{-1}_{fa}\frac{\partial M_{ac}}{\partial M_{de}}M^{-1}_{cb} \\
    \delta_{fc} \frac{\partial M^{-1}_{cb}}{\partial M_{de}} = & - M^{-1}_{fa}\frac{\partial M_{ac}}{\partial M_{de}}M^{-1}_{cb} \label{eq:trace-break-si}
\end{align}
where $\delta_{ij}$ represents the Kronecker delta, where the value of $\delta_{ij}$ is 1 when $i = j$ and 0 otherwise. It is important at this point to note the derivative of a matrix with respect to itself:

\begin{align}
    \frac{\partial M_{ij}}{\partial M_{kl}} = & \delta_{il}\delta_{jk} \label{eq:matrix-derivative-with-itself-si}
\end{align}
where the result is a fourth-order tensor where the element is 1 if $i = l$ and $j = k$. Importantly, if using a column-wise definition of the derivative operator, the indices will change for the right-hand side expression. Also, an important definition for tensor products using summation notation (where the sum can be implied by matching the inner-most indices) is the following:

\begin{align}
    \boldsymbol{A}\boldsymbol{B} = \sum_j A_{ij}B_{jk} = A_{ij}B_{jk}
\end{align}

Continuing from Eq.~\ref{eq:trace-break-si} with Eq.~\ref{eq:matrix-derivative-with-itself-si} in mind, we have that:

\begin{align}
    \delta_{fc}\frac{\partial M^{-1}_{cb}}{\partial M_{de}} = & - M^{-1}_{fa}\delta_{ae}\delta_{cd}M^{-1}_{cb} \\
    \frac{\partial M^{-1}_{fb}}{\partial M_{de}} = & - M^{-1}_{fe}M^{-1}_{db}
\end{align}
or in a more useful format:
\begin{align}
    \frac{\partial M^{-1}_{ij}}{\partial M_{kl}} = - M^{-1}_{il}M^{-1}_{kj}
\end{align}

Now that this derivative has been established, we can continue from Eq.~\ref{eq:trace-sum-deriv-si}:

\begin{align}
    \frac{\partial \Psi_\text{A}}{\partial M_{ij}} = & \frac{\partial}{\partial M_{ij}} \left(\sum_g \left(\boldsymbol{M}^{-1}\right)_{gg}\right)\\
    = & \sum_g \frac{\partial M^{-1}_{gg}}{\partial M_{ij}} \\
    = & \sum_g \left(- M^{-1}_{gj}M^{-1}_{ig}\right) \\
    = &- M^{-1}_{ig}M^{-1}_{gj} = -\boldsymbol{M}^{-1}\boldsymbol{M}^{-1}
\end{align}
Now, taking the second derivative of A optimality is not too challenging. For the second derivatives, we use a scalar representation of the $\left(i, j, k, l\right)$-th element of the fourth order tensor. We also assume that the FIM, $\boldsymbol{M}$, is full rank and symmetric, by definition, and thus the transpose of its inverse is its inverse. With this in mind, we define the second derivative as follows:

\begin{align}
    \frac{\partial^2 \Psi_\text{A}}{\partial M_{ij} \partial M_{kl}} = & - \frac{\partial}{\partial M_{kl}} M^{-1}_{ig}M^{-1}_{gj} \\
    = & -\left(\frac{\partial M^{-1}_{ig}}{\partial M_{kl}}M^{-1}_{gj} + M^{-1}_{ig}\frac{\partial M^{-1}_{gj}}{\partial M_{kl}}\right) \\
    = & M^{-1}_{il}M^{-1}_{kg}M^{-1}_{gj} + M^{-1}_{ig}M^{-1}_{gl}M^{-1}_{kj} \\
    = & M^{-1}_{il}\left(\boldsymbol{M}^{-1}\boldsymbol{M}^{-1}\right)_{kj} + \left(\boldsymbol{M}^{-1}\boldsymbol{M}^{-1}\right)_{il}M^{-1}_{kj} \label{eq:trace-second-derivative-easier-to-see-si}
\end{align}
where Eq.~\ref{eq:trace-second-derivative-easier-to-see-si} is shown to highlight that there are matrix inverse terms and matrix inverse squared terms in this second derivative. This makes sense as the scalar derivative of $x^n$ reduces the power by 1 for each derivative and Eq.~\ref{eq:trace-second-derivative-easier-to-see-si} contains terms close to $\boldsymbol{M}^{-3}$.

\subsection{D-optimality derivatives}
\begin{align}
    \Psi_\text{D} = &\text{log}\left(\left|\boldsymbol{M}\right|\right) \label{eq:obj-determinant-si} \\
    \frac{\partial \Psi_\text{D}}{\partial \boldsymbol{M}} = & \frac{\partial \,\text{log}\left(\left|\boldsymbol{M}\right|\right)}{\partial \boldsymbol{M}} \\
    = & \frac{1}{\left|\boldsymbol{M}\right|}\frac{\partial \,\left|\boldsymbol{M}\right|}{\partial \boldsymbol{M}}
\end{align}
At this point, it is good to use some definitions from linear algebra relating the adjoint matrix to the original matrix and the determinant as well as the inverse of the matrix.

\begin{align}
    \left|\boldsymbol{M}\right| = & \sum_j M_{ij}\alpha_{ij} \\
    \text{adj}\left(\boldsymbol{M}\right) = & \left|\boldsymbol{M}\right| \boldsymbol{M}^{-1}
\end{align}
where $\alpha_{ij}$ is the $\left(i, j\right)$-th element of the cofactor matrix, whose transpose is the adjugate matrix, represented by the function $\text{adj}\left(\right)$. Using these definitions along with the symmetry of the FIM, we have the following continuation of the derivative of D-optimality:

\begin{align}
    \frac{\partial \Psi_\text{D}}{\partial \boldsymbol{M}} = & \frac{1}{\left|\boldsymbol{M}\right|}\frac{\partial \,}{\partial \boldsymbol{M}}\left(\sum_j M_{ij}\alpha_{ij}\right) \\
    = & \frac{1}{\left|\boldsymbol{M}\right|}\frac{\partial \,}{\partial \boldsymbol{M}}\left(\frac{1}{2}\sum_j M_{ij}\alpha_{ij} + M_{ji}\alpha_{ij}\right) \\
    = & \frac{1}{2\left|\boldsymbol{M}\right|} \sum_j \alpha_{ij}\left(\frac{\partial M_{ij}\,}{\partial \boldsymbol{M}} + \frac{\partial M_{ji}\,}{\partial \boldsymbol{M}}\right) \\
    = & \frac{1}{2\left|\boldsymbol{M}\right|} \sum_j \alpha_{ij}\left(\delta_{il}\delta_{jk} + \delta_{jl}\delta_{ik}\right) \\
    = & \frac{1}{2\left|\boldsymbol{M}\right|} \left(\alpha_{ik}\delta_{il} + \alpha_{il}\delta_{ik}\right) \\
    = & \frac{1}{2\left|\boldsymbol{M}\right|} \left(\alpha_{lk} + \alpha_{kl}\right) \\
    = & \frac{1}{2\left|\boldsymbol{M}\right|} \left(\text{adj}\left(\boldsymbol{M}\right)^T + \text{adj}\left(\boldsymbol{M}\right)\right) \\
    = & \frac{1}{2\left|\boldsymbol{M}\right|} \left(\left|\boldsymbol{M}\right|\boldsymbol{M}^{-T} + \left|\boldsymbol{M}\right|\boldsymbol{M}^{-1}\right) \\
    = & \frac{1}{2} \left(\boldsymbol{M}^{-T} + \boldsymbol{M}^{-1}\right)
\end{align}

Now, since the first derivative for D-optimality is in terms of the inverse of the matrix, and assuming that the inverse of the matrix is symmetric (by definition), we can continue with the second derivative below:

\begin{align}
    \frac{\partial^2 \Psi_\text{D}}{\partial M_{ij}\partial M_{kl}} = & \frac{\partial}{\partial M_{kl}} \left( \boldsymbol{M}^{-1}\right) \\
    = & -M^{-1}_{il}M^{-1}_{kj}
\end{align}

\subsection{E-optimality derivatives} \label{sec:E-opt-SI}
For E-optimality, a few steps are more in-depth than the previous derivatives. First, we pose the first derivative as follows:
\begin{align}
    \Psi_\text{E} = &\min_i \lambda_i \label{eq:obj-minimum-eigenvalue-si} \\
    \frac{\partial \Psi_\text{E}}{\partial \boldsymbol{M}} = & \frac{\partial \,\lambda_\text{min}}{\partial \boldsymbol{M}}
\end{align}
Importantly, we will establish a derivative for any eigenvalue $\lambda_s$ and thus can simply utilize the minimum operator to find the minimum eigenvalue and corresponding eigenvector. It is possible to pose these types of max-min problems as a single program, but with the Grey Box abstraction, we do not need to over complicate the process and can embed the inner minimization problem within the Grey Box. For any eigenvalue $\lambda_s$, the following is true (assuming full rank):

\begin{align}
    \boldsymbol{M}\boldsymbol{v}_s& = \lambda_s\boldsymbol{v}_s &\forall s \in \left\{1, \ldots, p\right\} \label{eq:eigenvalue-problem-si} \\
    \boldsymbol{v}_s\cdot \boldsymbol{v}_s & = 1 &\forall s \in \left\{1, \ldots, p\right\}\label{eq:eigenvector_orthnormality-si} \\
    \boldsymbol{v}_i\cdot \boldsymbol{v}_j & = 0 &\forall i \neq j; \left(i, j\right) \in \left\{1, \ldots, p\right\}^2\label{eq:eigenvector_orthonogonality-si}
\end{align}
Therefore, the important equation for E-optimality becomes:

\begin{align}
    \boldsymbol{M}\boldsymbol{v}_\text{min} = & \lambda_\text{min}\boldsymbol{v}_\text{min} = \Psi_E \boldsymbol{v}_\text{min} \\
    \frac{\partial \boldsymbol{v}_\text{min}\Psi_\text{E}}{\partial \boldsymbol{M}} = & \frac{\partial \boldsymbol{M} \boldsymbol{v}_\text{min}}{\partial \boldsymbol{M}} \\
    \lambda_\text{min} \frac{\partial \boldsymbol{v}_\text{min}}{\partial \boldsymbol{M}} + \boldsymbol{v}_\text{min}\frac{\partial \Psi_\text{E}}{\partial \boldsymbol{M}} = &  \frac{\partial \boldsymbol{M}}{\partial \boldsymbol{M}}\boldsymbol{v}_\text{min} + \boldsymbol{M} \frac{\partial \boldsymbol{v}_\text{min}}{\partial \boldsymbol{M}} \label{eq:E-opt-deriv-pause-si}
\end{align}
Importantly, we now introduce the derivative of Eq.~\ref{eq:eigenvector_orthnormality-si} with respect to $\boldsymbol{M}$:
\begin{align}
    \frac{\partial \boldsymbol{v}_s}{\partial \boldsymbol{M}} \cdot \boldsymbol{v}_s + \boldsymbol{v}_s\cdot\frac{\partial \boldsymbol{v}_s}{\partial \boldsymbol{M}} & = 0 &\forall s \in \left\{1, \ldots, p\right\}\\
    \frac{\partial \boldsymbol{v}_s}{\partial \boldsymbol{M}} \cdot \boldsymbol{v}_s = \boldsymbol{v}_s\cdot\frac{\partial \boldsymbol{v}_s}{\partial \boldsymbol{M}}  &= 0  &\forall s \in \left\{1, \ldots, p\right\} \label{eq:E-opt-zero-value-derivative-relationship-si}
\end{align}
Note that since the dot product is commutative, we can arrive at Eq.~\ref{eq:E-opt-zero-value-derivative-relationship-si}. Notably, this is useful when multiplying the entirety of Eq.~\ref{eq:E-opt-deriv-pause-si} by $\boldsymbol{v}_\text{min}^T$:

\begin{align}
    \boldsymbol{v}_\text{min}^T\lambda_\text{min} \frac{\partial \boldsymbol{v}_\text{min}}{\partial \boldsymbol{M}} + \boldsymbol{v}_\text{min}^T\boldsymbol{v}_\text{min}\frac{\partial \Psi_\text{E}}{\partial \boldsymbol{M}} = &  \boldsymbol{v}_\text{min}^T\frac{\partial \boldsymbol{M}}{\partial \boldsymbol{M}}\boldsymbol{v}_\text{min} + \boldsymbol{v}_\text{min}^T\boldsymbol{M} \frac{\partial \boldsymbol{v}_\text{min}}{\partial \boldsymbol{M}}
\end{align}
Now, using the definition in Eq.~\ref{eq:eigenvalue-problem-si} and the fact the $\boldsymbol{M}$ is symmetric, we can replace $\boldsymbol{v}_\text{min}^T\boldsymbol{M}$ with $\lambda_\text{min} \boldsymbol{v}_\text{min}^T$: 

\begin{align}
    \boldsymbol{v}_\text{min}^T\lambda_\text{min} \frac{\partial \boldsymbol{v}_\text{min}}{\partial \boldsymbol{M}} + \boldsymbol{v}_\text{min}^T\boldsymbol{v}_\text{min}\frac{\partial \Psi_\text{E}}{\partial \boldsymbol{M}} = &  \boldsymbol{v}_\text{min}^T\frac{\partial \boldsymbol{M}}{\partial \boldsymbol{M}}\boldsymbol{v}_\text{min} + \lambda_\text{min} \boldsymbol{v}_\text{min}^T \frac{\partial \boldsymbol{v}_\text{min}}{\partial \boldsymbol{M}}
\end{align}
Then, the definition in Eq.~\ref{eq:E-opt-zero-value-derivative-relationship-si} allows us to equate two terms in the expression to zero. When combined with the understanding that $\frac{\partial \Psi_\text{E}}{\partial \boldsymbol{M}}$ is a scalar, we arrive at a rather simple expression:

\begin{align}
    \boldsymbol{v}_\text{min}^T \boldsymbol{v}_\text{min}\frac{\partial \Psi_\text{E}}{\partial \boldsymbol{M}} = &  \boldsymbol{v}_\text{min}^T\frac{\partial \boldsymbol{M}}{\partial \boldsymbol{M}}\boldsymbol{v}_\text{min} \\
    \frac{\partial \Psi_\text{E}}{\partial \boldsymbol{M}} = &  v_{\text{min},i} \delta_{il}\delta_{jk}v_{\text{min},k} \\
    \frac{\partial \Psi_\text{E}}{\partial \boldsymbol{M}} = &  v_{\text{min},l}v_{\text{min},k} \\
    \frac{\partial \Psi_\text{E}}{\partial \boldsymbol{M}} = &\boldsymbol{v}_\text{min}\boldsymbol{v}_\text{min}^T
\end{align}
With this conclusion, we can now pose the second derivative of E-optimality as follows:

\begin{align}
    \frac{\partial^2 \Psi_\text{E}}{\partial M_{ij} \partial M_{kl}} = & \frac{\partial}{\partial M_{kl}}\left(v_{\text{min}, i}v_{\text{min}, j}\right) \\
    = & \frac{\partial v_{\text{min}, i}}{\partial M_{kl}}v_{\text{min}, j} + v_{\text{min}, i}\frac{\partial v_{\text{min}, j}}{\partial M_{kl}}
\end{align}
From this result, it becomes clear that to get the second derivative for E-optimality, we need the derivative of each eigenvector $\boldsymbol{v}_i$ with respect to $\boldsymbol{M}$. To do this, we return to Eq.~\ref{eq:E-opt-deriv-pause-si} and instead multiply by an eigenvector $\boldsymbol{v}_s$ where $s \neq i$:

\begin{align}
    \boldsymbol{v}_s^T \lambda_\text{min} \frac{\partial \boldsymbol{v}_\text{min}}{\partial \boldsymbol{M}} + \boldsymbol{v}_s^T \boldsymbol{v}_\text{min}\frac{\partial \Psi_\text{E}}{\partial \boldsymbol{M}} = &  \boldsymbol{v}_s^T \frac{\partial \boldsymbol{M}}{\partial \boldsymbol{M}}\boldsymbol{v}_\text{min} + \boldsymbol{v}_s^T \boldsymbol{M} \frac{\partial \boldsymbol{v}_\text{min}}{\partial \boldsymbol{M}} \label{eq:E-opt-second-deriv-pause-si} \\
    & \qquad\qquad\qquad\forall s \neq \text{min} \in \left\{1, \ldots, p\right\} \notag
\end{align}
Once again using Eq.~\ref{eq:eigenvalue-problem-si}, we can simplify the previous equation to the following:

\begin{align}
    \boldsymbol{v}_s^T \lambda_\text{min} \frac{\partial \boldsymbol{v}_\text{min}}{\partial \boldsymbol{M}} + \boldsymbol{v}_s^T \boldsymbol{v}_\text{min}\frac{\partial \Psi_\text{E}}{\partial \boldsymbol{M}} = &  \boldsymbol{v}_s^T \frac{\partial \boldsymbol{M}}{\partial \boldsymbol{M}}\boldsymbol{v}_\text{min} + \lambda_s \boldsymbol{v}_s^T \frac{\partial \boldsymbol{v}_\text{min}}{\partial \boldsymbol{M}} \\
    & \qquad\qquad\qquad\forall s \neq \text{min} \in \left\{1, \ldots, p\right\} \notag
\end{align}
We can then use the identity in Eq.~\ref{eq:eigenvector_orthnormality-si} to get rid of one term on the left-hand side:

\begin{align}
    \boldsymbol{v}_s^T \lambda_\text{min} \frac{\partial \boldsymbol{v}_\text{min}}{\partial \boldsymbol{M}} = &  \boldsymbol{v}_s^T \frac{\partial \boldsymbol{M}}{\partial \boldsymbol{M}}\boldsymbol{v}_\text{min} + \lambda_s \boldsymbol{v}_s^T \frac{\partial \boldsymbol{v}_\text{min}}{\partial \boldsymbol{M}} & \forall s \neq \text{min} \in \left\{1, \ldots, p\right\}
\end{align}
Finally, we organize like terms and simplify the $\frac{\partial \boldsymbol{M}}{\partial \boldsymbol{M}}$ term on the right-hand side to get the following:

\begin{align}
    \left(\lambda_\text{min} - \lambda_s\right)\boldsymbol{v}_s^T \frac{\partial \boldsymbol{v}_\text{min}}{\partial \boldsymbol{M}} = &  \boldsymbol{v}_s\boldsymbol{v}_\text{min}^T & \forall s \neq \text{min} \in \left\{1, \ldots, p\right\} \\
    \boldsymbol{v}_s^T \frac{\partial \boldsymbol{v}_\text{min}}{\partial \boldsymbol{M}} = &  \frac{1}{\lambda_\text{min} - \lambda_s}\boldsymbol{v}_s\boldsymbol{v}_\text{min}^T & \forall s \neq \text{min} \in \left\{1, \ldots, p\right\} \label{eq:eigenvector-derivative-before-simplification-si}
\end{align}
We now have a system of equations with $s$ equations and $s$ unknowns with $s-1$ equations coming from Eq.~\ref{eq:eigenvector-derivative-before-simplification-si} and one extra coming from Eq.~\ref{eq:E-opt-zero-value-derivative-relationship-si} when $s = \text{min}$. The solution to this system when $\boldsymbol{M}$ is symmetric is a conclusion from perturbation theory of linear operators and is an important conclusion in Kato's book \cite{kato2013perturbation}. Using these facts, we can reduce this equation system to a sum:

\begin{align}
    \frac{\partial \boldsymbol{v}_s}{\partial \boldsymbol{M}} = &  \sum_{r \neq s} \left(\left(\frac{1}{\lambda_s - \lambda_r}\boldsymbol{v}_r\boldsymbol{v}_s^T\right)\boldsymbol{v}_r\right) & \forall s \in \left\{1, \ldots, p\right\} \label{eq:eigenvector-derivative-si}
\end{align}
Another option to solve this system is to solve the square, linear system using algebra and the solution is identical to the formula in Eq.~\ref{eq:eigenvector-derivative-si}. Now that we have the derivative formula for any eigenvector $\boldsymbol{v}_r$, we can finish the second derivative of any eigenvalue, but, more specifically, the minimum eigenvalue, as shown below:

\begin{align}
    \frac{\partial^2 \Psi_\text{E}}{\partial M_{ij} \partial M_{kl}} = &\sum_{r \neq \text{min}} \left(\frac{1}{\lambda_\text{min} - \lambda_r}v_{r, l}v_{\text{min}, k}v_{r,i}\right)v_{\text{min}, j} \notag\\
    &+ v_{\text{min}, i}\sum_{r \neq \text{min}} \left(\frac{1}{\lambda_\text{min} - \lambda_r}v_{r, k}v_{\text{min}, l}v_{r, j}\right)
\end{align}

\subsection{ME-optimality derivatives}
Given that ME-optimality is represented completely as a function of eigenvalues, one can construct the overall formulae for the first and second derivatives using those in section~\ref{sec:E-opt-SI} and using Eq.~\ref{eq:ME-opt-first-deriv} and~\ref{eq:ME-opt-second-deriv} as follows:

\begin{align}
    \Psi_{ME} = &\;\text{ln}\left(\frac{\lambda_\text{max}}{\lambda_\text{min}}\right)\\
    \frac{\partial\,\Psi_{ME}}{\partial \boldsymbol{M}} = &\frac{1}{\lambda_\text{max}} \frac{\partial \lambda_\text{max}}{\partial \boldsymbol{M}} - \frac{1}{\lambda_\text{min}} \frac{\partial \lambda_\text{min}}{\partial \boldsymbol{M}} \\
    \frac{\partial\,\Psi_{ME}}{\partial \boldsymbol{M}}= & \frac{1}{\lambda_\text{max}} \boldsymbol{v}_\text{max} \boldsymbol{v}_\text{max}^T - \frac{1}{\lambda_\text{min}} \boldsymbol{v}_\text{min} \boldsymbol{v}_\text{min}^T \label{eq:ME-opt-first-deriv-si}\\ 
    \frac{\partial^2 \Psi_\text{E}}{\partial M_{ij} \partial M_{kl}} = &\frac{1}{\lambda_\text{max}}  \left(\sum_{r \neq s} \left(\frac{1}{\lambda_\text{max} - \lambda_r}v_{r, l}v_{\text{max}, k}v_{r,i}\right)v_{\text{max}, j}\right. \notag\\
    &\qquad\;\;\;+ v_{\text{max}, i}\left.\sum_{r \neq s} \left(\frac{1}{\lambda_\text{max} - \lambda_r}v_{r, k}v_{\text{max}, l}v_{r, j}\right)\right) \notag\\
    & - \frac{1}{\lambda_\text{max}^2} v_{\text{max}, l} v_{\text{max}, k} v_{\text{max}, j} v_{\text{max}, i}\notag \\ 
    & + \frac{1}{\lambda_\text{min}^2} v_{\text{min}, l} v_{\text{min}, k} v_{\text{min}, j} v_{\text{min}, i} \label{eq:ME-opt-second-deriv-si} \\ 
    & - \frac{1}{\lambda_\text{min}} \left(\sum_{r \neq s} \left(\frac{1}{\lambda_\text{min} - \lambda_r}v_{r, l}v_{\text{min}, k}v_{r,i}\right)v_{\text{min}, j}\right. \notag\\
    &\qquad\;\;\;\;\;\;\;+ v_{\text{min}, i}\left.\sum_{r \neq s} \left(\frac{1}{\lambda_\text{min} - \lambda_r}v_{r, k}v_{\text{min}, l}v_{r, j}\right)\right) \notag
\end{align}

\section{Numerical confirmation of analytical derivatives} \label{sec:si-numerical-things}
Although these formulas are mathematically consistent, we also confirm that these results align with a numerical approximation of the derivative. We employ a finite difference perturbation of each element of a randomly generated, symmetric matrix and perform an element-wise comparison between the numerical derivative and the respective formula above.

\begin{figure}[t]
\includegraphics[width = \textwidth]{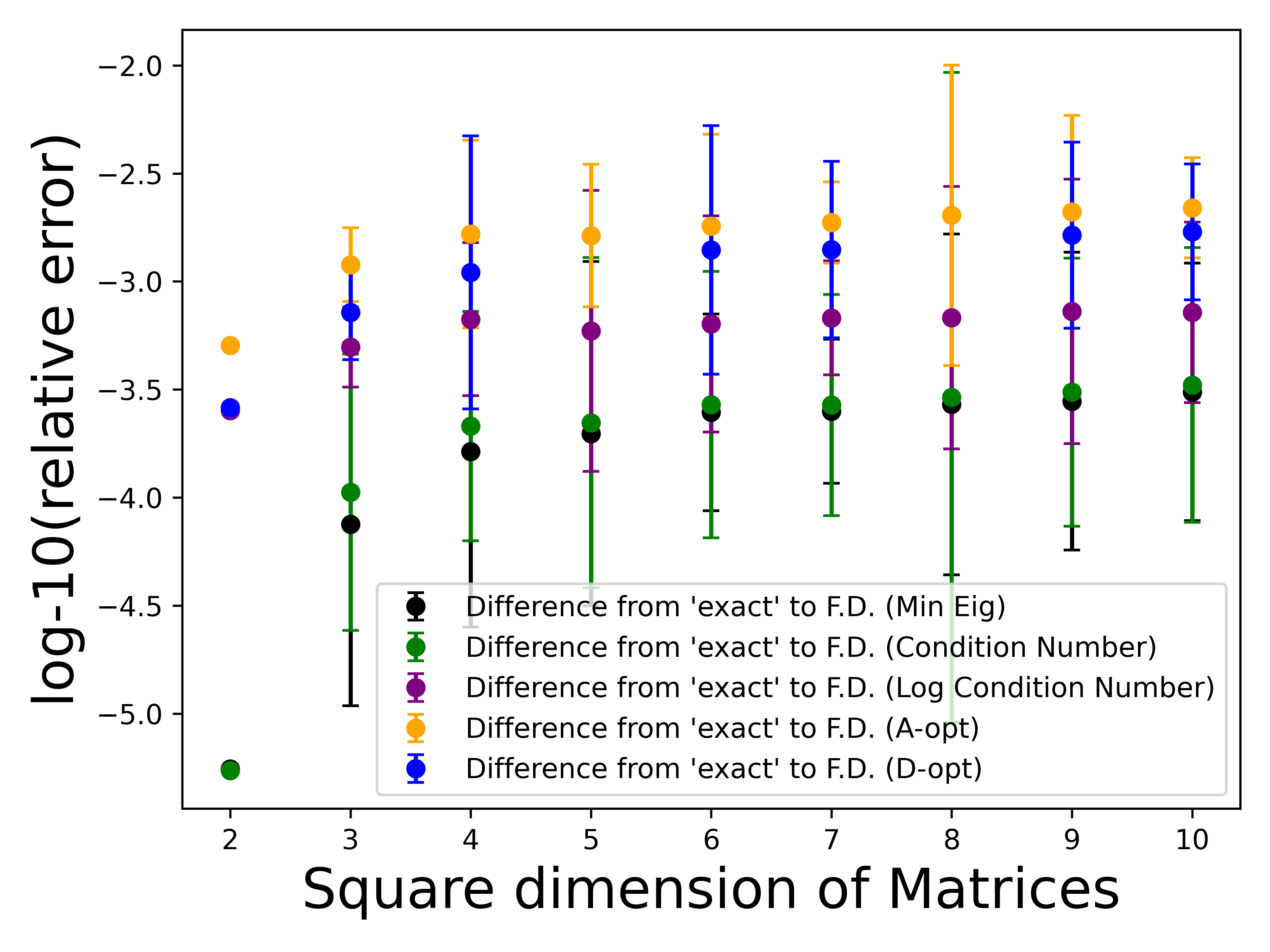}
\caption{The mean (points) and standard error (error bars) plotted for the difference between the finite difference and exact derivative evaluations for the specified criteria. The differences are plotted as log of the relative error versus the size of the square matrix. For each size, 100 randomly generated matrices were analyzed to generate the error bars and mean values.} \label{fig:si-first-deriv-numerics}
\end{figure}

For the first derivative, we used 100 random samples of square matrices ranging from 2-by-2 to 10-by-10. Using a forward difference with a perturbation step size of $10^{-4}$, the  difference between the numerical and exact derivatives are shown in Figure~\ref{fig:si-first-deriv-numerics}. As seen, the derivative values are close to the tolerance used, indicating that the exact derivative is a correct representation of the derivative of each criterion.

\begin{figure}[t]
\includegraphics[width = \textwidth]{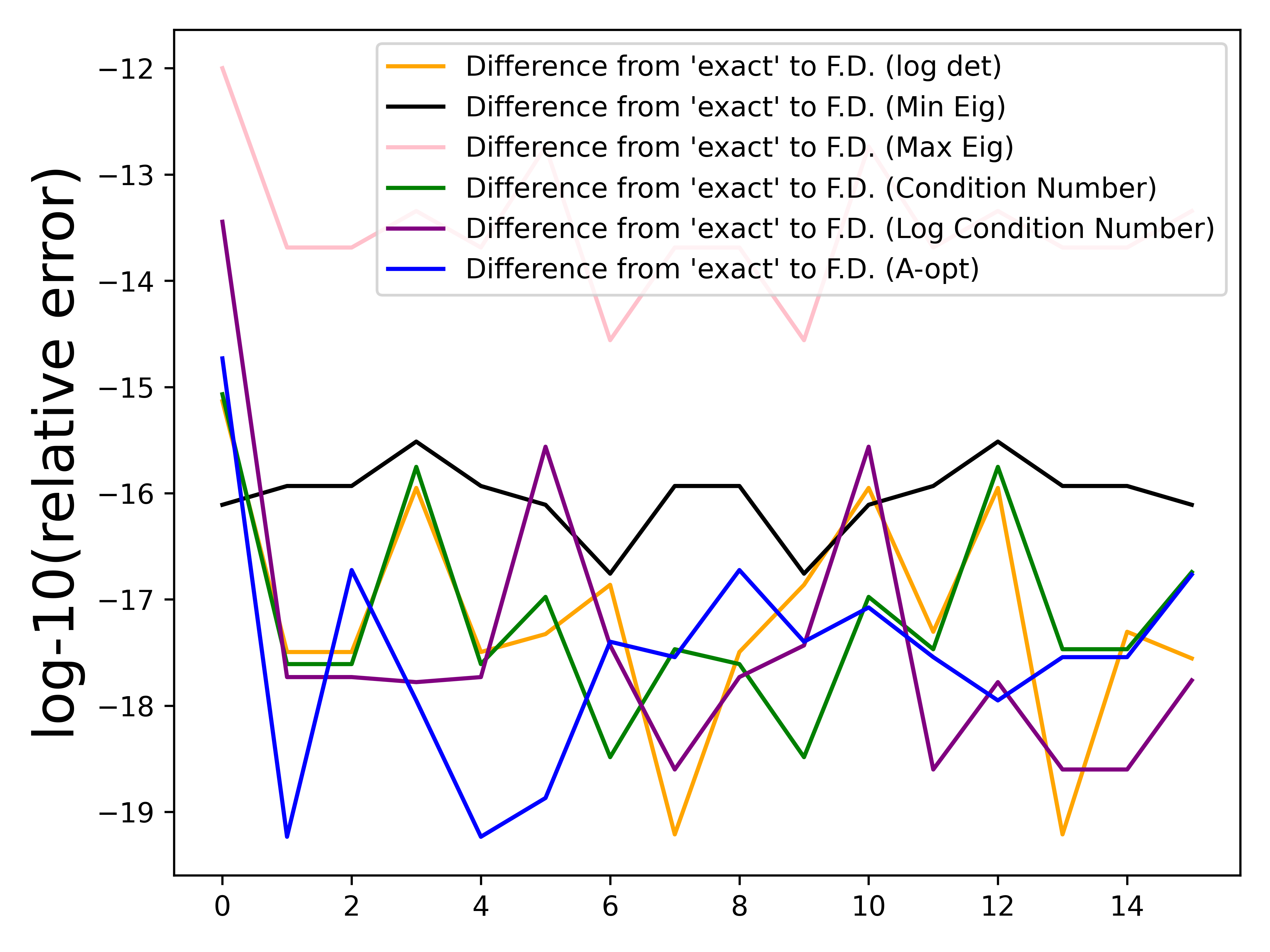}
\caption{Difference between finite difference and exact values for the second derivatives of each condition. Here, we show the log of the relative error for each component of the 4-th order tensor, with 16 total elements.} \label{fig:si-second-deriv-numerics}
\end{figure}

For the second derivative, we use a single, randomly generated square matrix of size 2-by-2 to test the differences between the numerical and exact representations. Here, the result is a 2-by-2-by-2-by-2, 4-th order tensor, and each element is compared. For this case, we utilized a central difference formula for second derivatives and a perturbation step size of $10^{-6}$. As shown in Figure~\ref{fig:si-second-deriv-numerics}, we can see that all the second derivative criteria are well within the expected error and are an indicator that the exact second derivative is a correct representation.

All these derivatives can be tested using the files at the following repository for the curious reader: \url{https://github.com/djlaky/eigenvalue_derivatives}. If desired, the user can adjust the size of the square matrix and the perturbation step size to confirm for themselves that these derivative formulas are an adequate representation both mathematically and numerically.

\subsection{Condition number numerical intricacies} \label{sec:si-numerical-problems}
An early version of Figure~\ref{fig:si-second-deriv-numerics} led us to question the formula derived in Eq.~\ref{eq:ME-opt-second-deriv-si}. However, the calculus is correct, and for educational purposes, we include Figure~\ref{fig:si-second-deriv-numerics-bad} below.

\begin{figure}[t]
\includegraphics[width = \textwidth]{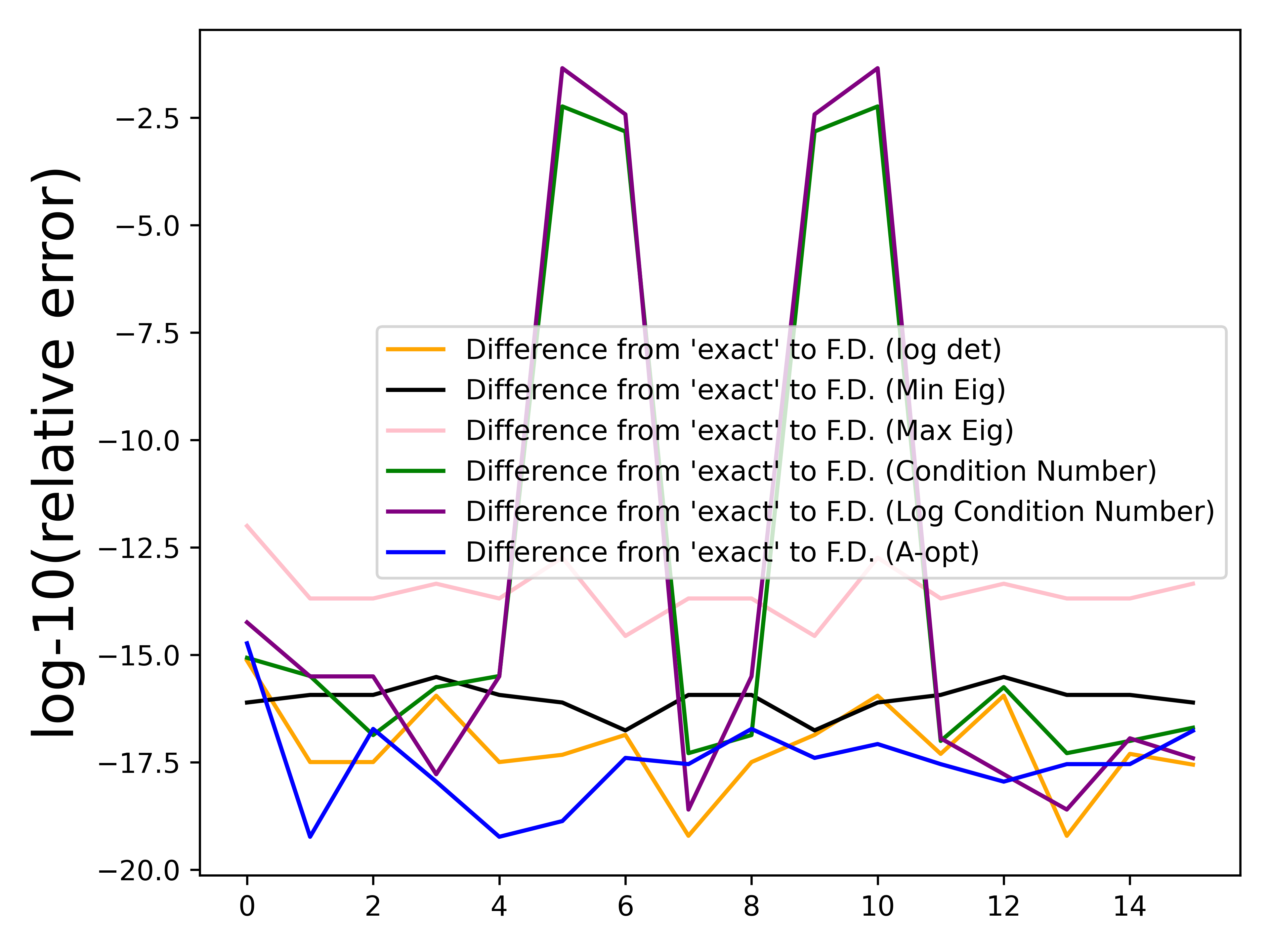}
\caption{Difference between finite difference and exact values for the second derivatives of each condition using the \texttt{numpy.linalg.cond} function instead of the maximum and minimum eigenvalues. Comparing with Figure~\ref{fig:si-second-deriv-numerics}, indices 5, 6, 9, and 10 have significant numerical deviation from the numerical second derivative.} \label{fig:si-second-deriv-numerics-bad}
\end{figure}

Figure~\ref{fig:si-second-deriv-numerics-bad} shows that some of the off-diagonal elements of the Hessian have larger numerical error than expected. This error exceeds the acceptable numerical threshold while using a step size of $10^{-6}$. This error results from using \texttt{numpy.linalg.cond} to calculate the condition number instead of the formula specified in Eq.~\ref{eq:obj-cond-SI} using the maximum and minimum eigenvalues from \texttt{numpy.linalg.eig}. The slight difference is that \texttt{numpy.linalg.cond} defines the condition number as follows: ``The condition number of $\boldsymbol{M}$ is defined as the norm of $\boldsymbol{M}$ times the norm of the inverse of $\boldsymbol{M}$.'' However, for symmetric real matrices, Eq.~\ref{eq:obj-cond-SI} is efficient as we need the entire eigenvalue-eigenvector solution to compute the derivatives in Eqs.~\ref{eq:ME-opt-first-deriv-si} and ~\ref{eq:ME-opt-second-deriv-si}. The difference between the condition number from the minimum and maximum eignevalues from \texttt{numpy.linalg.eig} and \texttt{numpy.linalg.cond} was significant enough at a step size of $10^{-6}$ to cause the Hessian to differ (Figure~\ref{fig:si-second-deriv-numerics-bad}). Therefore, we ultimately utilize \texttt{numpy.linalg.eig} to compute the condition number of the matrix using Eq.~\ref{eq:obj-cond-SI} which has acceptable numerical error, as shown in Figure~\ref{fig:si-second-deriv-numerics}.

We have included this short analysis to emphasize how important it is that the numerical method is consistent throughout all related computations, especially when utilizing contributed scientific computing tools. These tools are extremely useful, but when mathematics need to be precise and within solver tolerances at or below $10^{-5}$, a difference as small as using a different method to find eigenvalues and the condition number of a matrix can corrupt your exact derivatives enough to be unusable in practice.

\section{Symmetric representations of the FIM within the Grey Box formulation} \label{sec:si-symmetry}
One problem that could be faced when solving formulation~\ref{eq:formulation-gb} is that symmetry is not guaranteed at each iteration when modeling the entire FIM. For example, the following 2-by-2 matrix does not guarantee that $M_{12}$ and $M_{21}$ are the same at each iteration, rather that equality is enforced at a feasible solution.

\begin{align}
    &\boldsymbol{M} = 
        \begin{bmatrix}
            M_{11} & M_{12}  \\[0.35cm]
            M_{21} & M_{22}  \\[0.1cm]
        \end{bmatrix}
    &
\end{align}
The problem here is that evaluating these metrics, especially the eigenvalue metrics, rely on the symmetry of the matrix to be well-behaved. Therefore, we need only to model a triangular version of the matrix, and, within the Grey Box, enforce this symmetry directly:

\begin{align}
    \text{Full Matrix} & \qquad\text{Grey Box Matrix} \notag \\
    \notag\\
    \boldsymbol{M} = 
        \begin{bmatrix}
            M_{11} & M_{12}  \\[0.35cm]
            M_{21} & M_{22}  \\[0.1cm]
        \end{bmatrix}
    & \qquad\boldsymbol{M}^* = 
        \begin{bmatrix}
            M_{11} & {\color{blue}M_{12}}  \\[0.35cm]
            {\color{blue}M_{12}} & M_{22}  \\[0.1cm]
        \end{bmatrix}
\end{align}
We do not need to represent $M_{21}$ in the base formulation because we know, by definition, the matrix is symmetric and the Grey Box can fill in the rest of the matrix itself. Since this is the case, we need only send an upper (or lower) triangular version of the FIM to the Grey Box and then symmetry can be enforced. This reduces the number of constraints in the model and also ensures the $\boldsymbol{M}$ is symmetric at each iteration. However, one issue arises in that the inputs to the Grey Box are reduced, meaning the expected shape of both the Jacobian and the Hessian during optimization will be changed. This is clear when observing the shape of the vectorized Jacobian for the full versus the triangular inputs:

\begin{align}
    \text{Full Matrix Input} & \qquad\text{Upper Triangle Matrix Input} \notag \\
    \notag\\
    \boldsymbol{J} = 
        \begin{bmatrix}
            \frac{\partial\,\Psi}{\partial M_{11}} \\[0.2cm]
            \frac{\partial\,\Psi}{\partial M_{12}} \\[0.2cm]
            \frac{\partial\,\Psi}{\partial M_{21}} \\[0.2cm]
            \frac{\partial\,\Psi}{\partial M_{22}} \\[0.1cm]
        \end{bmatrix} \;\;\;\;
    & \qquad\qquad\ \boldsymbol{J}^* = 
        \begin{bmatrix}
            \frac{\partial\,\Psi}{\partial M_{11}}  \\[0.2cm]
            \frac{\partial\,\Psi}{\partial M_{12}} \\[0.2cm]
            \frac{\partial\,\Psi}{\partial M_{22}} \\[0.1cm]
        \end{bmatrix}
\end{align}
The formulas we developed are in the complete matrix space, not the symmetric matrix space, which is much more convenient for mathematical operations, but now requires additional care while constructing the Jacobian. Here, since we assume that $M_{12}$ and $M_{21}$ are the same, we also assume that the {\emph{changes}} to $M_{12}$ and $M_{21}$ are the same. Thus, we can get the correct Jacobian by including the $M_{21}$ derivative with the $M_{12}$ term, as follows:

\begin{align}
    \boldsymbol{J} = 
        \begin{bmatrix}
            \frac{\partial\,\Psi}{\partial M_{11}} \\[0.2cm]
            {\color{red}\frac{\partial\,\Psi}{\partial M_{12}}} \\[0.2cm]
            {\color{red}\frac{\partial\,\Psi}{\partial M_{21}}} \\[0.2cm]
            \frac{\partial\,\Psi}{\partial M_{22}} \\[0.1cm]
        \end{bmatrix}
    \rightarrow \boldsymbol{J}^* = 
        \begin{bmatrix}
            \frac{\partial\,\Psi}{\partial M_{11}} \\[0.2cm]
            {\color{red}\frac{\partial\,\Psi}{\partial M_{12}} + \frac{\partial\,\Psi}{\partial M_{21}}} \\[0.2cm]
            \frac{\partial\,\Psi}{\partial M_{22}} \\[0.1cm]
        \end{bmatrix}
\end{align}
Since we utilize the triangular matrix to complete the full matrix $\boldsymbol{M}$, we have access to the full-space derivatives and can easily find the augmented Jacobian and take the correct step. We highlight the relevant terms in red in the previous equation for emphasis.

With the Hessian ($\boldsymbol{H}$), or second derivative, this is slightly more complicated. With vectorized notation, the Hessian matrix for a 2-by-2 system can be represented as follows:

\begin{align}
    \boldsymbol{M} = &
        \begin{bmatrix}
            M_{11} & M_{12}  \\[0.2cm]
            M_{21} & M_{22}  \\[0.1cm]
        \end{bmatrix} \\
        \notag \\
    \boldsymbol{J} = &
        \begin{bmatrix}
            \frac{\partial\,\Psi}{\partial M_{11}} \\[0.2cm]
            \frac{\partial\,\Psi}{\partial M_{12}} \\[0.2cm]           \frac{\partial\,\Psi}{\partial M_{21}} \\[0.2cm]
            \frac{\partial\,\Psi}{\partial M_{22}} \\[0.1cm]
        \end{bmatrix}\\
        \notag \\
    \boldsymbol{H} = &
        \begin{bmatrix}
            \frac{\partial^2\,\Psi}{\partial M_{11}\partial M_{11}} & \frac{\partial^2\,\Psi}{\partial M_{11}\partial M_{12}} & \frac{\partial^2\,\Psi}{\partial M_{11}\partial M_{21}} &
            \frac{\partial^2\,\Psi}{\partial M_{11}\partial M_{22}}\\[0.2cm]
            \frac{\partial^2\,\Psi}{\partial M_{12}\partial M_{11}} & \frac{\partial^2\,\Psi}{\partial M_{12}\partial M_{12}} & \frac{\partial^2\,\Psi}{\partial M_{12}\partial M_{21}} &
            \frac{\partial^2\,\Psi}{\partial M_{12}\partial M_{22}}\\[0.2cm]
            \frac{\partial^2\,\Psi}{\partial M_{21}\partial M_{11}} & \frac{\partial^2\,\Psi}{\partial M_{21}\partial M_{12}} & \frac{\partial^2\,\Psi}{\partial M_{21}\partial M_{21}} &
            \frac{\partial^2\,\Psi}{\partial M_{21}\partial M_{22}}\\[0.2cm]
            \frac{\partial^2\,\Psi}{\partial M_{22}\partial M_{11}} & \frac{\partial^2\,\Psi}{\partial M_{22}\partial M_{12}} & \frac{\partial^2\,\Psi}{\partial M_{22}\partial M_{21}} &
            \frac{\partial^2\,\Psi}{\partial M_{22}\partial M_{22}}\\[0.1cm]
        \end{bmatrix}
\end{align}
Typically, the Hessian only requires a triangular representation to the solver as the Hessian is also a symmetric matrix. However, when considering the contraction to the symmetric space, the Hessian in the symmetric space misses some elements of the full Hessian (marked in red):

\begin{align}
    \boldsymbol{H} = &
        \begin{bmatrix}
            \frac{\partial^2\,\Psi}{\partial M_{11}\partial M_{11}} & \frac{\partial^2\,\Psi}{\partial M_{11}\partial M_{12}} & {\color{red}\frac{\partial^2\,\Psi}{\partial M_{11}\partial M_{21}}} &
            \frac{\partial^2\,\Psi}{\partial M_{11}\partial M_{22}}\\[0.2cm]
            \frac{\partial^2\,\Psi}{\partial M_{12}\partial M_{11}} & \frac{\partial^2\,\Psi}{\partial M_{12}\partial M_{12}} & {\color{red}\frac{\partial^2\,\Psi}{\partial M_{12}\partial M_{21}}} &
            \frac{\partial^2\,\Psi}{\partial M_{12}\partial M_{22}}\\[0.2cm]
            {\color{red}\frac{\partial^2\,\Psi}{\partial M_{21}\partial M_{11}}} & {\color{red}\frac{\partial^2\,\Psi}{\partial M_{21}\partial M_{12}}} & {\color{red}\frac{\partial^2\,\Psi}{\partial M_{21}\partial M_{21}}} &
            {\color{red}\frac{\partial^2\,\Psi}{\partial M_{21}\partial M_{22}}}\\[0.2cm]
            \frac{\partial^2\,\Psi}{\partial M_{22}\partial M_{11}} & \frac{\partial^2\,\Psi}{\partial M_{22}\partial M_{12}} & {\color{red}\frac{\partial^2\,\Psi}{\partial M_{22}\partial M_{21}}} &
            \frac{\partial^2\,\Psi}{\partial M_{22}\partial M_{22}}\\[0.1cm]
        \end{bmatrix}\\
        \notag \\
    \boldsymbol{H}^* = &
        \begin{bmatrix}
            \frac{\partial^2\,\Psi}{\partial M_{11}\partial M_{11}} & \frac{\partial^2\,\Psi}{\partial M_{11}\partial M_{12}} &
            \frac{\partial^2\,\Psi}{\partial M_{11}\partial M_{22}}\\[0.2cm]
            \frac{\partial^2\,\Psi}{\partial M_{12}\partial M_{11}} & \frac{\partial^2\,\Psi}{\partial M_{12}\partial M_{12}} &
            \frac{\partial^2\,\Psi}{\partial M_{12}\partial M_{22}}\\[0.2cm]
            \frac{\partial^2\,\Psi}{\partial M_{22}\partial M_{11}} & \frac{\partial^2\,\Psi}{\partial M_{22}\partial M_{12}} &
            \frac{\partial^2\,\Psi}{\partial M_{22}\partial M_{22}}\\[0.1cm]
        \end{bmatrix}
\end{align}
Using mapping, we can then achieve the following, correct representation of the Hessian that considers all terms that are missed.

\begin{align}
    \boldsymbol{H}^* = &
        \begin{bmatrix}
            \frac{\partial^2\,\Psi}{\partial M_{11}\partial M_{11}} & \frac{\partial^2\,\Psi}{\partial M_{11}\partial M_{12}} + {\color{red}\frac{\partial^2\,\Psi}{\partial M_{11}\partial M_{21}}} &
            \frac{\partial^2\,\Psi}{\partial M_{11}\partial M_{22}}\\[0.3cm]
            \frac{\partial^2\,\Psi}{\partial M_{12}\partial M_{11}} + {\color{red}\frac{\partial^2\,\Psi}{\partial M_{12}\partial M_{21}}} & \frac{\partial^2\,\Psi}{\partial M_{12}\partial M_{12}} + {\color{red}\frac{\partial^2\,\Psi}{\partial M_{12}\partial M_{21}}} &
            \frac{\partial^2\,\Psi}{\partial M_{12}\partial M_{22}} + {\color{red}\frac{\partial^2\,\Psi}{\partial M_{21}\partial M_{22}}}\\
            & + {\color{red}\frac{\partial^2\,\Psi}{\partial M_{21}\partial M_{12}}} + {\color{red}\frac{\partial^2\,\Psi}{\partial M_{21}\partial M_{21}}}& \\[0.3cm]
            \frac{\partial^2\,\Psi}{\partial M_{22}\partial M_{11}} & \frac{\partial^2\,\Psi}{\partial M_{22}\partial M_{12}} + {\color{red}\frac{\partial^2\,\Psi}{\partial M_{22}\partial M_{21}}} &
            \frac{\partial^2\,\Psi}{\partial M_{22}\partial M_{22}}\\[0.1cm]
        \end{bmatrix}
\end{align}
Some terms are obvious to map; for instance, the $\frac{\partial^2\,\Psi}{\partial M_{11}\partial M_{21}}$ element is simply mapped to the $\frac{\partial^2\,\Psi}{\partial M_{11}\partial M_{12}}$ following the same logic that the changes to $M_{12}$ and $M_{21}$ are the same. This is also easy to consider during the construction of a triangular Hessian, as the symmetry is held from $\boldsymbol{H}$ to $\boldsymbol{H}^*$. However, when mapping an element $\frac{\partial^2\,\Psi}{\partial M_{ij}\partial M_{ji}}$ with $j\neq i$ to the reduced counterpart, we now map onto the diagonal of the Hessian, meaning we must count both components, not just one (as in the reduced symmetry-holding case). 

Special care must be taken when utilizing full-space matrix representation on a symmetric-space object. Recently, it has been published that a long-standing derivative formula that has been around for over 60 years is incorrect, as this symmetric-to-full mapping is not taken into account \cite{srinivasan2023bad-symm-deriv}.

\section{TCLab: Eigenvalue and Eigenvector Analysis} \label{sec:si-tclab-eigenanalysis}
This section will describe a brief eigenvalue and eigenvector analysis of the uncertainty reduction for the TCLab system. We first take a look at the differences in the orders of magnitude of the covariance matrix with preliminary data only (Eq.~\ref{eq:si-tclab-cov-before}) and including the optimal experimental data (Eq.~\ref{eq:si-tclab-cov-after})

\begin{align}
&\;\;\boldsymbol{V}_{\boldsymbol{\theta},\text{before}} = 
    \begin{bmatrix}
        2.28 \cdot 10^{-5} & 1.65 \cdot 10^{1} & 2.37 & -1.89 \cdot 10^{1}   \\[0.35cm]
        1.65 \cdot 10^{1}  & 489 \cdot 10^{2} & 1.89 \cdot 10^{6} & -1.51 \cdot 10^{7}    \\[0.35cm]
        2.37 & 1.89 \cdot 10^{6} & 2.72 \cdot 10^{5} & -2.16 \cdot 10^{6}   \\[0.35cm]
        -1.89 \cdot 10^{1} & -1.51 \cdot 10^{7} & -2.16 \cdot 10^{6} & 1.72 \cdot 10^{7}   \\
    \end{bmatrix} \label{eq:si-tclab-cov-before} \\ \notag \\
&\;\;\boldsymbol{V}_{\boldsymbol{\theta},\text{after}} = 
    \begin{bmatrix}
        1.47 \cdot 10^{-6} & 4.63 \cdot 10^{-3} & 6.74 \cdot 10^{-4} & -5.32 \cdot 10^{-3}   \\[0.35cm]
        4.63 \cdot 10^{-3} & 1.5 \cdot 10^{3} & 2.16 \cdot 10^{2} & -1.72 \cdot 10^{3}   \\[0.35cm]
        6.74 \cdot 10^{-4} & 2.16 \cdot 10^{2} & 3.10 \cdot 10^{1} & -2.47 \cdot 10^{2}    \\[0.35cm]
        -5.32 \cdot 10^{-3} & -1.72 \cdot 10^{3} & -2.47 \cdot 10^{2} & 1.96 \cdot 10^{3}   \\
    \end{bmatrix} \label{eq:si-tclab-cov-after}
\end{align}
At first glance, these matrices appear as a wall of numbers, but generally we look for two things: (i) magnitude of the values of the covariance matrix and (ii) the eigendecomposition. From the perspective of magnitude, the covariance matrix before (Eq.~\ref{eq:si-tclab-cov-before}) has significantly larger entries, indicating that there is likely higher uncertainty with the parameters. Also, there are wildly different orders of magnitude, indicating numerical instability or poor conditioning of the matrix. However, how these uncertainties are realized must be analyzed visually (Figure~\ref{fig:tclab-pairwise-uncertainty-combined}) or using the eigenvalues, which are used to plot the pairwise uncertainties. The eigendecomposition consists of eigenvalues and eigenvectors corresponding to the solution of Eq.~\ref{eq:eigenvalue-problem-si}. For the before and after covariance matrices, we have the following:

\begin{align}
&\;\;\boldsymbol{\lambda}_{\boldsymbol{\theta},\text{before}} = 
    \begin{bmatrix}
        3.07 \cdot 10^{7} & 6.44 \cdot 10^{-4} & 1.90 \cdot 10^{-5} & 1.22 \cdot 10^{-6}   \\
    \end{bmatrix} \label{eq:si-tclab-eigval-before} \\ \notag \\
&\;\;\boldsymbol{\Lambda}_{\boldsymbol{\theta},\text{before}} = 
    \begin{bmatrix}
        0. & -0.028 & -0.137 & \boldsymbol{-0.990}   \\[0.35cm]
        \boldsymbol{0.655} & \boldsymbol{0.720} & -0.228 & 0.011   \\[0.35cm]
        0.094 & -0.377 & \boldsymbol{-0.911} & 0.137   \\[0.35cm]
        \boldsymbol{-0.749} & \boldsymbol{0.582} & -0.314 & 0.027   \\
    \end{bmatrix} \label{eq:si-tclab-eigvec-before} \\ \notag \\
&\;\;\boldsymbol{\lambda}_{\boldsymbol{\theta},\text{after}} = 
    \begin{bmatrix}
        3.49 \cdot 10^{3} & 5.67 \cdot 10^{-4} & 1.48 \cdot 10^{-5} & 8.31 \cdot 10^{-7}   \\
    \end{bmatrix} \label{eq:si-tclab-eigval-after} \\ \notag \\
&\;\;\boldsymbol{\Lambda}_{\boldsymbol{\theta},\text{after}} = 
    \begin{bmatrix}
        0. & 0.025 & 0.140 & \boldsymbol{0.990}   \\[0.35cm]
        \boldsymbol{0.655} & \boldsymbol{-0.723} & 0.218 & -0.013   \\[0.35cm]
        0.094 & 0.364 & \boldsymbol{0.916} & -0.139  \\[0.35cm]
        \boldsymbol{-0.749} & \boldsymbol{-0.587} & 0.306 & -0.029   \\
    \end{bmatrix} \label{eq:si-tclab-eigvec-after}
\end{align}
where the eigenvectors are column-wise with contributions according to the following directions:
\begin{align}
    &\;\;\boldsymbol{v}_\text{direction} = 
    \begin{bmatrix}
        U_a \\[0.35cm]
        U_b \\[0.35cm]
        \frac{1}{Cp_H} \\[0.35cm]
        \frac{1}{Cp_S}  \\
    \end{bmatrix} \label{eq:si-tclab-directions}
\end{align}
For emphasis, the largest contributor(s) to the direction of uncertainty are bolded in the eigenvector matrix. The interpretation is as follows: A larger value for an eigenvalue means a larger amount of uncertainty. The direction of that uncertainty is specified by the eigenvector corresponding to that eigenvalue. For example, the direction of the largest uncertainty before running an optimal experiment has a magnitude on the order of $10^7$ whereas the largest direction of uncertainty after conducting the experiment is predicted to be on the order of $10^3$. This is a large uncertainty direction, but we can see that the direction of uncertainty is almost identical, pointing slightly more in the direction of $\frac{1}{Cp_S}$ but also in the direction of $U_b$. This is directly represented in Figure~\ref{fig:tclab-pairwise-uncertainty-combined}b, by looking at the pairwise uncertainty for $\frac{1}{Cp_S}$ vs. $U_b$. Although difficult to make out, the direction of uncertainty is nearly identical but the black border is smaller than the gray border, indicating reduction in uncertainty. Similar results can be demonstrated for the other parameters. For instance, the smallest direction of uncertainty (highest direction of confidence) is the smallest eigenvalue, whose eigenvector in both cases points almost entirely in the direction of $U_a$. In Figure~\ref{fig:tclab-pairwise-uncertainty-combined}b, it is clear that $U_a$ is the only parameter that likely is estimable with physical bounds and realistic uncertainty both before and after the optimal experiment.

In general, both Figure~\ref{fig:tclab-pairwise-uncertainty-combined}b and an eigenanalysis give visual and quantitative insight, respectively, to the conditioning of a system. Although the experiment decreases uncertainty in model parameters, we can see through the eigenanalysis that there remains a consistent problem with the system. The reparameterization helps with identifying good experiments, but does not fix potential structural identifiability issues with the model. These analyses provide insight that the FIM and MBDoE are alone not a complete tool for model-building, and successful predictive model building may require different methods and will be addressed in future work.

\end{document}